\documentclass[titlepage,11pt]{article}
\usepackage{latexsym}
\usepackage{amsfonts}
\usepackage{amsmath}
\usepackage{tikz}
\usepackage{float}
\usetikzlibrary{arrows}
\usetikzlibrary{decorations.markings}
\newcommand\blackslug{\hbox{\hskip 1pt \vrule width 4pt height 8pt depth 1.5pt
        \hskip 1pt}}
\newcommand\bbox{\hfill \quad \blackslug \bigbreak}
\def\d{\hbox{-}}
\def\c{\hbox{-}\cdots\hbox{-}}
\def\ll{,\ldots,}
\newcommand{\vare}{\varepsilon}

\title{Concatenating bipartite graphs}
\author{Maria Chudnovsky\thanks{Supported by NSF grant DMS 1763817 and US Army Research Office Grant W911NF-16-1-0404.}\\
Princeton University, Princeton, NJ 08544, USA
\\
\\
Patrick Hompe\\
Princeton University, Princeton, NJ 08544, USA
\\
\\
Alex Scott\thanks{Supported by a Leverhulme Trust Research
Fellowship.}\\
Mathematical Institute, University of Oxford, Oxford OX2 6GG, UK
\\
\\
Paul Seymour\thanks{Supported by ONR grant N00014-14-1-0084 and NSF
grants DMS-1265563 and DMS-1800053.}\\
Princeton University, Princeton, NJ 08544, USA
\\
\\
Sophie Spirkl\thanks{This material is based upon work supported by the National Science
Foundation under Award No. DMS-1802201.}\\
Rutgers University, New Brunswick, NJ 08901, USA}

\date{September 21, 2018; revised \today}

\newtheorem{thm}{}[section]

\newcommand{\Proof}{\noindent{\bf Proof.}\ \ }

\begin{document}
\maketitle
\begin{abstract}
Let $x,y\in (0,1]$; and let $A,B,C$ be disjoint nonempty subsets of a graph $G$, where every vertex in $A$ has at least $x|B|$
neighbours in $B$, and every vertex in $B$ has at least $y|C|$ neighbours in $C$. We denote by $\phi(x,y)$
the maximum $z$ such that, in all such graphs $G$,
there is a vertex $v\in C$ that is joined to at least $z|A|$ vertices in $A$ by two-edge paths. 
The function $\phi$ is interesting, and we investigate some of its properties. For instance, we show that
\begin{itemize}
\item $\phi(x,y)=\phi(y,x)$ for all $x,y$; and 
\item for each integer $k>1$, there is a discontinuity in $\phi(x,x)$ when $x=1/k$: $\phi(x,x)\le 1/k$ when $x\le 1/k$, and
$\phi(x,x)\ge \frac{2k-1}{2k(k-1)}$ when $x>1/k$.
\end{itemize}
We raise several questions and
conjectures.

\end{abstract}

\section{Introduction}
All graphs in this paper are finite, and have no loops or multiple edges.
We denote the semi-open interval $\{x:0< x\le 1\}$ of real numbers by $(0,1]$.
Let $x,y\in (0,1]$; and let $A,B,C$ be disjoint nonempty subsets of a graph $G$, where every vertex in $A$ has at least $x|B|$
neighbours in $B$, and every vertex in $B$ has at least $y|C|$ neighbours in $C$. If we ask for a real number $z$
such that we can guarantee that some vertex in $A$ can reach at least $z|C|$ vertices in $C$ by two-edge paths,
then $z$ must be at most $y$, since perhaps all the vertices in $B$ have the same neighbours in $C$. But in the reverse
direction the question becomes much more interesting; that is, we ask for $z$ such that some vertex in $C$ can reach 
at least 
$z|A|$ vertices in $A$ by two-edge paths. Then there might well be values of $z>\max(x,y)$ with this property.

Let us say this more precisely. A {\em tripartition} of a graph $G$ is a partition $(A,B,C)$ of $V(G)$ where $A,B,C$ are 
all nonempty stable sets. For $x,y\in (0,1]$, we say a graph $G$ is {\em $(x,y)$-constrained, via  a tripartition $(A,B,C)$,}
if 
\begin{itemize}
\item every vertex in $A$ has at least $x|B|$
neighbours in $B$; 
\item every vertex in $B$ has at least $y|C|$ neighbours in $C$; and
\item there are no edges between $A$ and $C$.
\end{itemize}
For $v\in V(G)$, $N(v)$ denotes its set of neighbours, and $N^2(v)$ is the set of vertices with distance exactly two 
from $v$. We write $N^2_A(v)$ for $N^2(v)\cap A$, and so on.
A first observation:
\begin{thm}\label{suptomax}
Let $x,y\in (0,1]$, and let $Z$ be the set of all $z\in (0,1]$ such that,
for every graph $G$, if $G$ is $(x,y)$-constrained via $(A,B,C)$
then $|N^2_A(v)|\ge z|A|$ for some $v\in C$. Then $\sup\{z\in Z\}$ belongs to $Z$.
\end{thm}
\Proof
Let $z'=\sup\{z\in Z\}$, and let $G$ be an $(x,y)$-constrained graph, via $(A,B,C)$. We must show that
$|N^2_A(v)|\ge z'|A|$ for some $v\in C$.
We may assume that $z'>0$; so there exists
$z$ with $0< z<z'$, such that $ \lceil z|A| \rceil=\lceil z'|A| \rceil$. Since $z'=\sup\{z\in Z\}$
and $z<z'$, and $Z$ is an initial interval of $(0,1]$, it follows that $z\in Z$,
and so $|N^2_A(v)|\ge z|A|$ for some $v\in C$. Consequently $|N^2_A(v)|\ge \lceil z|A|\rceil\ge z'|A|$,
as required. This proves \ref{suptomax}.~\bbox

We define $\phi(x,y)$ to be $\sup\{z\in Z\}$, as defined in \ref{suptomax}. The objective of this paper
is to study the properties of the function $\phi$.
We have a trivial lower bound:
\begin{thm}\label{maxbound}
$\phi(x,y)\ge \max(x,y)$ for all $x,y>0$.
\end{thm}
\Proof
Let $G$ be $(x,y)$-constrained, via $(A,B,C)$. Since every vertex in $A$ has at least $x|B|$ neighbours in $B$, and
$B\ne \emptyset$, there exists $u\in B$ with at least $x|A|$ neighbours in $A$; let $v\in C$ be adjacent to $u$
(this is possible since $y>0$), and then $|N^2_A(v)|\ge x|A|$. Consequently $\phi(x,y)\ge x$. Now every
vertex in $A$ can reach at least $y|C|$ vertices in $C$ by two-edge paths (since $x>0$); and so by averaging, some
vertex in $C$ can reach at least $y|A|$ vertices in $A$ by two-edge paths. Hence $\phi(x,y)\ge y$. 
This proves \ref{maxbound}.~\bbox

And a trivial upper bound:
\begin{thm}\label{cayleybound}
For all $x,y\in (0,1]$, 
$$\phi(x,y)\le \frac{\lceil kx\rceil +\lceil ky\rceil -1}{k}$$ 
for every integer $k\ge 1$.
\end{thm}
\Proof
Let $x,y\in (0,1]$, and let $k\ge 1$ be an integer. Let $A,B,C$ be three disjoint sets each of cardinality $k$, 
where $A=\{a_1\ll a_k\}$, $B=\{b_1\ll b_k\}$ and $C=\{c_1\ll c_k\}$. Make a graph $G$ with vertex set
$A\cup B\cup C$ as follows. Let $g=\lceil kx\rceil$, and for $1\le i\le k$ make $a_i$ adjacent to 
$b_i, b_{i+1}\ll b_{i+g-1}$ (reading subscripts modulo $k$). Now let $h=\lceil ky\rceil$, and for 
$1\le i\le k$ make $b_i$ adjacent to 
$c_i, c_{i+1}\ll c_{i+h-1}$ (reading subscripts modulo $k$). Then $G$ is $(x,y)$-constrained via $(A,B,C)$;
and for $1\le i\le k$, $N^2_A(c_i)=\{a_i, a_{i-1}\ll a_{i-g-h+2}\}$ (again, reading subscripts modulo $k$).
Consequently $\phi(x,y)\le (g+h-1)/k$. This proves \ref{cayleybound}.~\bbox

In particular, we have:
\begin{thm}\label{intk}
For every integer $k\ge 1$, if $x,y>0$ and $\max(x,y) = 1/k$ then $\phi(x,y) = 1/k$.
\end{thm}
\Proof
From \ref{maxbound}, $\phi(x,y)\ge 1/k$; and the graph consisting of $k$ disjoint three-vertex paths shows that $\phi(x,y)\le 1/k$.
(This also follows from \ref{cayleybound}, since $\lceil kx\rceil, \lceil ky\rceil =1$.)
This proves \ref{intk}.~\bbox

What makes the function $\phi$ interesting is that for some values of $x,y$, \ref{maxbound} is far from best possible,
and indeed \ref{cayleybound} seems closer to the truth. We were originally motivated by the hope of extending
Kneser's theorem from additive group theory~\cite{kneser} to a general graph-theoretic setting, and 
a corresponding wild conjecture that
the bound in \ref{cayleybound} is always best possible, that is, that for all $x,y\in (0,1]$, there
is an integer $k>0$ with $\phi(x,y)= \frac{\lceil kx\rceil +\lceil ky\rceil -1}{k}$. 
This turns out to be false, but perhaps not ridiculously false; maybe something
like it is true.

There are two other related problems:
\begin{itemize}
\item 
Let us say $G$ is {\em $(x,y)$-biconstrained} ({\em via} $(A,B,C)$) if $G$ is $(x,y)$-constrained via $(A,B,C)$, and in addition
\begin{itemize}
\item every vertex in $B$ has at least $x|A|$ neighbours in $A$, and
\item every vertex in $C$ has at least $y|B|$ neighbours in $B$.
\end{itemize}
\item Say $G$ is {\em $(x,y)$-exact} ({\em via} $(A,B,C)$) if
$G$ is $(x,y)$-constrained via $(A,B,C)$, and in addition there exist $x'\ge x$ and $y'\ge y$ such that
\begin{itemize}
\item every vertex in $A$ has exactly $x'|B|$ neighbours in $B$;
\item every vertex in $B$ has exactly $x'|A|$ neighbours in $A$;
\item every vertex in $B$ has exactly $y'|C|$ neighbours in $C$; and
\item every vertex in $C$ has exactly $y'|B|$ neighbours in $B$.
\end{itemize}
\end{itemize}
We shall sometime use ``mono-constrained'' to clarify that we mean the $(x,y)$-constrained case and not the $(x,y)$-biconstrained
case.
Let $\psi(x,y)$ be the analogue of $\phi(x,y)$ for biconstrained graphs; that is, the maximum $z$ such that for
all $G$, if $G$ is $(x,y)$-biconstrained via $(A,B,C)$, then $|N^2_A(v)|\ge z|A|$ for some $v\in C$. (As before, this
maximum exists.) Similarly, let $\xi(x,y)$ be the analogue of $\phi$ and $\psi$ for the exact case.
Then we have
\begin{thm}\label{trivialbounds}
For all $x,y\in (0,1]$, 
$$\max(x,y)\le \phi(x,y)\le \psi(x,y)\le \xi(x,y)\le \frac{\lceil kx\rceil +\lceil ky\rceil -1}{k}$$
for every integer $k\ge 1$.
\end{thm}
The proof of the non-trivial part of this is the same as the proof of \ref{cayleybound}.
One might hope that $\psi$ (and even more $\xi$) are better-behaved than $\phi$.

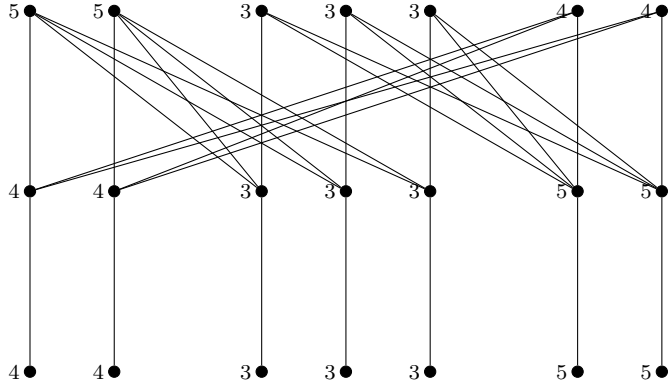
\begin{figure}[ht]
\centering

\begin{tikzpicture}[xscale=0.7,yscale=0.8,auto=left]
\tikzstyle{every node}=[inner sep=1.5pt, fill=black,circle,draw]

\def\s{.4}
\node (a1) at (6+\s,0) {};
\node (a2) at (8,0) {};
\node (a3) at (10-\s,0) {};
\node (a4) at (2,0) {};
\node (a5) at (4-\s,0) {};
\node (a6) at (12+\s,0) {};
\node (a7) at (14,0) {};
\def\r{-3}
\node (b1) at (6+\s,\r) {};
\node (b2) at (8,\r) {};
\node (b3) at (10-\s,\r) {};
\node (b4) at (12+\s,\r) {};
\node (b5) at (14,\r) {};
\node (b6) at (2,\r) {};
\node (b7) at (4-\s,\r) {};
\def\r{-6}
\node (c1) at (6+\s,\r) {};
\node (c2) at (8,\r) {};
\node (c3) at (10-\s,\r) {};
\node (c4) at (12+\s,\r) {};
\node (c5) at (14,\r) {};
\node (c6) at (2,\r) {};
\node (c7) at (4-\s,\r) {};

\foreach \from/\to in {a1/b1,a2/b2,a3/b3,a4/b1,a4/b2,a4/b3,a5/b1,a5/b2,a5/b3,a1/b4,a2/b4,a3/b4,a1/b5,a2/b5,a3/b5,a4/b6,a5/b7,a6/b4,a7/b5,a6/b6,a6/b7,a7/b6,a7/b7,b1/c1,b2/c2,b3/c3,b4/c4,b5/c5,b6/c6,b7/c7}
\draw [-] (\from) -- (\to);

\tikzstyle{every node}=[left]
\draw (a1) node []           {\scriptsize$3$};
\draw (a2) node []           {\scriptsize$3$};
\draw (a3) node []           {\scriptsize$3$};
\draw (a4) node []           {\scriptsize$5$};
\draw (a5) node []           {\scriptsize$5$};
\draw (a6) node []           {\scriptsize$4$};
\draw (a7) node []           {\scriptsize$4$};

\draw (b1) node []           {\scriptsize$3$};
\draw (b2) node []           {\scriptsize$3$};
\draw (b3) node []           {\scriptsize$3$};
\draw (b4) node []           {\scriptsize$5$};
\draw (b5) node []           {\scriptsize$5$};
\draw (b6) node []           {\scriptsize$4$};
\draw (b7) node []           {\scriptsize$4$};

\draw (c1) node []           {\scriptsize$3$};
\draw (c2) node []           {\scriptsize$3$};
\draw (c3) node []           {\scriptsize$3$};
\draw (c4) node []           {\scriptsize$5$};
\draw (c5) node []           {\scriptsize$5$};
\draw (c6) node []           {\scriptsize$4$};
\draw (c7) node []           {\scriptsize$4$};

\end{tikzpicture}

\caption{$\psi(13/27,1/9)=13/27$} \label{fig:exactcount}
\end{figure}

Let us see an example. Start with the graph of figure~\ref{fig:exactcount}. Each vertex has a number
written next to it in the figure; replace each vertex $v$ by a set $X_v$ of new vertices
of the specified cardinality, and for each edge $uv$ of the figure make every vertex in $X_u$ adjacent
to every vertex in $X_v$. This results in a graph with $81$ vertices, divided into three sets of $27$ corresponding to
the three rows of the figure; call these $A,B,C$. The graph produced is $(13/27,1/9)$-biconstrained via $(A,B,C)$, and
yet $|N^2_A(v)|=13$ for every vertex $v\in C$; so this proves that $\psi(13/27,1/9)\le 13/27$ (and therefore
equality holds, by \ref{maxbound}). This shows that there need not exist an integer $k$ with
$\psi(x,y)=\frac{\lceil kx\rceil +\lceil ky\rceil -1}{k}$. The same graph, used from bottom to top, shows that
$\psi(1/9,13/27)= 13/27$.

The example is not yet $(13/27,1/9)$-exact, because some vertices in $A$ have three, four or five neighbours in
$B$, and vice versa. We can make it exact as follows. For each edge $uv$ of the figure with $u$ in the second row
and $v$ in the third, the two sets $X_u,X_v$ have the same cardinality, one of three, four, five. Delete some edges
between $X_u$ and $X_v$ such that every vertex in $X_u$ has exactly three neighbours in $X_v$ and vice versa. Then
the modified graph is $(13/27,1/9)$-exact, and shows that $\xi(13/27,1/9)= 13/27$. Consequently, even for the
supposedly nicest function $\xi$ of our three functions, there is not always an integer $k$ with
$\xi(x,y)=\frac{\lceil kx\rceil +\lceil ky\rceil -1}{k}$.

So what can we prove about the functions $\phi$ and $\psi$? For which $x,y,z$ is $\phi(x,y)\ge z$, or $\psi(x,y)\ge z$?
In order to make the question a little more manageable,
we focus on seven special cases, $x=y$, and $z=1/2, 2/3, 1/3, 3/4, 2/5, 3/5$,  but in each case the results for $\phi$ and for $\psi$ are
quite different. The paper is organized as follows:
\begin{itemize}
\item We begin with a proof that $\phi(x,y)=\phi(y,x)$ for all $x,y$.
\item Then we give some general upper bounds on $\phi(x,y)$ and $\psi(x,y)$, particularly focussing on the
case when $x=y$. 
\item  Next we consider when $\phi(x,y)\ge 1/2$, or $\psi(x,y)\ge 1/2$. There are several theorems that this is true
for certain pairs $(x,y)$, and their union fills a good part of the $(x,y)$-square. We also give a number of constructions
that shows the statement is {\em not} true for certain pairs $(x,y)$. Ideally this would fill the complementary part
of the square, but there is an ``undecided'' band of varying width down the middle.
\item Then we do the same for $2/3$ instead of $1/2$; and then for $1/3, 3/4, 2/5,3/5$. 
\item Finally, we discuss some other questions and approaches.
\end{itemize}
Some of our results appear in \cite{hompe} and \cite{hompearxiv}.

\section{Weighted graphs and some linear programming}

In this section we prove that $\phi(x,y) = \phi(y,x)$ for all $x,y$. The argument uses linear programming, and we need 
some preparation. 
We denote the set of real numbers by $\mathbb{R}$, and the non-negative reals numbers by $\mathbb{R}_+$.
A {\em weighted graph} $(G,w)$ consists of a graph $G$ together with a function $w:V(G)\rightarrow\mathbb{R}_+$.
If $X\subseteq V(G)$, we denote $\sum_{v\in X}w(v)$ by $w(X)$.
Let $(G,w)$ be a weighted graph, and $(A,B,C)$ a tripartition of $G$. If $x,y\in (0,1]$, a weighted graph $(G,w)$ 
is {\em $(x,y)$-constrained} via $(A,B,C)$, if:
\begin{itemize}
\item $\sum_{v\in A}w(v)=\sum_{v\in B}w(v)=\sum_{v\in C}w(v) = 1$;
\item for each $v\in A$, $w(N(v)\cap B)\ge x$; and
\item for each $v\in B$, $w(N(v)\cap C)\ge y$.
\end{itemize}
Similarly, we say $(G,w)$ is {\em $(x,y)$-biconstrained} via $(A,B,C)$, if in addition:
\begin{itemize}
\item for each $v\in B$, $w(N(v)\cap A)\ge x$; and
\item for each $v\in C$, $w(N(v)\cap B)\ge y$.
\end{itemize}
To make the graph of figure \ref{fig:exactcount} into an appropriate weighted graph, divide
all the numbers by 27.
\begin{thm}\label{weighttonon}
For $x,y,z\in (0,1]$, the following are equivalent:
\begin{itemize}\item $\phi(x,y)\ge z$;
\item $w(N^2_A(v))\ge z$ for some $v\in C$, for every weighted graph $(G,w)$
that is $(x,y)$-constrained via a tripartition $(A,B,C)$.
\end{itemize}
Similarly, the following are equivalent:
\begin{itemize}
\item $\psi(x,y)\ge z$;
\item $w(N^2_A(v))\ge z$ for some $v\in C$, for every weighted graph $(G,w)$
that is $(x,y)$-biconstrained via a tripartition $(A,B,C)$.
\end{itemize}
\end{thm}
\Proof To prove the ``if'' direction of the first statement, let $G$ be $(x,y)$-constrained via $(A,B,C)$. Define
$w(v)=1/|A|$, for each $v\in A$, and $w(v)=1/|B|$ for $v\in B$ and similarly for $v\in C$. Then $(G,w)$
is an $(x,y)$-constrained weighted graph, and the claim follows. The ``if'' direction of the second statement is proved similarly.

For the ``only if'' direction, let $(G,w)$ be a weighted graph, $(x,y)$-constrained via $(A,B,C)$, and suppose such a weighted graph
can be chosen with $w(N^2_A(v))< z$ for each $v\in C$. Consequently we may choose $(G,w)$ such that in addition,
$w$ is rational-valued. Choose an integer $N>0$ such that $Nw(v)$ is an integer for each $v\in G$. 
For each $v\in V(G)$, take a set $X_v$ of $Nw(v)$ new vertices; and make a graph $G'$ with vertex set
$\bigcup_{v\in V(G)}X_v$, by making every vertex of $X_u$ adjacent to every vertex of $X_v$ for all adjacent $u,v\in V(G)$.
Let $A'=\bigcup_{v\in A}X_v$, and define $B', C'$ similarly; then $(A', B', C')$ is a tripartition of $G'$, and
$G'$ is $(x,y)$-constrained via $(A', B', C')$. Since in $G$, $w(N^2_A(v))< z$ for each $v\in C$, it follows that in $G'$,
$|N^2_{A'}(v')|< z|A'|$ for each $v'\in C'$, a contradiction. The ``only if'' direction of the second statement is similar.
This proves \ref{weighttonon}.~\bbox

Let $G$ be a graph with a bipartition $(A,B)$,
and let $w:B\rightarrow\mathbb{R}_+$ be some function. We define
$w(A\rightarrow B)$ to mean the minimum, over all $u\in A$, of $w(N(u))$ (taking
$w(A\rightarrow B)=0$ if $A=\emptyset$).

\begin{thm}\label{lpbip}
Let $G$ be a graph with a bipartition $(A,B)$,
and let $w:B\rightarrow\mathbb{R}_+$ be some function
such that $w(B)=1$.
Then either
\begin{itemize}
\item there is a function $w':B\rightarrow\mathbb{R}_+$, such that $w'(B)=1$ and
$w'(A\rightarrow B)\ge  w(A\rightarrow B)$, and such that $w'(v)=0$ for some $v\in B$; or
\item there is a function $f:A\rightarrow\mathbb{R}_+$, such that $f(A)=1$ and
$f(B\rightarrow A)\ge w(A\rightarrow B)$.
\end{itemize}
\end{thm}
\Proof We may assume that $A\ne \emptyset$. If some vertex in $A$ has no neighbour in $B$, then 
$w(A\rightarrow B)=0$ and the second bullet holds; so we assume that each vertex in $A$ has a neighbour in $B$.

Let $x=w(A\rightarrow B)$. The function $w'$, defined by $w'(v)=1/|B|$ for each $v\in B$, satisfies
$w'(A\rightarrow B)>0$, since
every vertex in $A$ has a neighbour in $B$. Thus we may assume that $x>0$, replacing $w$ by $w'$ if necessary.

Let $M$ be the $0/1$-matrix $(a_{uv}:u\in A, v\in B)$, where $a_{uv}=1$ if and only if $u,v$ are adjacent.
Let ${\bf 1}_A\in \mathbb{R}^A$ be the vector of all $1$'s, and define
${\bf 1}_B$ similarly. Then $w\in \mathbb{R}_+^B$ satisfies:
\begin{itemize}
\item ${\bf 1}_B^Tw = 1$; and
\item $Mw\ge x{\bf 1}_A$.
\end{itemize}
Consequently $b=w/x$ satisfies $b\in \mathbb{R}_+^B$, and
\begin{itemize}
\item ${\bf 1}_B^Tb = 1/x$; and
\item $Mb\ge {\bf 1}_A$.
\end{itemize}
Choose $q\in \mathbb{R}_+^B$ with $Mq\ge {\bf 1}_A$, with ${\bf 1}_B^Tq$ minimum. (This is possible by compactness.)
Thus ${\bf 1}_B^Tq \le 1/x$.
Since $Mq\ge {\bf 1}_A$ and $G$ has an edge, it follows that ${\bf 1}_B^Tq>0$; let
$1/y={\bf 1}_B^Tq$, and define
$w'=yq$. Then  $y\ge x$, and ${\bf 1}_B^Tw'=1$ and $Mw'\ge y{\bf 1}_A$, and so we may assume that $w'(v)>0$ for each $v\in B$,
because otherwise the first bullet holds.

Now $q$ minimizes ${\bf 1}_B^Tq$ subject to the linear programme
$q\in \mathbb{R}_+^B$ and $Mq\ge {\bf 1}_A$. From the
linear programming duality theorem,
there exists $p\in \mathbb{R}_+^A$ such that $p^TM\le {\bf 1}_B^T$, and $p^T{\bf 1}_A={\bf 1}_B^Tq=1/y$.
Define $f=yp$. Then $f:A\rightarrow\mathbb{R}_+$ satisfies $f(A)=1$, and
$f(N(v))\le y$ for each $v\in B$.

Let $v'\in B$;
we claim that $f(N(v'))=y$. This follows from the ``complementary slackness'' principle, but we give the argument
in full, as follows. Let $s=w'(v')(y-f(N(v')))$. Thus $s\ge 0$, and we will show $s=0$. We have
$$y=\sum_{v\in B}yw'(v)\ge s+\sum_{v\in B}\sum_{u\in N(v)}w'(v)f(u)=s+\sum_{u\in A}\sum_{v\in N(u)}f(u)w'(v)\ge s+\sum_{u\in A}yf(u)=s+y.$$
Consequently $s=0$, as claimed. Hence $f$ satisfies the second bullet. This proves \ref{lpbip}.~\bbox

From \ref{lpbip} we deduce a very useful result.

\begin{thm}\label{permute}
If $x,y\in (0,1]$ then $\phi(x,y)=\phi(y,x)$.
\end{thm}
\Proof Let $z=\phi(x,y)$, and choose a weighted graph $(G,w)$ that is $(x,y)$-constrained via $(A,B,C)$, such that
$w(N^2_A(v))\le z$ for each $z\in C$. Moreover, choose $G$ with $|V(G)|$ minimum. If there is a function 
$w':B\rightarrow\mathbb{R}_+$, such that $w'(B)=1$ and
$w'(A\rightarrow B)\ge  w(A\rightarrow B)$, and such that $w'(v)=0$ for some $v\in B$, then we may replace
$w$ by a new weight function, changing $w$ to $w'$ on $B$ and otherwise keeping $w$ unchanged, and then we may delete the vertex
$v\in B$ with $w'(v)=0$, contrary to the minimality of $|V(G)|$. Thus there is no such $w'$, and so by \ref{lpbip},
there is a function $f:A\rightarrow\mathbb{R}_+$, such that $f(A)=1$ and
$f(B\rightarrow A)\ge w(A\rightarrow B)\ge x$. Similarly, there is a function
$g:B\rightarrow\mathbb{R}_+$, such that $g(B)=1$ and
$g(C\rightarrow B)\ge y$. Let $H$ be the graph with bipartition $(A,C)$ in which $u\in A$ and $v\in C$ are adjacent
if $u\notin N^2_A(v)$ in $G$. Thus, in $H$, $w(C\rightarrow A)\ge 1-z$; and so from \ref{lpbip} and the minimality of $|V(G)|$,
there is a function $h:C\rightarrow\mathbb{R}_+$, such that $h(C)=1$ and
(in $H$) $h(A\rightarrow C)\ge 1-z$. Let $w'$ be defined by the union of $f,g$ and $h$ in the natural sense;
then $(G,w')$ is a weighted graph and is $(y,x)$-constrained via $(C,B,A)$, and $w'(N^2_C(v))\le z$ for each $v\in A$.
This proves that $\phi(y,z)\le z$, and so proves \ref{permute}.~\bbox

We remark that we have not been able to prove an analogue of \ref{permute} for the biconstrained case, or for the exact case,
although
we have no counterexample for either one.

There is another useful application of \ref{lpbip}, the following:

\begin{thm}\label{nodom}
Let $(G,w)$ be an $(x,y)$-constrained weighted graph, via $(A,B,C)$, such that $w(N^2_A(v))\le z$ for each $v\in C$.
Suppose that there exists $X\subseteq A$ with $|X|<z^{-1}$ such that $\bigcup_{v\in X}N^2_C(v) = C$. Then there exists $u\in A$
and a weighted graph $(G', w')$ such that
\begin{itemize}
\item $G'$ is obtained from $G$ by deleting $u$;
\item $(G', w')$ is $(x,y)$-constrained via $(A',B,C)$, where $A'=A\setminus \{u\}$;
\item in $G'$, $w'(N^2_{A'}(v))\le z$ for all $v\in C$; and
\item $w'(u)=w(u)$ for all $u\in B\cup C$.
\end{itemize}
\end{thm}
\Proof
Suppose not. Let $H$ be the graph with bipartition $(A,C)$, in which $u\in A$ and $v\in C$ are adjacent
if $u\notin N^2_A(v)$ in $G$. Then by \ref{lpbip}, applied to $H$, there is a  function
$h:C\rightarrow\mathbb{R}_+$, such that $h(C)=1$ and
(in $H$) $h(A\rightarrow C)\ge 1-z$. Consequently, in $G$, $h(N^2_C(v))\le z$ for each $v\in A$. In particular,
$h(N^2_C(v))\le z$ for each $v\in X$, and so $h(C)\le z|A|<1$, a contradiction. This proves \ref{nodom}.~\bbox

Let $x,y,z\in (0,1]$. We say that $(x,y,z)$ is {\em triangular} if no triangle-free graph $G$ admits a tripartition
$A,B,C$
of $V(G)$ with the following properties:
\begin{itemize}
\item $A,B,C$ are nonempty stable sets;
\item every vertex in $A$ has at least $x|B|$ neighbours in $B$;
\item every vertex in $B$ has at least $y|C|$ neighbours in $C$; and
\item every vertex in $C$ has at least $z|A|$ neighbours in $A$.
\end{itemize}

It is possible to reformulate results about $\phi(x,y)$ in terms of triangular triples, because we have:
\begin{thm}\label{rotate}
For $x,y,z\in (0,1]$, $\phi(x,y)> 1-z$ if and only if $(x,y,z)$ is triangular.
Consequently the three statements $\phi(x,y)\le  1-z$, $\phi(z,x)\le 1-y$, and $\phi(y,z)\le 1-x$ are equivalent.
\end{thm}
\Proof
Suppose that $(x,y,z)$ is not triangular. Then there is a triangle-free graph $G$ with a tripartition $(A,B,C)$,
satisfying the three bullets in the definition of ``triangular''. Let $H$ be the subgraph of $G$
with $V(H)=V(G)$, obtained by deleting all edges between $A$ and $C$. If $v\in C$, then $N^2_A(v)$ (defined
with respect to $H$) contains only vertices in $A$ that are nonadjacent to $v$ in $G$, since $G$ is triangle-free;
and so $|N^2_A(v)|\le |A|-z|A|$, since in $G$, $v$ has at least $z|A|$ neighbours in $A$. Consequently $\phi(x,y)\le 1-z$.

For the reverse implication, suppose that $\phi(x,y)\le 1-z$, and let $H$ be $(x,y)$-constrained via $(A,B,C)$,
such that $|N^2_A(v)|\le |A|-z|A|$ for each $v\in C$. Make a graph $G$ by adding certain edges to $H$, namely
for each $v\in C$ and $u\in A$, add an edge $uv$ if $u\notin N^2_A(v)$. Then $G$ is triangle-free, and every vertex
$v\in C$ is adjacent in $G$ to at least $|A|-(1-z)|A|=z|A|$ vertices in $A$; and so $(x,y,z)$ is not triangular.

In particular, $(x,y,z)$ is triangular if and only if $(z,x,y)$ is triangular; so it follows that
$\phi(x,y)\le  1-z$ if and only if $\phi(z,x)\le 1-y$, and similarly if and only if $\phi(y,z)\le 1-x$.
This proves \ref{rotate}.~\bbox

We call the equivalence of the second statement of \ref{rotate} ``rotating''.

\section{Constructions}
In this section we construct some graphs to prove upper bounds on $\phi(x,y)$ or $\psi(x,y)$ for certain values of $x,y$.
We begin with:
\begin{thm}\label{phiaddpath}
Let $x,y\in (0,1]$, and let $z\in (0,1]$ such that $z/(1-z)=\phi(x/(1-x),y/(1-y))$; then $\phi(x,y)\le z$.
\end{thm}
\Proof
Let $(G',w')$ be a weighted graph that is $(x/(1-x),y/(1-y))$-constrained via some tripartition $(A',B',C')$, such that $w'(N^2_{A'}(v))\le z/(1-z)$
for each $v\in C'$. Add three new vertices $a,b,c$ to $G'$, and two edges $ab$ and $bc$, forming $G$. Define $w$ by
\begin{eqnarray*}
w(a)&=&z\\
w(v)&=& (1-z)w'(v) \text{ for each }v\in A'\\
w(b)&=&x\\
w(v)&=& (1-x)w'(v) \text{ for each }v\in B'\\
w(c)&=& y\\
w(v)&=& (1-y)w'(v) \text{ for each }v\in C'.\\
\end{eqnarray*}
Then $G$ is $(x,y)$-constrained via $(A'\cup\{a\}, B'\cup \{b\}, C'\cup \{c\})$ and shows that $\phi(x,y)\le z$. 
This proves \ref{phiaddpath}.~\bbox

\begin{thm}\label{phi12curve}
Let $k\ge 0$ be an integer, and let $x,y\in (0,1]$ with $\frac{x}{1-kx}+\frac{y}{1-ky}\le 1$, with strict inequality if $x$ or $y$
is irrational; then $\phi(x,y)<\frac{1}{k+1}$.
\end{thm}
\Proof By increasing $x$ and $y$ if necessary, we may assume that $x,y$ are rational. 
Suppose first that $k=0$; then we may assume that $x+y=1$. Choose an integer $k\ge 1$
such that $kx$ (and hence $ky$) is an integer. By \ref{cayleybound}, 
$$\phi(x,y)\le \frac{\lceil kx\rceil +\lceil ky\rceil -1}{k}=x+y-1/k<1.$$
This completes the proof for $k=0$. For general $k$
we proceed by induction on $k$. We may assume that $k>0$;  let $x,y\in (0,1]$ with $\frac{x}{1-kx}+\frac{y}{1-ky}\le 1$,
with strict inequality if $x$ or $y$
is irrational. Let $x'=x/(1-x)$, and $y'=y/(1-y)$. Thus $x',y'\in (0,1]$ with 
$$\frac{x'}{1-(k-1)x'}+\frac{y'}{1-(k-1)y'}= \frac{x}{1-kx}+\frac{y}{1-ky}\le 1,$$
with strict inequality if $x'$ or $y'$
is irrational. From the inductive hypothesis, $\phi(x',y')<1/k$. Let $z$ satisfy
$z/(1-z)=\phi(x',y')$; then $z/(1-z)<1/k$, and so $z<1/(k+1)$. From \ref{phiaddpath}, $\phi(x,y)\le z<1/(k+1)$. 
This proves \ref{phi12curve}.~\bbox

\begin{thm}\label{psi12curve}
Let $k\ge 0$ be an integer, and let $x,y\in (0,1]$ with $x+ky\le 1$ and $kx+y\le 1$, with strict inequality in both if $x$ or $y$
is irrational; then $\psi(x,y)<\frac{1}{k}$.
\end{thm}
\Proof Again, we may assume that $x,y$ are rational. Let $s=\max(x,y)$; thus, $s<1/k$. Choose an integer $N\ge 1$
such that $p=xN/(1-(k-1)s)$ and $q=yN/(1-(k-1)s)$ are integers. (It follows that $p+q\le N$, from the hypothesis.)
Let $G$ be a graph with vertex set partitioned into three sets $A,B,C$,
with $|A|=N+k$ and $|B|=|C|=N+k-1$; let 
\begin{eqnarray*}
A &=& \{a_1, a_2\ll a_N, a_1'\ll a_{k-1}', a^*\},\\
B&=&\{b_1, b_2\ll b_N, b_1'\ll b_{k-1}'\},\\
C&=&\{c_1, c_2\ll c_N, c_1'\ll c_{k-1}'\}.
\end{eqnarray*}
Let $G$ have the following edges:
\begin{itemize}
\item for $1\le i\le N$, $a_i$ is adjacent to $b_i, b_{i+1}\ll b_{i+p-1}$ reading subscripts modulo $N$;
\item for $1\le i\le N$, $b_i$ is adjacent to $c_i, c_{i+1}\ll c_{i+q-1}$ reading subscripts modulo $N$;
\item for $1\le i \le k-1$, $a_i'$ is adjacent to $b_i'$, and $b_i'$ is adjacent to $c_i'$.
\item $a^*$ is adjacent to $b_i$ for $1\le i\le N$.
\end{itemize}
(Thus, this is the same as in the proof of \ref{phi12curve}, except for the extra vertex $a^*$.)
Let $r$ satisfy 
$$krN=1/k-x.$$ 
Thus $r>0$.
For each $v\in V(G)$, define $w(v)$ as follows:
\begin{itemize}
\item $w(v)=(k-1)r(N+1)/N$ for $v\in \{a_1\ll a_N\}$; $w(v)=1/k-r$ for $v\in \{a_1'\ll a_{k-1}'\}$; 
\item $w(a^*)=1/k-N(k-1)r$;
\item $w(v)=(1-(k-1)s)/N$ for $v\in \{b_1\ll b_N\}$; $w(v)=s$ for $v\in \{b_1'\ll b_{k-1}'\}$; and
\item $w(v)=(1-(k-1)y)/N$ for $v\in \{c_1\ll c_N\}$; $w(v)=y$ for $v\in \{c_1'\ll c_{k-1}'\}$.
\end{itemize}
Then $(G,w)$ is a weighted graph. We claim it is $(x,y)$-biconstrained via $(A,B,C)$, and
$w(N^2_A(v))<1/k$ for each $v\in C$. To see this we must verify:
\begin{eqnarray*}
x&\le&p(1-(k-1)s)/N\\
y&\le &q(1-(k-1)y)/N\\
y&\le &q(1-(k-1)s)/N\\
x&\le&1/k-r\\
x&\le &p(k-1)r(N+1)/N+1/k-N(k-1)r, \text{ and}\\
1/k&>&1/k-N(k-1)r+(p+q-1)(k-1)r(N+1)/N.
\end{eqnarray*}
The first and third hold with equality from the definitions of $p,q$, and the second follows since $y\le s$. The fourth follows from
the definition of $r$. For the fifth, on substituting for $p$ and simplifying, we need to show that
$r(k-1)(N-x(N+1)/(1-(k-1)s))\le 1/k-x$, and this follows from the definition of $r$. Finally, the sixth
simplifies to $(p+q-1)(N+1)/N<N$, and this is true since $p+q\le N$. Consequently $\psi(x,y)<\frac{1}{k}$, by \ref{weighttonon}.
This proves \ref{psi12curve}.~\bbox

We also need the next three results:
\begin{thm}\label{psilb1}
Let $x',y',z'\in (0,1)$, with $\psi(x',y')\le z'$. If $x,y,z \in (0,1]$ satisfy $x \le 1/(2-x')$, $y \le y'/(1+y')$,
$x + (1-x')y/y' \le 1$, $z \ge 1/(2-z')$, and
$$x \le \frac{z-z'+x'(1-z)}{1-z'},$$
then $\psi(x,y) \le z$.
%It follows that if $z'=\psi(x',y')$ and $x,y \in (0,1]$ satisfy $x \le 1/(2-x')$, 
%$y \le y'/(1+y')$, $x < z$, $x + (1-x')y/y' \le 1$, $z > 1/(2-z')$, and \[x < \frac{z-z'+x'(1-z)}{1-z'}\] then $\psi(x,y) < z$.
\end{thm}
\Proof
Since $x'\le z'$ (because $\psi(x',y')\le z'$) it follows that
$$x\le 1/(2-x')\le 1/(2-z')\le z.$$
Let $G'$ be $(x,y)$-biconstrained via $(A,B,C)$, such that $|N_A^2(w)| \le z'|A|$ for all $w \in C$.
Add three vertices $a,b,c$ to the graph, and edges from $a$ to every vertex in $B$, edges from $b$ to every vertex in $A$,
and an edge between $b$ and $c$. Let this new graph be $G$.
Assign weights as follows:
\begin{eqnarray*}
w(a)&=& p\\
w(v)&=&(1-p)/|A| \text{ for each }v\in A\\
w(b)&=&q\\
w(v)&=& (1-q)/|B| \text{ for each }v\in B\\
w(c)&=& y\\
w(v)&=& (1-y)/|C| \text{ for each }v\in C.
\end{eqnarray*}
We will choose $p,q$ such that the weighted graph $(G,w)$ is $(x,y)$-biconstrained via $(A\cup \{a\}, B\cup \{b\}, C\cup \{c\})$
and $w(N_{A\cup \{a\}}^2(w)) \le z$ for all $w \in C\cup \{c\}$. The conditions
are: $1-p \ge x$, $1-q \ge x$, $q \ge y$, $(1-q)x'+q \ge x$, $(1-p)x'+p \ge x$, $(1-y)y' \ge y$, $(1-q)y' \ge y$, $1-p \le z$,
and $(1-p)z'+p \le z$. These are equivalent to the following:
\begin{align*}
    \max\left(1-z,\frac{x-x'}{1-x'}\right)\le p& \le \min\left(1-x, \frac{z-z'}{1-z'}\right) \\
    \max\left(y,\frac{x-x'}{1-x'}\right) \le q& \le \min\left(1-x,1-\frac{y}{y'}\right)
\end{align*}
Thus, it suffices to show that the lower bound on $p$ is at most the upper bound on $p$, and the same for $q$. We obtain
eight conditions, which simplify to those given in the theorem statement. This proves \ref{psilb1}.~\bbox

\begin{thm}\label{psilb2}
Let $x', y', z'\in (0,1)$, with $\psi(x',y')\le z'$. If $x,y \in (0,1]$ satisfy $y \le 1/(2-y')$, $x \le x'/(1+x')$, $x \le x'z$,
$(1-y')x/x' + y \le 1$, $z \ge 1/(2-z')$, and $x \le (z-z')/(1-z')$, then $\psi(x,y) \le z$.
\end{thm}
\Proof
Let $G'$ be $(x,y)$-biconstrained via $(A,B,C)$, such that $|N_A^2(w)| \le z'|A|$ for all $w \in C$. Add three vertices
$a,b,c$ to $G'$, with an edge from $a$ to $b$, edges from $b$ to every vertex in $C$, and edges from $c$ to every vertex in $B$.
Let this new graph be $G$. Assign weights as follows:
\begin{eqnarray*}
w(a)&=& p\\
w(v)&=&(1-p)/|A| \text{ for each }v\in A\\
w(b)&=&q\\
w(v)&=& (1-q)/|B| \text{ for each }v\in B\\
w(c)&=& 1-y\\
w(v)&=& y/|C| \text{ for each }v\in C.
\end{eqnarray*}
The conditions that the weighted graph $(G,w)$ is $(x,y)$-biconstrained via $(A\cup \{a\}, B\cup \{b|], C\cup \{c\})$
with $w(N_{A\cup \{a\}}A^2(w)) \le z$ for all $w \in C\cup \{c\}$ can be written as follows:
\begin{align*}
    \max\left(x,1-z\right) \le p& \le \min\left(1-\frac{x}{x'}, \frac{z-z'}{1-z'}\right) \\
    \max\left(x,\frac{y-y'}{1-y'}\right) \le q& \le \min\left(1-y,1-\frac{x}{x'}\right).
\end{align*}
We need to check that the lower bound for $p$ is at most the upper bound for $p$, and the same for $q$. This gives
eight conditions, which simplify (using that $1-y'>x'$, since $\psi(x',y')<1$) to those given in the theorem. This proves \ref{psilb2}.~\bbox

\begin{thm}\label{3.3gen}
Let $s,t\ge 1$ be integers with $s/t \le 1/2$. Let $x,y \in (0,1]$, satisfying $tx/s+y \le 1$, $x+ty/s \le 1$, and either $sy \le x$
or $sx \le y$. Furthermore, if either $x$ or $y$ is irrational, let strict inequality hold in all of these, that is, $tx/s+y < 1$,
$x+ty/s < 1$, and either $sy < x$
or $sx < y$. Then $\psi(x,y) < s/t$.
%\\ \\ In addition, if at least one of $x,y$ is irrational, and $s/t \le 1/2$, $tx/s+y < 1$, $x+ty/s < 1$, and either $sy < x$ or $sx < y$, then $\psi(x,y) < s/t$.  
\end{thm}
\Proof By increasing $x$ or $y$ if necessary, we may assume that $x,y$ are both rational.
Let $k+1 = \frac{t}{s}$. In terms of $k$, the hypotheses become $k\ge 1$, $(k+1)x+y \le 1$, $x+(k+1)y \le 1$, and either $sy \le x$
or $sx \le y$.

Suppose first that $sy \le x$.
Choose an integer $N \ge 1$ such that $p = xN / (1-kx)$ and $q = yN / (1-kx)$ are integers, and thus $p + q \le (x+y)N/(x+y) = N$.
Let $G_1$ be the graph with vertices $\{a_1,...,a_N, a^*, b_1,...,b_N, c_1,...,c_N\}$ where (reading subscripts modulo $N$)
each $a_i$ is adjacent to
$b_i,...,b_{i+p-1}$, each $b_i$ is adjacent to $c_i,...,c_{i+q-1}$, and $a^*$ is adjacent to all of the $b_i$.

Let $m = t-s$. Let $G_2$ be the graph with vertex set $\{a_1',...,a_m', b_1',...,b_m', c_1',...,c_m'\}$, where each
$a_i'$ is adjacent to $b_i',...,b_{i+s-1}'$ (reading subscripts modulo $m$), and each $b_i'$ is adjacent to $c_i'$.
Let $G$ be the disjoint union of $G_1$ and $G_2$. Let $A=\{a_1,...,a_N, a^*,a_1',...,a_m'\}$, and
$B=\{b_1,...,b_N, b_1',...,b_m'\}$ and define $C$ similarly.

Let $r$ satisfy $(k+1)rN = \frac{1}{k+1}-x$. Assign weights as follows:
\begin{eqnarray*}
w(a_i)&=& kr(N+1)/N\\
w(a_i')&=&1/b-r/s\\
w(a^*)&=&1/(k+1)-Nkr\\
w(b_i)&=& (1-kx)/N\\
w(b_i')&=& x/s\\
w(c_i)&=& (1-ksy)/N\\
w(c_i') &=& y.
\end{eqnarray*}
This defines a weighted graph $(G,w)$, $(x,y)$-biconstrained via $(A,B,C)$, such that $w(N_A^2(v)) < s/t=1/(k+1)$ for all $v \in C$,
and so
$\psi(x,y) < s/t$, as desired.

Now suppose $sy > x$, and consequently $sx \le y$. Choose an integer $N \ge 1$ such that $p=xN/(1-ky)$, $q = yN/(1-ky)$ are integers,
and thus $p+q \le N$. Let $G_1$ be as before. Let $m = b-a$, and let $G_2$ be the graph with vertex set
$\{a_1',...,a_m', b_1',...,b_m', c_1',...,c_m'\}$, where each $a_i'$ is adjacent to $b_i'$, and each $b_i'$ is adjacent to
$c_i',...,c_{i+a-1}'$, reading subscripts modulo $m$ (thus, this is the earlier graph $G_2$ flipped).
Let $G$ be the disjoint union of $G_1$ and $G_2$, and define $A,B,C$ as before. Let $r>0$ satisfy $(k+1)rN \le 1/(k+1)-x$ and
$r \le 1/(k+1)-y$. Assign weights as follows:

\begin{eqnarray*}
w(a_i)&=& kr(N+1)/N\\
w(a_i')&=&1/t-r/s\\
w(a^*)&=&1/(k+1)-Nkr\\
w(b_i)&=& (1-ky)/N\\
w(b_i')&=& y/s\\
w(c_i)&=& (1-ky)/N\\
w(c_i') &=& y/s.
\end{eqnarray*}
Then $(G,w)$ is a weighted graph, and is $(x,y)$-biconstrained via $(A,B,C)$, and $w(N_A^2(v)) < s/t$ for all $v \in C$, showing
that
$\psi(x,y) < s/t$. This proves \ref{3.3gen}.~\bbox

\begin{thm}\label{philb}
Let $x,y,z\in (0,1]$, with $y \le 1/2<x$. If $\phi\left(2-1/x,y/(1-y)\right) \le 2-1/z$, then $\phi(x,y) \le z$.
\end{thm}
\Proof
Let $x'=(2x-1)/x$, and $y'=y/(1-y)$, and let $z'= \phi\left(2-1/x,y/(1-y)\right)$. Let $G'$ be a graph that is
$(x,y)$-constrained via $(A,B,C)$, such that
$|N_A^2(w)| \le z'|A|$ for all $w \in C$. Add three vertices $a,b,c$ to the graph, and edges from $a$ to every vertex
in $B$, edges from $b$ to every vertex in $A$, and an edge between $b$ and $c$. Let this new graph be $G$.

Assign weights $w(v)\;(v\in V(G))$ as follows:
\begin{eqnarray*}
w(a)&=& 1-z\\
w(v)&=&z/|A| \text{ for each }v\in A\\
w(b)&=&1-x\\
w(v)&=& x/|B| \text{ for each }v\in B\\
w(c)&=& y\\
w(v)&=& (1-y)/|C| \text{ for each }v\in C.
\end{eqnarray*}
Then the weighted graph $(G,w)$ is $(x,y)$-constrained via $(A\cup \{a\},B\cup \{b\},C\cup \{c\})$, since $xx'+(1-x) = x$ and $(1-y)y'=y$. Moreover,
$w(N_A^2(w)) \le (1-z)+zz'\le z$
for all $w\in C$; and $w(N_{A\cup \{a\}}^2(c))=z$.
Thus $\phi(x,y) \le z$. This proves \ref{philb}.~\bbox

\section{Biconstrained graphs}
In this section we prove some lower bounds on $\psi(x,y)$. On the diagonal $x=y$, $\psi(x,y)$ behaves 
perfectly; it turns out that for all $x$, $\psi(x,x) = 1/k$, where $k$ is the largest integer with $1/k\ge x$. 
That follows from:
\begin{thm}\label{bisym}
For all integers $k\ge 1$,
if $x,y\in (0,1]$ with $x+ky>1$ and $kx+\frac{x}{1-(k-1)y}\ge 1$, then $\psi(x,y)\ge 1/k$.
\end{thm}
\Proof
By \ref{maxbound} we may assume that $x,y<1/k$.
Let $G$ be $(x,y)$-biconstrained, via $(A,B,C)$. 
We must show that $|N^2_A(v)|\ge |A|/k$ for some $v\in C$. Suppose not.
Choose $K\subseteq C$ with $|K|\le k$, and subject to that with $|K|$ maximum such that the sets $N(v)\;(v\in K)$ are
pairwise disjoint. Let
$I\subseteq A$ be the union of the sets $N^2_A(v)\;(v\in K)$, and let $J\subseteq B$ be the union of the sets $N(v)\;(v\in K)$.
It follows that
\\
\\
(1) {\em $|A\setminus I|> (1-|K|/k)|A|$, and $|B\setminus J|\le (1-|K|y)|B|$.}

\bigskip
If $|K|=k$, then by (1), $|B\setminus J|\le (1-ky)|B|< x|B|$, and since every vertex in $A$ has $x|B|$ neighbours in $B$,
it follows that every vertex in $A$ has a neighbour in $J$, that is, $I=A$, contrary to~(1). Thus $|K|<k$.

Since each vertex in $A\setminus I$ has at least
$x|B|$ neighbours in $B$, and they all belong to $B\setminus J$, some vertex $t\in B\setminus J$ has at least 
$$x|B|\frac{|A\setminus I|}{|B\setminus J|}\ge x\frac{1-|K|/k}{1-|K|y}|A|$$ 
neighbours
in $A\setminus I$ by (1). Since $|K|\le k-1$ and $ky<1$, and therefore $|K|y<1$, it follows that
$$\frac{1-|K|/k}{1-|K|y}\ge \frac{1-(k-1)/k}{1-(k-1)y}= \frac{1}{k(1-(k-1)y)},$$
and so $t$ has at least $\frac{x|A|}{k(1-(k-1)y)}$ neighbours in $A\setminus I$.
Let $u\in C$ be adjacent to $t$.
From the maximality of $K$,
$u$ has a neighbour $w\in N(v)$ for some $v\in K$. Since $w$ has at least $x|A|$ neighbours
in $I$, it follows that
$$|N^2_A(u)|\ge x|A|+\frac{x|A|}{k(1-(k-1)y)}\ge |A|/k,$$
a contradiction.
This proves \ref{bisym}.~\bbox

We deduce:
\begin{thm}\label{bisym2}
For all integers $k\ge 1$,
if $x,y\in (0,1]$ with $x+ky>1$ and $kx+y\ge 1$, then $\psi(x,y)\ge 1/k$.
\end{thm}
\Proof
If $k=1$ the result is easy (and follows from \ref{trivial} below), so we assume that $k\ge 2$; and hence we may assume that 
$x,y< 1/k\le 1/2$ by \ref{maxbound}. By \ref{bisym} we may assume (for a contradiction) that $kx+\frac{x}{1-(k-1)y}< 1$.
Consequently $kx+\frac{x}{1-(k-1)(1-kx)}< 1$. Let $t=1-kx$. Then $\frac{1-t}{1-(k-1)t}< kt$, and so
$k(k-1)t^2 -(k+1)t+1<0$. This is quadratic in $t$, with discriminant $(k+1)^2-4k(k-1)$, and the latter is negative if $k>2$;
so we may assume that $k=2$. Then $2t^2 -3t+1<0$, so $(2t-1)(t-1)<0$, that is, $1/2<t<1$. But $t=1-2x$, so $1/2<1-2x<1$, that is,
$x<1/4$. But  $2x+y\ge 1$ and $y<1/2$, a contradiction. This proves \ref{bisym2}.~\bbox

Consequently we have:
\begin{thm}\label{bisym3}
For all $x\ge 0$, $\psi(x,x) = 1/k$, where $k$ is the largest integer with $1/k\ge x$.
\end{thm}
\Proof
Certainly $\psi(x,x)\le 1/k$, since by \ref{trivialbounds},
$$\psi(x,x)\le \frac{\lceil kx\rceil +\lceil kx\rceil -1}{k}= 1/k.$$
Equality holds by \ref{bisym2}. This proves \ref{bisym3}.~\bbox

Next we need a lemma:
\begin{thm}\label{backpsi}
Let $k\ge 1$ be an integer, let $(k-1)/k^2\le y\le 1$, and let $(A,B,C)$ be a tripartition of a graph $G$, such that:
\begin{itemize}
\item every vertex in $B$ has at least $y|C|$ neighbours in $C$; and
\item $|N^2_A(v)|<|A|/k$ for each $v\in C$. 
\end{itemize}
Then there exist $v_1\ll v_k\in A$ such that $N(v_i)\cap N(v_j)=\emptyset $ for $1\le i<j\le k$.
\end{thm}
\Proof If some vertex $v$ in $A$ has degree zero, then we may take $v_1=\cdots=v_k=v$. So we assume that
every vertex in $A$ has a neighbour in $B$.
For each $v\in A$, let $c(v)=|N^2_C(v)|$, and let $A(v)\subseteq A$ be the set of vertices in $A$ that have a neighbour in
$N(v)$. Let $|A(v)|=a(v)$. 
\\
\\
(1) {\em For each $v\in A$, $c(v)>kya(v)|C|/|A|$.}
\\
\\
If we choose $u\in N^2_C(v)$ independently at random, then since every vertex in $A(v)$
has at least $y|C|$ second neighbours in $N^2_C(v)$, the probability that a given vertex $w\in A(v)$ belongs to $N^2_A(u)$
is at least $y|C|/c(v)$, and so the expectation of $|N^2_A(u)|$ is at least $(y|C|/c(v))a(v)$. On the other hand, the 
expectation of $|N^2_A(u)|$ is less than $|A|/k$.
This proves (1).

\bigskip
Let $H$ be the graph with vertex set $A$, in which distinct $u,v$ are adjacent if (in $G$) $u,v$ have a common neighbour in $B$.
Thus every vertex $v$ has degree $a(v)-1$ in $H$. So
$2|E(H)|=\sum_{v\in A} (a(v)-1)$;
but
$$(ky|C|/|A|)\sum_{v\in A}a(v)\le \sum_{v\in A}c(v)=\sum_{v\in A}|N^2_C(v)|=\sum_{u\in C}|N^2_A(u)|<|A|\cdot|C|/k.$$
Consequently 
$$2|E(H)|<(|A|\cdot|C|/k)/(ky|C|/|A|)- |A|=|A|^2/(k^2y)- |A|\le |A|^2/(k-1)- |A|.$$
By Tur\'an's theorem, $H$ has a stable set of cardinality $k$. This proves \ref{backpsi}.~\bbox

\begin{thm}\label{semibi}
Let $k\ge 1$ be an integer, and let $x,y\in (0,1]$ where $y\ge (k-1)/k^2$ and $kx+y>1$.
Let $G$ be $(x,y)$-constrained via $(A,B,C)$, such that
every vertex in $C$ has at least $y|B|$ neighbours in $B$.
Then $|N^2_A(v)|\ge |A|/k$ for some  $v\in C$. Consequently $\psi(x,y)\ge 1/k$.
\end{thm}
\Proof Suppose not; then there is a weighted graph $(G',w)$, $(x,y)$-constrained via some tripartition $(A',B',C')$,
such that 
\begin{itemize}
\item for each $v\in C'$, $w(N(v))\ge y|B'|$; and
\item for each $v\in C'$, $w(N^2_{A'}(v))<1/k$.
\end{itemize}
Choose such a weighted graph $(G',w)$ with $|V(G')|$ minimum, and let $z<1/k$ such that
$w(N^2_{A'}(v))\le z$ for each $v\in C'$. By \ref{backpsi}, there exist $v_1\ll v_k\in A'$
such that $N(v_1)\ll N(v_k)$ are pairwise disjoint. Consequently 
$w(N(v_1)\cup\cdots\cup N(v_k))\ge kx$; 
and since $w(N(u))\ge y>1-kx$ for each $u\in C'$, it follows that
$\bigcup_{v\in X}N^2_{C'}(v) = C'$ where $X=\{v_1\ll v_k\}$. But $|X|<z^{-1}$, contrary to \ref{nodom} and the minimality of $|V(G')|$.
This proves the first claim, and the second follows. This proves \ref{semibi}.~\bbox

It is awkward to express the biconstrained problem in the language of triangular triples, but we can do so as follows.
For $x,y,z\in (0,1]$ we say that $(x^*,y,z)$ is {\em triangular} if no triangle-free graph $G$
admits a tripartition $(A,B,C)$ that satisfies the three bullets of the previous definition, and in addition
satisfies
\begin{itemize}
\item every vertex in $B$ has at least $x|A|$ neighbours in $A$.
\end{itemize}
Similarly, we say $(x^*,y^*,z)$ is triangular if no triangle-free graph $G$
admits a tripartition $(A,B,C)$ that satisfies the three bullets of the previous definition, and in addition
satisfies
\begin{itemize}
\item every vertex in $B$ has at least $x|A|$ neighbours in $A$; and
\item every vertex in $C$ has at least $y|B|$ neighbours in $B$;
\end{itemize}
and so on.
Then we have:
\begin{thm}\label{bireformtri}
For $x,y,z\in (0,1]$, $\psi(x,y)> 1-z$ if and only if $(x^*,y^*,z)$ is triangular.
\end{thm}

This is not the same as saying that $(x^*,y^*,z^*)$ is triangular, so we need to keep track of the asterisks
if we rotate; but still it can be useful, as we shall see.

\section{The mono-constrained case}\label{sec:mono}

In this section we are mostly concerned with $\phi(x,y)$ when $x=y$.
We know that $\psi$ behaves well on the diagonal $x=y$, because of \ref{bisym3}, so what about $\phi$? 
More generally, what about an analogue of \ref{bisym} or \ref{bisym2}
with $\psi$ replaced by $\phi$?

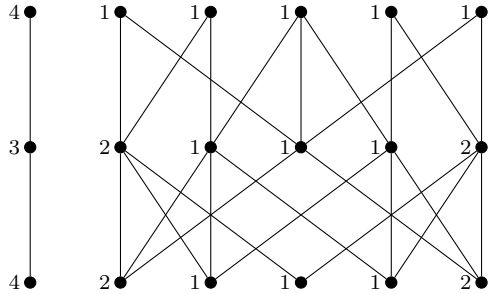
\begin{figure}[ht]
\centering

\begin{tikzpicture}[yscale=0.6,xscale=0.6,auto=left]
\tikzstyle{every node}=[inner sep=1.5pt, fill=black,circle,draw]

\def\r{0}
\node (a1) at (10,\r) {};
\node (a2) at (8,\r) {};
\node (a3) at (6,\r) {};
\node (a4) at (4,\r) {};
\node (a5) at (2,\r) {};
\node (a6) at (0,\r) {};
\def\r{-3}
\node (b1) at (10,\r) {};
\node (b2) at (8,\r) {};
\node (b3) at (6,\r) {};
\node (b4) at (4,\r) {};
\node (b5) at (2,\r) {};
\node (b6) at (0,\r) {};
\def\r{-6}
\node (c4) at (10,\r) {};
\node (c5) at (8,\r) {};
\node (c3) at (6,\r) {};
\node (c1) at (4,\r) {};
\node (c2) at (2,\r) {};
\node (c6) at (0,\r) {};
\foreach \from/\to in {a1/b1,a1/b3,a2/b1,a2/b2,a3/b2,a3/b3,a3/b4,a4/b4,a4/b5,a5/b5,a5/b3,a6/b6,b1/c5,b1/c4,b1/c3,b2/c5,b2/c4,b2/c1,b3/c2,b3/c4,b4/c5,b4/c2,b4/c1,b5/c3,b5/c2,b5/c1,b6/c6}
\draw [-] (\from) -- (\to);
\tikzstyle{every node}=[left]
\draw (b1) node []           {\scriptsize$2$};
\draw (b2) node []           {\scriptsize$1$};
\draw (b3) node []           {\scriptsize$1$};
\draw (b4) node []           {\scriptsize$1$};
\draw (b5) node []           {\scriptsize$2$};
\draw (b6) node []           {\scriptsize$3$};
\draw (c1) node []           {\scriptsize$1$};
\draw (c2) node []           {\scriptsize$2$};
\draw (c3) node []           {\scriptsize$1$};
\draw (c4) node []           {\scriptsize$2$};
\draw (c5) node []           {\scriptsize$1$};
\draw (c6) node []           {\scriptsize$4$};
\draw (a1) node []           {\scriptsize$1$};
\draw (a2) node []           {\scriptsize$1$};
\draw (a3) node []           {\scriptsize$1$};
\draw (a4) node []           {\scriptsize$1$};
\draw (a5) node []           {\scriptsize$1$};
\draw (a6) node []           {\scriptsize$4$};
\end{tikzpicture}

\caption{$\phi(3/10,4/11)\le 4/9$.} \label{fig:3-10:4-11}
\end{figure}
If we replace $\psi$ by $\phi$ in \ref{bisym}, it
becomes false, even with $k=2$, because $\phi(3/10, 4/11)\le 4/9$,
as the graph of figure \ref{fig:3-10:4-11} shows (the sets $A,B,C$ are the rows,
and the numbers on the vertices are used as in figure \ref{fig:exactcount}).
But as far as we know, \ref{bisym2} might hold with $\psi$ replaced by $\phi$.
Let us state this as a conjecture:
\begin{thm}\label{monosymconj}
{\bf Conjecture: }
For all integers $k\ge 1$,
if $x,y\in (0,1]$ with $x+ky>1$ and $kx+y\ge 1$, then $\phi(x,y)\ge 1/k$.
\end{thm}
On the other hand, we have not even been able to prove what is presumably the simplest nontrivial case of this, namely that $\phi(x,y)\ge 1/2$
for all $x,y$ with $x,y>1/3$.  But we do have several results approaching \ref{monosymconj}.
First, it is true with $k=1$; we have the trivial:

\begin{thm}\label{trivial}
For $x,y\in (0,1]$, if $x+y>1$, or $x+y=1$ and $x$ is irrational, then $\phi(x,y)=1$.
\end{thm}
\Proof
Let $G$ be $(x,y)$-constrained via $(A,B,C)$. Then some vertex $v\in C$ has at least $y|B|$ neighbours in $B$, and strictly
more if $y$ is irrational; and so $N_A^2(v)=A$, as every
vertex in $A$ has at least $x|B|$ neighbours in $B$. This proves \ref{trivial}.~\bbox

This is tight, in that if $x+y=1$ and $x,y$ are rational, then $\phi(x,y)<1$. We omit the proof, which is easy.

\ref{monosymconj} implies that $\phi(x,x)\ge 1/2$ if $x>1/3$. We have not been able to prove this, but we can show that
$\phi(x,x)>3/7$ if $x>1/3$. That is implied by the following:

\begin{thm}\label{3-7}
Let $k\ge 2$ be an integer; then for $x,y\in (0,1]$, if $y>1/k$ then 
$$\phi(x,y)\ge \frac{x(2-3x)}{kx(1-x)+x^2-3x+1}.$$ 
Indeed if $k=2$, then $\phi(x,y)\ge 2x-x^2$ (which is larger).
\end{thm}
\Proof
Let $G$ be $(x,y)$-constrained via $(A,B,C)$. If $x$ is irrational 
then $G$ is $(x,y)$-constrained via $(A,B,C)$,
for some rational $x'>x$; so we may assume that $x$ is rational, by increasing $x$ if necessary.
Suppose that $k=2$, and choose $v_1,v_2\in B$ independently and uniformly  at random. For each $u\in A$,
the probability that $u$ is adjacent to at least one of $v_1,v_2$ is at least $2x-x^2$, since $u$ has at least $x|B|$
neighbours in $B$; and so we may choose $v_1,v_2$ such that at least $(2x-x^2)|A|$ vertices in $A$ are adjacent to at least
one of them. But $v_1,v_2$ have a common neighbour in $C$, since $y>1/2$, and the claim follows.

Thus we may assume that $k\ge 3$. By \ref{maxbound}, $\phi(x,y)\ge x$, and so we may assume that
$$\frac{x(2-3x)}{kx(1-x)+x^2-3x+1}>x,$$ 
that is, $x<1/(k-1)$. Consequently $x\le (k-2)/(k-1)$ since $k\ge 3$.
Define 
\begin{eqnarray*}
p&=&\frac{x(1-x)}{kx(1-x)+x^2-3x+1},\\
s&=&\frac{x}{(k-2)(1-x)}, \text { and }\\
m&=&\frac{x(2-3x)}{kx(1-x)+x^2-3x+1}.
\end{eqnarray*}
These are all non-negative, and 
$p$ is rational with denominator $T$ say; and by replacing each vertex by $T$ copies,
we may assume that $p|A|$
is an integer. Since $x\le (k-2)/(k-1)$ it follows that $s\le 1$.

For $1\le i\le k-1$, we define $v_i\in B$, and a subset $P_i$ of $N_A(v_i)$ with $|P_i|=p|A|$, inductively, as follows. 
Let $Q=P_1\cup\cdots\cup P_{i-1}$. 
\\
\\
(1) {\em There exists $v_i\in B$ such that $sa+b\ge x(s|Q|+|A|-|Q|)$, where $a=|N_A(v_i)\cap Q|$ and
$b=|N_A(v_i)\setminus Q|$.}
\\
\\
Suppose not; then summing over all $v\in B$, we deduce that
$$\sum_{v\in B}s|N_A(v)\cap Q| +\sum_{v\in B}|N_A(v)\setminus Q| < x(s|Q|+|A|-|Q|)|B|.$$
But the first sum is $s$ times the number of edges between $Q$ and $B$, and so at least $xs|B|\cdot |Q|$; and the second
is similarly at least $x|B|(|A|-|Q|)$, a contradiction. This proves (1).

\bigskip
Let $v_i$ be as in (1). Thus $sa+b\ge x(s|Q|+|A|-|Q|)\ge x(1-(1-s)(k-2)p)|A|$. In particular, since
$$a+b\ge sa+b\ge x(1-(1-s)(k-2)p)|A|= p|A|,$$ 
there exists $P_i\subseteq N_A(v_i)$ of cardinality $p|A|$.
Also, since $a\le (k-2)p|A|$, and so 
$$s(k-2)p|A|+b\ge sa+b\ge x(1-(1-s)(k-2)p)|A|,$$
it follows that 
$$b\ge x(1-(1-s)(k-2)p)|A|-s(k-2)p|A|=(m-p)|A|,$$
and so
$$|N_A(v_h)\cup N_A(v_i)|\ge m|A|,$$
for $1\le h< i$. This completes the inductive definition of $v_1\ll v_{k-1}$ and $P_1\ll P_{k-1}$.

Let $P=P_1\cup\cdots\cup P_{k-1}$. Then $|P|\le (k-1)p|A|$.
Since every vertex in $A\setminus P$ has at least $x|B|$ neighbours in $B$, there exists $v_k\in B$
with at least $x(|A|-|P|)\ge x(1-(k-1)p)|A|$ neighbours in $A\setminus P$. Let $P_k$ be its set of neighbours in
$A\setminus P$. Then for all $i$ with $1\le i\le k-1$, 
$$|P_i|+|P_k|\ge (x(1-(k-1)p)+p)|A|=m|A|.$$ 
Consequently, for all distinct
$v,v'\in \{v_1\ll v_k\}$, $|N_A(v)\cup N_A(v')|\ge m|A|$. But since $y>1/k$, some two of $v_1\ll v_k$
have a common neighbour $u\in C$, and so $|N^2_A(u)|\ge m$. This proves \ref{3-7}.~\bbox

We deduce from \ref{3-7} a version of \ref{bisym3} for the mono-constrained case:
\begin{thm}\label{5-12}
For $y\in (0,1]$, if $y>1/k$ where $k\ge 2$ is an integer, then 
$\phi(1/k,y)\ge \frac{2k-3}{2k^2-4k+1}$.
\end{thm}
Consequently $\phi(1/k,y)\ge 1/k + 1/(2k^2) +\Omega(k^{-3})$.

\ref{5-12} tells us in particular that $\phi(x,x)\ge \frac{2k-3}{2k^2-4k+1}>1/k$ when $x>1/k$ (if $k\ge 2$ is an integer), and since 
$\phi(1/k,1/k) = 1/k$, there is a discontinuity in $\phi(x,x)$ when $x=1/k$, 
and the limit of $\phi(x,x)$ as $x\rightarrow 1/k$ from above is different from $\phi(1/k,1/k)$.
What happens when $x\rightarrow 1/k$ from below? The next results investigate this. We will show that 
if $x$ is sufficiently close to $1/k$ from below, then $\phi(x,x)=1/k$.

\begin{thm}\label{limk1}
If $k>0$ is an integer and $x\in (0,1]$ satisfies $(1-x)^k<x$, then $\phi(x,x)\ge 1/k$. In particular,
if $x>0.382$ then $\phi(x,x)\ge 1/2$, and if $x>0.318$ then $\phi(x,x)>1/3$.
\end{thm}
\Proof
Let $G$ be $(x,x)$-constrained via $(A,B,C)$. If we choosing $k$ vertices from $C$ uniformly at random,
the number of vertices in $B$ nonadjacent to all of them is at most $(1-x)^k|B|$ in expectation; and so there
exist $v_1\ll v_k\in C$ such that at most $(1-x)^k|B|$ vertices in $B$ are nonadjacent to all of them.
Since $(1-x)^k|B|<x|B|$, it follows that the sets $N^2_A(v_i)\;(1\le i\le k)$ have union $A$, and so one of them
has cardinality at least $|A|/k$. This proves \ref{limk1}.~\bbox

The proof of \ref{limk1} is very simple, but the result is not of any value.
It is of no use when $k\ge 4$ because then $(1-x)^k<x$ implies $x>1/k$; and
we will prove in \ref{improvedphi12} and \ref{phi13monster} that $\phi(x,x)\ge 1/2$ when $x>0.352202$, and $\phi(x,x)\ge 1/3$ 
when $x\ge 0.28231$, which are stronger than
\ref{limk1} when $k=2,3$.
Here is another approach to the same question, more successful for larger values of $k$.

\begin{thm}\label{limk}
Let $k\ge 1$ be an integer, and let $x\ge 1/k-\vare$ where $\vare=1/(13k^3)$. Then $\phi(x,x)\ge 1/k$.
%Assume $(y-c)/k>(c+\vare)$
%Assume $1/(k+1)> 3(c+\vare)k+c$
%Assume $(k+1)y>1$
%Assume $(y-c)/(1-c)> k\vare/(1+k\vare)$
\end{thm}
\Proof
We may assume that $x=1/k-\vare$. By \ref{trivial} we may assume that $k\ge 2$. We leave the reader to check that
\begin{itemize}
\item $1/(2k)-\vare> 6k^2\vare$;
\item $x>1/(k+1)$; and
\item $(2kx-1)/(2k-1)> (k\vare)/(x+k\vare)$.
\end{itemize}
(These are inequalities we will need later.)
Let $G$ be $(x,x)$-constrained via $(A,B,C)$, and suppose that $|N^2_A(v)|< |A|/k$ for each $v\in C$.
Let
$P$ be the set of vertices in $B$ that have at most $(1/k-2k\vare)|A|$
neighbours in $A$. 
\\
\\
(1) {\em $|P|\le |B|/(2k)$.}
\\
\\
Every vertex in $B$ has fewer than $|A|/k$ neighbours in $A$, and so the number of edges between $A$ and $B$
is at most $|P|(1/k-2k\vare)|A|+(|B|-|P|)|A|/k$. On the other hand, the number of such edges is at least $(1/k-\vare)|A|\cdot |B|$; and so
$$|P|(1/k-2k\vare)|A|+(|B|-|P|)|A|/k\ge (1/k-\vare)|A|\cdot |B|,$$ 
which simplifies to 
$2k|P|\le |B|.$ This proves (1).
\\
\\
(2) {\em If $u,v\in B\setminus P$ have a common neighbour in $C$, then $|N_A(u)\setminus N_A(v)|\le 2k\vare|A|$.}
\\
\\
Since $u,v\in B\setminus P$ have a common neighbour in $C$, it follows that $|N_A(u)\cup N_A(v)|\le |A|/k$. But
$|N_A(u)|\ge (1/k-2k\vare)|A|$ since $u\in B\setminus P$, and so $|N_A(u)\setminus N_A(v)|\le 2k\vare|A|$. This proves (2).
\\
\\
(3) {\em There exist $v_1\ll v_k\in B\setminus P$ such that for $1\le i<j\le k$, there
are at least $(1/(2k)-\vare)|A|/k$ vertices in $A$ that are adjacent to $v_j$ and not to $v_i$.}
\\
\\
Choose $v_1\ll v_k\in B\setminus P$ as follows. Choose $v_1\in B\setminus P$ arbitrarily.
Inductively, suppose we have defined $v_1\ll v_i$ where $i<k$. Each has at most $|A|/k$
neighbours in $A$, and so the set of vertices in $A$ adjacent to
one of $v_1\ll v_i$ has cardinality at most $(i/k)|A|\le (1-1/k)|A|$. Let $D$ be the set of vertices in $A$ nonadjacent to each of $v_1\ll v_i$;
then $|D|\ge |A|/k$. Since, by (1), each vertex in $D$ has at least $x|B|-|P|\ge (1/(2k)-\vare)|B|$
neighbours in $B\setminus P$, there exists $v_{i+1}\in B\setminus P$ with at least
$(1/(2k)-\vare)|A|/k$ neighbours in $D$. This completes the inductive definition. We see that for $1\le i<j\le k$, there
are at least $(1/(2k)-\vare)|A|/k$ vertices in $A$ that are adjacent to $v_j$ and not to $v_i$. This proves (3).

\bigskip
Let $H$ be the bipartite graph $G[(B\setminus P)\cup C]$.
\\
\\
(4) {\em For $1\le i<j\le k$, $v_i$ and $v_j$ belong to distinct components of $H$.}
\\
\\
From (2), the sets $N_C(v_1)\ll N_C(v_k)$ are pairwise disjoint, because $(1/(2k)-\vare)|A|/k>2k\vare|A|$.
Suppose that there is a path of $H$ joining some two of $v_1\ll v_k$, and take the shortest such path $Q$;
between $v_i$ and $v_j$ say, where $j>i$. Let $Q$ have
$m$ vertices in $B$, say $u_1\ll u_m$ in order where $u_1=v_i$. We claim that $m\le 4$. For suppose that $m\ge 5$.
From the minimality of the length of $Q$, $u_3$ has no common neighbour in $C$
with any of $v_1\ll v_k$, and so the sets $N_C(u_3), N_C(v_1)\ll N_C(v_k)$ are pairwise disjoint, which is impossible
since $x>1/(k+1)$. Thus $m\le 4$.
By applying (2) to each pair of consecutive members of $V(Q)\cap B$, we deduce that
$$|N_A(v_j)\setminus N_A(v_i)|\le (m-1)2k\vare|A|\le 6k\vare|A|.$$
But $|N_A(v_j)\setminus N_A(v_i)|\ge (1/(2k)-\vare)|A|/k$, and so
$(1/(2k)-\vare)|A|/k\le 6k\vare|A|$,
a contradiction. This proves~(4).

\bigskip

For $1\le i\le k$, let $H_i$ be the component of $H$ containing $v_i$, and let $V(H_i)\cap B=B_i$
and $V(H_i)\cap C=C_i$. If there exists $v\in B\setminus P$ that does not belong to any of $B_1\ll B_k$, then
the sets $N_C(v), N_C(v_1)\ll N_C(v_k)$ are pairwise disjoint, which is impossible since they all have
cardinality at least $x|C|$, and $(k+1)x>1$. Consequently the sets $B_1\ll B_k$ and $P$ form a partition of $B$.
\\
\\
(5) {\em For $1\le i\le k$ there exists $u_i\in C_i$ adjacent to at least $\frac{(1-k\vare)}{1+ k(k-1)\vare}|B_i|$ 
vertices in $B_i$.}
\\
\\
For $1\le i\le k$, since $v_i$ has at least $x|C|$ neighbours in $C$, it follows that $|C_i|\ge x|C|$.
Let $1\le i\le k$. Since $C_1\ll C_k$ are pairwise disjoint, and the union of the sets $C_j\;(j\in \{1\ll k\}\setminus \{i\})$
has cardinality at least $(k-1)x|C|$, it follows that 
$$|C_i|\le |C|-(k-1)x|C|= x|C|+ k\vare|C|.$$
There are at least $x|B_i|\cdot |C|$ edges between $B_i$ and $C_i$, and so some vertex in $C_i$
has at least
$$x(|C|/|C_i|)|B_i|\ge x(|C|/( x|C|+ k\vare|C|))|B_i|= (x/( x+ k\vare))|B_i|$$
neighbours in $B_i$. By substituting $x=1/k-\vare$, this proves (5).

\bigskip
For $1\le i\le k$, let $A_i=N^2_A(u_i)$. Since $|A_i|<|A|/k$ for $1\le i\le k$, there exists $v\in A$
that is in none of $A_1\ll A_k$. Now $v$ has at least $x|B|$ neighbours in $B$, and they all
belong to $B_1\cup\cdots\cup B_k$ except for at most $|P|$ of them. Consequently there exists $i\in \{1\ll k\}$
such that $v$ has at least $(x|B|-|P|)|B_i|/|B\setminus P|$ neighbours in $B_i$. Since $v\notin A_i$,
it follows that
$$(x|B|-|P|)|B_i|/|B\setminus P|+(x/(x+ k\vare))|B_i|\le  |B_i|.$$
Since $x|B|\le |B|$ and $|P|\le |B|/(2k)$ by (1), it follows that 
$$(x|B|-|P|)|B_i|/|B\setminus P|\ge (x-1/(2k))|B_i|/(1-1/(2k)) = (2kx-1)|B_i|/(2k-1),$$
and so
$(2kx-1)/(2k-1)\le  k\vare/(x+k\vare)$,
a contradiction.
This proves \ref{limk}.~\bbox

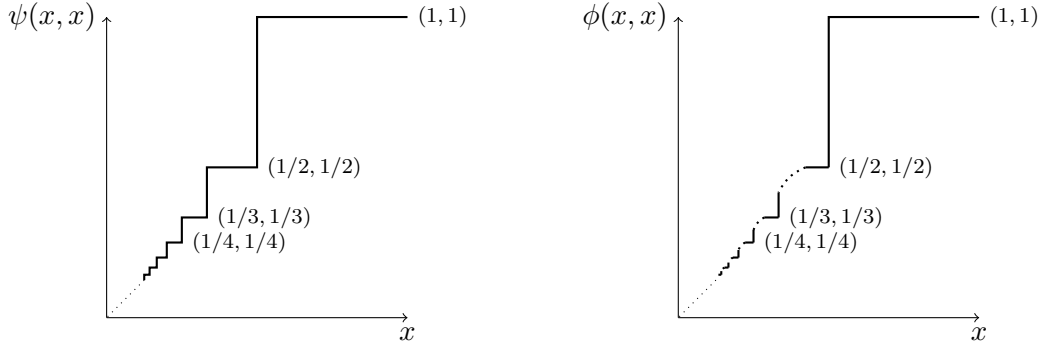
\begin{figure}[ht]
\centering

\begin{tikzpicture}[scale=.4,auto=left]

\def\l{-15}
\draw[->] (\l,0) -- (\l+10,0) node[anchor=north]{$x$};
\draw[->] (\l,0) -- (\l,10) node[anchor=east] {$\psi(x,x)$};
\node[right] at (\l+10,10) {\begin{scriptsize}$(1,1)$\end{scriptsize}};
\node[right] at (\l+5,5) {\begin{scriptsize}$(1/2,1/2)$\end{scriptsize}};
\node[right] at (\l+10/3,10/3) {\begin{scriptsize}$(1/3,1/3)$\end{scriptsize}};
\node[right] at (\l+10/4,10/4) {\begin{scriptsize}$(1/4,1/4)$\end{scriptsize}};
\draw[thick] (\l+10,10) -- (\l+5,10) -- (\l+5,5) -- (\l+10/3,5)-- (\l+10/3,10/3) -- (\l+10/4,10/3) -- (\l+10/4,10/4) -- (\l+10/5,10/4) -- (\l+10/5,10/5) -- (\l+10/6,10/5) -- (\l+10/6,10/6)
 -- (\l+10/7,10/6) -- (\l+10/7,10/7) -- (\l+10/8,10/7) -- (\l+10/8,10/8);
\draw[dotted] (\l+10/8,10/8) -- (\l,0);

\def\l{4}
\draw[->] (\l,0) -- (\l+10,0) node[anchor=north]{$x$};
\draw[->] (\l,0) -- (\l,10) node[anchor=east] {$\phi(x,x)$};
\node[right] at (\l+10,10) {\begin{scriptsize}$(1,1)$\end{scriptsize}};
\node[right] at (\l+5,5) {\begin{scriptsize}$(1/2,1/2)$\end{scriptsize}};
\node[right] at (\l+10/3,10/3) {\begin{scriptsize}$(1/3,1/3)$\end{scriptsize}};
\node[right] at (\l+10/4,10/4) {\begin{scriptsize}$(1/4,1/4)$\end{scriptsize}};
\draw[thick] (\l+10,10) -- (\l+5,10) -- (\l+5,5);
\draw[thick] (\l+5,5) -- (\l+4.2265,5);
\draw[thick] (\l+10/3,10/3) -- (\l+10/3,50/12);
\draw[thick] (\l+10/3,10/3) -- (\l + 10/3.5,10/3);
\draw[thick] (\l+10/4,10/4) -- (\l+10/4.5,10/4);
\draw[thick] (\l+10/4,10/4) -- (\l+10/4,70/24);
\draw[thick] (\l+10/5,10/5) -- (\l+10/5.5,10/5);
\draw[thick] (\l+10/5,10/5) -- (\l+10/5,9/4);
\draw[thick] (\l+10/6,10/6) -- (\l+10/6.5,10/6);
\draw[thick] (\l+10/6,10/6) -- (\l+10/6,11/6);
\draw[thick] (\l+10/7,10/7) -- (\l+10/7.5,10/7);
\draw[thick] (\l+10/7,10/7) -- (\l+10/7,130/84);
\draw[dotted] (\l+10/8,10/8) -- (\l,0);
\draw[thick, dotted](\l+4.2265,5) to [bend right=20] (\l+10/3,50/12);
\draw[thick, dotted] (\l + 10/3.5,10/3) to [bend right=20] (\l+10/4,70/24);
\draw[thick, dotted] (\l+10/4.5,10/4)  to [bend right=20] (\l+10/5,9/4);
\draw[thick, dotted] (\l+10/5.5,10/5)  to [bend right=20] (\l+10/6,11/6);
\draw[thick, dotted] (\l+10/6.5,10/6)  to [bend right=20] (\l+10/7,130/84);

\end{tikzpicture}

\caption{Graphs of $\psi(x,x)$ and $\phi(x,x)$} \label{fig:psi&phi}
\end{figure}

For comparison, in figure \ref{fig:psi&phi} we give graphs of the function $\psi(x,x)$ (which we know completely, because of \ref{bisym3}),
and the function $\phi(x,x)$ (which we only know partially, from \ref{limk}, \ref{5-12} and \ref{intk}.)

The next result is a useful general lower bound on $\phi(x,y)$.

\begin{thm}\label{hompe4}
For $x,y,z\in (0,1]$, if $y> 1/2$ and $4x^2y(1-z)\ge (z-x)^2$ then $\phi(x,y)\ge z$. If in addition $4x^2y(1-z)> (z-x)^2$ then $\phi(x,y)> z$.
\end{thm}
\Proof
Let $G$ be $(x,y)$-constrained via $(A,B,C)$, and suppose that $|N^2_A(w)|< z|A|$ for each $w\in C$.
There are at least $xy|A|\cdot |B|\cdot |C|$ two-edge paths between $A$ and $C$, and so there is a vertex
$w\in C$ that is an end of at least $xy|A|\cdot |B|$ such paths. Let $w$ be an end of exactly $xq|A|\cdot |B|$ such paths; thus $y\le q$.
Let $B_1=N_B(w)$, and let $t=|B_1|/|B|$. Since $|N^2_A(w)|\le z|A|$,  
there
exists $A_1\subseteq A$ including $N^2_A(w)$ with $|A_1|=z|A|$ (we may assume the latter is an integer.) For each $u\in A_1$
let $u$ have exactly $d(u)|B|$ neighbours in $B_1$, and therefore at least $(x-d(u))|B|$ neighbours in $B\setminus B_1$.
It follows that
$$\sum_{u\in A_1}d(u) = qx|A|.$$

Let $v_1\in B_1$ and $v_2\in B\setminus B_1$, and let $A(v_1,v_2)=N_A(v_1)\cup N_A(v_2)$. For every such choice of $v_1,v_2$, since
$y>1/2$, there is a vertex $w$ in $C$ adjacent to both $v_1,v_2$, and since $|N^2_A(w)|<z|A|$, it follows that  
$|A(v_1,v_2)|<z|A|$. Let us choose $v_1\in B_1$ and $v_2\in B\setminus B_1$, uniformly at random.
It follows that the expected value of $|A(v_1,v_2)|$ is less than $z|A|$. The expected value of $|A(v_1,v_2)\cap A_1|$
is at least
$$\sum_{u\in A_1}\left(\frac{d(u)}{t} + \frac{x-d(u)}{1-t}-\frac{d(u)(x-d(u))}{t(1-t)}\right)$$
and the expected value of $|A(v_1,v_2)\setminus A_1|$ is at least
$$\sum_{u\in A\setminus A_1}\frac{x}{1-t}.$$
Consequently the sum of these two is less than $z|A|$, and so  
$$\sum_{u\in A_1}\left(\frac{d(u)}{t} + \frac{x-d(u)}{1-t}-\frac{d(u)(x-d(u)}{t(1-t)}\right) + \sum_{u\in A\setminus A_1}\frac{x}{1-t}<z|A|.$$
Since $\sum_{u\in A_1}d(v) = xq|A|$, this simplifies to
$$xq|A|(1-2t -x)+\sum_{u\in A_1}d(u)^2 + xt|A|<zt(1-t)|A|.$$
Now since $\sum_{u\in A_1}d(v) = xq|A|$ and $|A_1|=z|A|$, it follows from the Cauchy-Schwarz inequality that $\sum_{u\in A_1}d(u)^2\ge x^2q^2|A|/z$.
Consequently
$$xq|A|(1-2t -x)+x^2q^2|A|/z + xt|A|<zt(1-t)|A|.$$
This can be rewritten as:
$$(zt-xq+x/2-z/2)^2 + x^2q(1-z)- (z-x)^2/4< 0.$$
Since the first term above is a square, it is nonnegative, and so, since $q\ge y$, it follows that
$$x^2y(1-z)- (z-x)^2/4< 0,$$      
contrary to the hypothesis.
This proves the first statement of the theorem, and the second is immediate by slightly increasing $z$. This proves \ref{hompe4}.~\bbox

\section{When is $\phi(x,y)$ or $\psi(x,y)\ge 1/2$?}

Another way to approach the problem is to ask, given some value $z$, for which $x,y\in (0,1]$ is $\phi(x,y)\ge z$? Or we could ask the
same question for $\psi$, or ask when $\phi(x,y)>z$. For instance:
\begin{thm}\label{levelk+}
If $k\ge 1$ is an integer, then for $x,y\in (0,1]$, $\phi(x,y)>1/k$ if and only if $\max(x,y)>1/k$.
\end{thm}
This follows trivially from \ref{cayleybound} and \ref{maxbound}. And the same holds with $\phi$ replaced by $\psi$.
But deciding when $\psi(x,y)\ge 1/k$ or $\phi(x,y)\ge 1/k$ seems to be much less obvious. In this section
we discuss when $\psi(x,y)$ or $\phi(x,y)$ is at least $1/2$; and in later sections we look at when they are at least $2/3$,
and at least $1/3$.

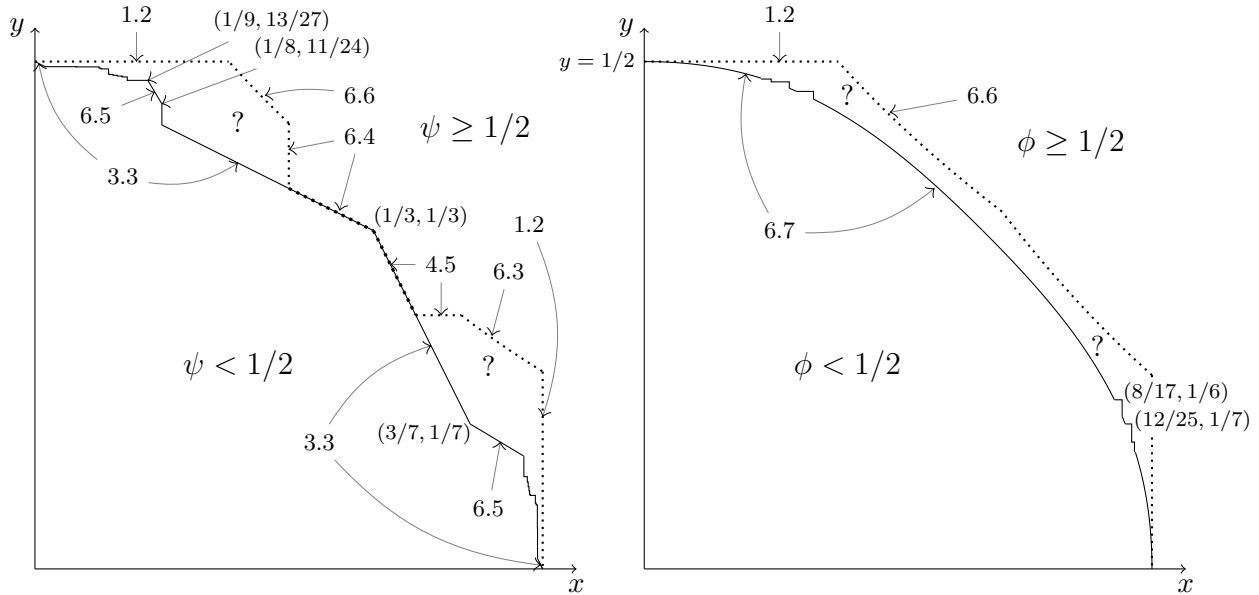
\begin{figure}[ht]
\centering

\begin{tikzpicture}[scale=.9,auto=left]

\begin{scope}[shift ={(-9,0)}]
\draw[->] (0,0) -- (8,0) node[anchor=north]{$x$};
\draw[->] (0,0) -- (0,8) node[anchor=east] {$y$};
\node at (5.7,5.2) {\begin{scriptsize}$(1/3,1/3)$\end{scriptsize}};
\node (n1) at (3.5,8.1) {\begin{scriptsize}$(1/9, 13/27)$\end{scriptsize}};

\draw[
    gray, ultra thin,decoration={markings,mark=at position 1 with {\arrow[black,scale=2]{>}}},
    postaction={decorate},
    ]
(n1) to (15/9,15*13/27);
\node (n2) at (4.1,7.7) {\begin{scriptsize}$(1/8,11/24)$\end{scriptsize}};
\draw[
    gray, ultra thin,decoration={markings,mark=at position 1 with {\arrow[black,scale=2]{>}}},
    postaction={decorate},
    ]
(n2) to (15/8,15*11/24);

\node at (5.75,2) {\begin{scriptsize}$(3/7,1/7)$\end{scriptsize}};

\node at (3,3) {\large{$\psi<1/2$}};
\node at (6.5,6.5) { \large{$\psi\ge 1/2$}};
\node at (3.0,6.6) {$?$};
\node at (6.7,3) {$?$};
\draw[domain=0:.19108,smooth,variable=\x,thick, dotted] plot ({15*\x},{15/2});
\draw[thick, dotted] (15/4, 15*.43934) -- (15/4,15*3/8);
%\draw[domain=.25:.2808,smooth,variable=\x,thick, dotted] plot ({15*\x},{15*(1-3*\x)/(1-2*\x)});
\draw[domain=1/4:1/3,smooth,variable=\x, very thick, dotted] plot ({15*\x},{7.5-15*\x/2});
\draw[domain=1/3:3/8,smooth,variable=\x, very thick, dotted] plot ({15*\x},{15*(1-2*\x)});
\draw[domain=3/8:.41991,smooth,variable=\x,thick, dotted] plot ({15*\x},{15*1/4});
\draw[domain=1/4:.1945, smooth, variable=\y, thick, dotted] plot ({15*(2+\y-4*\y*\y -\y*sqrt(9-16*\y-8*\y*\y))/(2+6*\y)},{15*\y});
%\draw[domain=1/4: .2113, smooth, variable=\y, thick, dotted] plot ({15*(1-3*\y+3*\y*\y)},{15*\y});
\draw[thick, dotted] (7.5,15*.1945) -- (7.5,0);

%\draw[domain=.19098:0.3522, smooth,variable=\y,thick, dotted] plot ({15*(1-\y)/(1+sqrt(2*\y))},{15*\y});
\draw[domain=1/4:.19098, smooth,variable=\x,thick, dotted] plot ({15*\x},{15*(1-\x)/(1+sqrt(2*\x))});

\draw[domain=1/3:3/7,smooth,variable=\x] plot ({15*\x},{15*(1-2*\x)});
\draw[domain=67/135:1/2,smooth,variable=\x] plot ({15*\x},{15*(1-2*\x)});
\draw[domain=3/7:13/27,smooth,variable=\x] plot ({15*\x},{15*(2-3*\x)/5});

\draw (15/2,0)--
(15*49/99, 15*1/99) -- (15*49/99, 15*4/99) --
(15*48/97, 15*4/99) -- (15*48/97, 15*4/97) --
(15*47/95, 15*4/97) -- (15*47/95, 15*6/95) --
(15*38/77, 15*6/95) -- (15*38/77, 15*5/77) --
(15*34/69, 15*5/77) -- (15*34/69, 15*5/69) --
(15*20/41, 15*5/69) -- (15*20/41, 15*3/41) --
(15*19/39, 15*3/41) -- (15*19/39, 15*1/13) --
(15*18/37, 15*1/13) -- (15*18/37, 15*3/37) --
(15*17/35, 15*3/37) -- (15*17/35, 15*3/35) --
(15*16/33, 15*3/35) -- (15*16/33, 15*1/11) -- (15*13/27,15*1/11)-- (15*13/27,15*1/9);

\draw 
(15*1/3,15*1/3) -- (15*1/8, 15*7/16) -- (15*1/8, 15*11/24) -- (15*1/9, 15*13/27) -- (15*1/11, 15*13/27)
-- (15*1/11, 15*16/33) -- (15*3/35, 15*16/33)
-- (15*3/35, 15*17/35) -- (15*3/37, 15*17/35)
-- (15*3/37, 15*18/37) -- (15*1/13, 15*18/37)
-- (15*1/13, 15*19/39) -- (15*3/41, 15*19/39)
-- (15*3/41, 15*20/41) -- (15*5/69, 15*20/41)
-- (15*5/69, 15*34/69) -- (15*5/77, 15*34/69)
-- (15*5/77, 15*38/77) -- (15*6/95, 15*38/77)
-- (15*6/95, 15*47/95) -- (15*4/97, 15*47/95)
-- (15*4/97, 15*48/97) -- (15*4/99, 15*48/97)
-- (15*4/99, 15*49/99) -- (15*1/99, 15*49/99)--(0,15/2);

\node (r1) at (1.5, 8.2) {\begin{footnotesize}$\ref{maxbound}$\end{footnotesize}};
\draw[
    gray, ultra thin, decoration={markings,mark=at position 1 with {\arrow[scale=2, black]{>}}},
    postaction={decorate},
    ]
(r1) to (1.5,7.5);

\node (r2) at (4.8,6.4) {\begin{footnotesize}$\ref{23topsi12}$\end{footnotesize}};
\draw[
    gray, ultra thin,decoration={markings,mark=at position 1 with {\arrow[black,scale=2]{>}}},
    postaction={decorate},
    ]
(r2) to ({4.5},{5.3});
\draw[
     gray, ultra thin,decoration={markings,mark=at position 1 with {\arrow[black,scale=2]{>}}},
    postaction={decorate},
    ]
(r2) to ({15/4},{6.2});

\node (r3) at (6,4.5) {\begin{footnotesize}$\ref{semibi}$\end{footnotesize}};
\draw[
    gray, ultra thin, decoration={markings,mark=at position 1 with {\arrow[black,scale=2]{>}}},
    postaction={decorate},
    ]
(r3) to ({15*.35},{15*.3});
\draw[
    gray, ultra thin,decoration={markings,mark=at position 1 with {\arrow[black,scale=2]{>}}},
    postaction={decorate},
    ]
(r3) -- ({15*.4},{15*.25});

\node (r4) at (7,4.4) {\begin{footnotesize}$\ref{bite}$\end{footnotesize}};
\draw[
    gray, ultra thin,decoration={markings,mark=at position 1 with {\arrow[black,scale=2]{>}}},
    postaction={decorate},
    ]
(r4) to ({15*(.45)},{15*.23});

\node (r5) at (7.3,5.1) {\begin{footnotesize}$\ref{maxbound}$\end{footnotesize}};
\draw[
    gray, ultra thin,decoration={markings,mark=at position 1 with {\arrow[black,scale=2]{>}}},
    postaction={decorate},
    ]
(r5) to [bend left = 20] ({7.5},{15*.15});

\node (r6) at (4.8,7) {\begin{footnotesize}$\ref{improvedphi12}$\end{footnotesize}};
\draw[
    gray, ultra thin,decoration={markings,mark=at position 1 with {\arrow[black,scale=2]{>}}},
    postaction={decorate},
    ]
(r6) to [bend right=10] ({3.4},{7});

\node (s1) at (.9,6.7) {\begin{footnotesize}$\ref{psi12extracurve}$\end{footnotesize}};
\draw[
    gray, ultra thin,decoration={markings,mark=at position 1 with {\arrow[black,scale=2]{>}}},
    postaction={decorate},
    ]
(s1) to ({7.5*(1/9+1/8)},{7.5*(13/27+11/24)});

\node (s2) at (1.3,5.8) {\begin{footnotesize}$\ref{psi12curve}$\end{footnotesize}};
\draw[
    gray, ultra thin,decoration={markings,mark=at position 1 with {\arrow[black,scale=2]{>}}},
    postaction={decorate},
    ]
(s2) to [bend right=20] ({15/5},{15*2/5});
\draw[
    gray, ultra thin,decoration={markings,mark=at position 1 with {\arrow[black,scale=2]{>}}},
    postaction={decorate},
    ]
(s2) to [bend left = 20] ({15*1/270},{15*269/540});

\node (s3) at (4.2,15/8) {\begin{footnotesize}$\ref{psi12curve}$\end{footnotesize}};
\draw[
    gray, ultra thin,decoration={markings,mark=at position 1 with {\arrow[black,scale=2]{>}}},
    postaction={decorate},
    ]
(s3) to [bend left = 20] ({15*(1-.22)/2},{15*.22});

\draw[
    gray, ultra thin,decoration={markings,mark=at position 1 with {\arrow[black,scale=2]{>}}},
    postaction={decorate},
    ]
(s3) to [bend right = 20] ({15*269/540},{15*1/270});

\node (s4) at (6.7,.9) {\begin{footnotesize}$\ref{psi12extracurve}$\end{footnotesize}};
\draw[
    gray, ultra thin,decoration={markings,mark=at position 1 with {\arrow[scale=2, black]{>}}},
    postaction={decorate},
    ]
(s4) to ({15*11/24},{15/8});

\end{scope}
%%%%%%%%%%%%%%%%%%%%%%%%%%%%%%%%%%%%%%%%%%%%%%%%%%%%%%%%%%%%%%%%%%%%%%%%%%%%%%%%%%%%%%%%%%%%%%%%%%%%%%%%%%%%%%%%%%%%%%%%%%%
\draw[->] (0,0) -- (8,0) node[anchor=north]{$x$};
\draw[->] (0,0) -- (0,8) node[anchor=east] {$y$};
\draw[thick,dotted] (0,15/2)--(15*.19098,15/2);
\draw[thick, dotted] (15/2,15*.19098) -- (15/2,15*1/6+.25);
\draw[thick,dotted] (15/2,15*1/7-.2) -- (15/2,0);

\draw (15*14/29,15*0.115) -- (15*14/29, 15*1/8) --
(15*12/25, 15*1/8) -- (15*12/25,15*1/7) -- (15*36/76, 15*1/7);
\draw (15*8/17, 15*.15022) -- (15*8/17, 15*1/6) -- (15*25/54,15*1/6);

\draw
(15*1/6, 15*25/54) --
(15*1/6, 15*8/17) -- (15*.15022, 15*8/17);
\draw (15*1/7, 15*36/76) -- (15*1/7, 15*12/25) -- (15*1/8, 15*12/25)
-- (15*1/8, 15*14/29) -- (15*0.115, 15*14/29);

\node at (15*12/25+0.9,15/7+.05) {\begin{scriptsize}$(12/25,1/7)$\end{scriptsize}};
\node at (15*8/17+.8,15*1/6+.1) {\begin{scriptsize}$(8/17,1/6)$\end{scriptsize}};
\node at (6.3,6.3) { \large{$\phi\ge 1/2$}};
\node at (3,3) { \large{$\phi< 1/2$}};
\node at (6.7,3.3) {$?$};
\node at (3.0,7) {$?$};
\draw[domain=.19098:0.3522, smooth,variable=\y,thick, dotted] plot ({15*(1-\y)/(1+sqrt(2*\y))},{15*\y});
\draw[domain=.3522:.19098, smooth,variable=\x,thick, dotted] plot ({15*\x},{15*(1-\x)/(1+sqrt(2*\x))});

\draw[domain=0:0.115,smooth,variable=\x] plot ({15*\x},{15*(1-\x)*(1-\x)/(2-4*\x+6*\x*\x)});
\draw[domain=1/7:.15022,smooth,variable=\x] plot ({15*\x},{15*(1-\x)*(1-\x)/(2-4*\x+6*\x*\x)});
\draw[domain=1/6:1/3,smooth,variable=\x] plot ({15*\x},{15*(1-\x)*(1-\x)/(2-4*\x+6*\x*\x)});
\draw[domain=0:0.115,smooth,variable=\y] plot ({15*(1-\y)*(1-\y)/(2-4*\y+6*\y*\y)},{15*\y});
\draw[domain=1/7:.15022,smooth,variable=\y] plot ({15*(1-\y)*(1-\y)/(2-4*\y+6*\y*\y)},{15*\y});
\draw[domain=1/6:1/3,smooth,variable=\y] plot ({15*(1-\y)*(1-\y)/(2-4*\y+6*\y*\y)},{15*\y});

\node at (-.7, 7.5) {\begin{scriptsize}$y=1/2$\end{scriptsize}};
\node (t1) at (2, 8.2) {\begin{footnotesize}$\ref{maxbound}$\end{footnotesize}};
\draw[
    gray, ultra thin,decoration={markings,mark=at position 1 with {\arrow[black,scale=2]{>}}},
    postaction={decorate},
    ]
(t1) to (2,7.5);
\node (t2) at (5, 7) {\begin{footnotesize}$\ref{improvedphi12}$\end{footnotesize}};
\draw[
    gray, ultra thin,decoration={markings,mark=at position 1 with {\arrow[black,scale=2]{>}}},
    postaction={decorate},
    ]
(t2) to ({15*.242},{15*.45});

\node (u1) at (2,5) {\begin{footnotesize}$\ref{phi12bettercurve}$\end{footnotesize}};
\draw[
    gray, ultra thin,decoration={markings,mark=at position 1 with {\arrow[black,scale=2]{>}}},
    postaction={decorate},
    ]
(u1) to [bend left=20] (15*1/10, 15*81/166);
\draw[
    gray, ultra thin,decoration={markings,mark=at position 1 with {\arrow[black,scale=2]{>}}},
    postaction={decorate},
    ]
(u1) to[bend right = 20]  (15*2/7,15*3/8);

\end{tikzpicture}

\caption{In the left-hand figure, $\psi(x,y)<1/2$ for pairs $(x,y)$ below the solid line, and $\psi(x,y)\ge 1/2$ 
above the dotted one;
between we don't know. The right-hand figure does the same for $\phi$.} \label{fig:phi12}
\end{figure}

For $x,y\in (0,1]$, we say (temporarily) that
$(x,y)$ is {\em good} if $\psi(x,y)\ge 1/2$, and {\em bad} otherwise.
The ``map'' of good and bad points is shown in the left half of figure~\ref{fig:phi12}. The solid black curve borders
the known bad points, and the dotted curve borders the good points; between them is undecided. The borders are complicated, and we
have indicated in the figure which theorem is responsible for each stretch of border.

Let us explain some of the details. First, if $\max(x,y)\ge 1/2$, then
$(x,y)$ is good; and all pairs $(x,y)$ with
$x+2y,2x+y\le 1$ are bad, by \ref{psi12curve}.
We searched by computer to find other examples of bad pairs $(x,y)$,
and found about 12 maximal such pairs of rationals, with numerator and denominator at most 100.
In fact we only searched for pairs $(x,y)$ where the corresponding $(x,y)$-biconstrained graph
is similar to the graph obtained from figure \ref{fig:exactcount}, that is, it is obtained by ``blowing up''
the vertices of another graph in which the graph between two of the three parts is a matching. All these examples
not only show that $\psi(x,y)<1/2$, but also that $\psi(y,x)<1/2$, and $\xi(x,y)<1/2$. In particular,
for every bad pair $(x,y)$ we found by computer search, $(y,x)$ is another.
This is just an artifact of our method of search, and is not evidence
that the set of all bad pairs is closed under switching $x$ and $y$ (though it might be; it is for $\phi$, by \ref{permute}).
Anyway, for each bad pair the computer found, all pairs it dominates are also bad, and that gave us a step function
bordering the area of the known bad points. We improved on this; we were able to smooth out some of the
steps of the step function, by means of
\ref{psi12curve} and  \ref{psi12extracurve}, so the step function the computer found now only survives towards the ends of the solid
black curve in the figure. (These ``fills'' are not invariant under switching $x$ and $y$.)
We give the coordinates of some bad pairs that we find particularly interesting.
The apparent asymmetry between $x$ and $y$ in the left half of the figure is just asymmetry among what we have been able to prove;
we have no proof of asymmetry. The right half of figure~\ref{fig:phi12} does the same for $\phi$. Here
there is symmetry exchanging $x$ and $y$, by \ref{permute}, and so
we only ``explain'' half of the border.

A graph has {\em radius} at most $r$ if there is a vertex $u$ such that every vertex has distance at most $r$ from $u$.
We will need the following theorem of Erd\H{o}s, Saks and S\'os~\cite{radius}:
\begin{thm}\label{radius}
Let $G$ be a connected graph with radius at least $r$, where $r\ge 1$ is an integer. Then $G$ has an induced path with $2r-1$ vertices,
and consequently has a stable set of cardinality at least $r$.
\end{thm}
When we have more than one graph defined using the same vertices, we speak of ``$H$-distance'' to mean distance in the graph $H$,
and so on.

\begin{thm}\label{bite}
Let $x,y\in (0,1]$, such that 
$$x^2(1+3y)+x(4y^2-y-2)+1-2y+2y^3<0.$$
Then $\psi(x,y)\ge 1/2$.
\end{thm}
%and $x/(1-y)\ge 3-6x$ and $y/(4(1-y))+x\ge 1/2$ and $3y>1-x$ and $y/(6(1-2y))>1/2-x$
%and let $s$ be real with $1\le s\le 2$, such that $x+(s+1)y>1$ and $x/(1-y)\ge 3-6x$ and $2sy(2x+y-1)>(1-x-2y)(s-1+x)$ and $y/(4(1-y))+x\ge 1/2$.  Then $\psi(x,y)\ge 1/2$.  \end{thm}
\Proof
Let $G$ be $(x,y)$-biconstrained via $(A,B,C)$, and suppose that $|N^2_A(v)|<|A|/2$ for each $v\in C$.
Then \ref{maxbound} implies that $x,y<1/2$.
Suppose that $y>1/4$. The given inequality implies that $56x^2-64x+17<0$, and so $x>.41$. 
Since $2x+y> 1$, \ref{semibi} implies that $y\le 1/4$, a contradiction. Thus $y\le 1/4$.
We leave the reader to verify that the following are consequences of the given inequality:
\begin{itemize}
\item $\frac{x}{3-3y}> 1-2x$; and in particular, $\frac{x}{1-y}> 1-2x$, so from \ref{bisym} it follows that $x+2y\le 1$;
\item $\frac{y}{3-6y}>1-2x$; and so $\frac{y}{2-2y}> 1-2x$, since $y\le 1/4$; and
\item $x+3y>1$.
\end{itemize}
(We found the easiest way to check these is to have a computer plot the various curves.)
Let $H$ be the bipartite graph $G[B\cup C]$.
\\
\\
(1) {\em If $v,v'\in C$ have $H$-distance at most $2t$ where $t>0$ is an integer, then 
$$|N^2_A(v')\setminus N^2_A(v)|<t(1/2-x)|A|.$$}
Take a path $P$ of $H$ joining $v$ and $v'$, of length at most $2t$.
Let the vertices of $P$ in $C$
be 
$$v=v_0\ll v_t=v',$$ in order. For $1\le i\le t$ let $u_i\in B$ be adjacent to $v_{i-1}$ and $v_i$.
Then for $1\le i\le t$,  $N^2_A(v_{i-1})\cap N^2_A(v_i)$ includes $N_A(u_i)$ and hence has cardinality at least $x|A|$;
and since 
$|N^2_A(v_i)|<|A|/2$, it follows that $|N^2_A(v_{i})\setminus N^2_A(v_{i-1})|<(1/2-x)|A|$. But the union of the $t$
sets $N^2_A(v_{i})\setminus N^2_A(v_{i-1})$ includes $N^2_A(v')\setminus N^2_A(v)$, and so the latter has cardinality less than 
$t(1/2-x)|A|$. This proves (1).
\\
\\
(2) {\em There do not exist $v_1\ll v_4\in C$, pairwise with no common neighbour in $B$.}
\\
\\
Then every three of $N(v_1)\ll N(v_4)$ have union of cardinality at least $3y$; and since $3y>1-x$, 
every vertex in $A$ has a neighbour in at least
two of $N(v_1)\ll N(v_4)$. Consequently every vertex in $A$ belongs to at least two of $N^2_A(v_1)\ll N^2_A(v_4)$, and so
one of $N^2_A(v_1)\ll N^2_A(v_4)$ has cardinality at least $|A|/2$, a contradiction. This proves (2).
\\
\\
(3) {\em $H$ has at least two components.}
\\
\\
Suppose not, and let $H'$ be the graph with vertex set $C$ in which $v,v'$ are adjacent if they have a common neighbour
in $H$. By (2), it follows that $H'$ has no stable set of cardinality four, and so has radius at most three by \ref{radius}.
Choose $v\in C$ such that every vertex in $C$ has $H'$-distance at most three from $v$. Let $B_1=N_B(v)$ and $A_1=N^2_A(v)$.
Every vertex in $A\setminus A_1$ has at least $x|B|$ neighbours in $B\setminus B_1$, and so some vertex $u\in B\setminus B_1$
has at least $x(|B|/|B\setminus B_1|)|A\setminus A_1|$ neighbours in $A\setminus A_1$. Let $A_2$ be the set of neighbours of $u$
in $A\setminus A_1$. Since $|B\setminus B_1|\le (1-y)|B|$ and $|A\setminus A_1|>|A|/2$, and $\frac{x}{3-3y}> 1-2x$,
it follows that
$$|A_2|\ge (x/(1-y))|A|/2\ge 3(1/2-x)|A|.$$ 
Let $v'\in C$ be adjacent to $u$.
Since the $H'$-distance from $v$ to $v'$ is at most three, the $H$-distance from $v$ to $v'$
is at most six. By (1), $|N^2_A(v')\setminus N^2_A(v)|<3(1/2-x)|A|$, a contradiction.
This proves (3).
\\
\\
(4) {\em If $H'$ is a component of $H$ then $|V(H')\cap B|\le (1-x)|B|$.}
\\
\\
Suppose that $|V(H')\cap B|>(1-x)|B|$; then every vertex in $A$ has a neighbour in $V(H')$. By (2), there do not exist three vertices
in $C\cap V(H')$ pairwise with no common neighbour, and so by \ref{radius}, it follows that there is a vertex $v\in C\cap V(H')$
with $H'$-distance at most four from every vertex in $C\cap V(H')$. Let
$A'=N^2_A(v)$; then $|A'|< |A|/2$. Since every vertex in $A\setminus A'$ has at least $y|C|$ second neighbours in $C\cap V(H')$,
and $|C\cap V(H')|\le (1-y)|C|$, some vertex $v'\in C\cap V(H')$ has at least $(y/(1-y))|A\setminus A'|$ second neighbours
in $A\setminus A'$. By (1),
$(y/(1-y))|A\setminus A'|<2(1/2-x)|A|$. But $|A'|<|A|/2$, so $y/(4(1-y)) \le 1/2-x$, a contradiction.
This
proves (4).
\\
\\
(5) {\em Some component $H'$ of $H$ satisfies $(1-x)|B|\ge |V(H')\cap B|\ge x|B|$.}
\\
\\
By (2) and (3), $H$ has either two or three components. If $H$ has only two components, then they both satisfy (5), by (4); so we assume there are three.
Let the components of $H$ be $H_1,H_2,H_3$, and for $1\le i\le 3$, let
$V(H_i)\cap B=B_i$ and $V(H_i)\cap C=C_i$;
and let $|B_i|/|B|=b_i$ and $|C_i|/|C|=c_i$. Suppose that $b_1,b_2,b_3< x$.
Consequently every vertex in $A$
has neighbours in at least two of $B_1,B_2,B_3$. For $1\le i\le 3$, let $A_{i}$ be the set of vertices in $A$
with a neighbour in $B_i$. Thus every vertex in $A$ belongs to at least two of $A_1,A_2,A_3$, so from the symmetry
we may assume that $|A_1|\ge 2|A|/3$.
By (2), every two vertices in $C_1$
have a common neighbour in $B$. Choose $v\in C_1$, and let
$A'=N^2_A(v)$; then $|A'|\le |A|/2$. Since every vertex in $A_1$ has at least $y|C|$ second neighbours in $C_1$,
some vertex $v'$ in $C_1$ has at least $(y/c_1)|A_1\setminus A'|$ second neighbours
in $A_1\setminus A'$. By (1),
$(y/c_1)|A_1\setminus A'|<(1/2-x)|A|$. But $|A'|<|A|/2$, so $|A_1\setminus A'|\ge |A|/6$; and $c_1\le 1-2y$,
so $y/(6(1-2y)) \le 1/2-x$, a contradiction. This proves (5).
\\
\\
(6) {\em Every vertex in $C$ has at most $(1-x-y)|B|$ neighbours in $B$. Consequently
$$\frac{|V(H')\cap B|}{|B|}\le \frac{1-x-y}{y}\frac{|V(H')\cap C|}{|C|}$$ 
for each component $H'$ of $H$.}
\\
\\
Suppose that $v\in C$ has more than $(1-x-y)|B|$ neighbours in $B$. Choose $v'\in C$ in a different component of $H$; so
$v,v'$ have no common neighbour in $B$. Consequently
$$|N(v)\cup N(v')|> ((1-x-y)+y))|B|,$$
and so every vertex in $A$
has a neighbour in $N(v)\cup N(v')$. But then one of $|N^2_A(v)|,|N^2_A(v')|\ge |A|/2$, a contradiction. This proves
the first assertion. Let $H'$ be a component of $H$. Then $H'$ has at least $y|C|\cdot |V(H')\cap B|$ edges, and at most
$(1-x-y)|B|\cdot |V(H')\cap C|$ edges, so the second claim follows. This proves (6).

\bigskip

Let $H'$ be as in (5), and take the union of the other (one or two) components of $H'$. We obtain
nonnull subgraphs $H_1,H_2$ of $H$, pairwise vertex-disjoint and with union $H$, such that 
$|V(H_i)\cap B|\ge x|B|$ for $i = 1,2$.
For $i = 1,2$, let
$V(H_i)\cap B=B_i$ and $V(H_i)\cap C=C_i$;
and let $|B_i|/|B|=b_i$ and $|C_i|/|C|=c_i$. Thus $b_1,b_2\ge x$. From (6), $b_i\le (1-x-y)c_i/y$ for $i = 1,2$; and
$c_1,c_2\ge y$, since every vertex
in $B_i$ has at least $y|C|$ neighbours in $C_i$. Also $b_1+b_2=c_1+c_2=1$.

For $i = 1,2$ let $A_i$ be the set of vertices $v\in A$ that have more than $(b_i-y)|B|$ neighbours in $B_i$.
Let $A_0 = A\setminus (A_1\cup A_2)$. Hence if $u\in A_0$, then since $u$ has at least $x|B|$ neighbours in $B$,
$u$ has at least $(x+y-b_2)|B|$ neighbours in $B_1$, and at least $(x+y-b_1)|B|$ neighbours in $B_2$.

Since $A_1, A_2$ and $A_0$ have union $A$, we may assume that $|A_1|+|A_0|/2\ge |A|/2$.
Now $A_1\subseteq N^2_A(v)$ for each $v\in C_1$, since if $u\in A_1$, then $u$
has more than $(b_1-y)|B|$ neighbours in $B_1$, and $v$ has at least $y|B|$ neighbours in $B$. Consequently
$|N^2_A(v)\cap A_0|<|A_0|/2$ for each $v\in C_1$.

Let us choose $v\in C_1$ uniformly at random; then the expected
number of second neighbours of $v$ in $A_0$ is less than $|A_0|/2$, and so for some vertex $u\in A_0$, the probability
that $u\in N^2_A(v)$ is less than $1/2$. Let $D$ be the set of neighbours of $u$ in $B_1$.
Then $|D|\ge (x+y-b_2)|B|$,
and the probability that $v$ has a neighbour in $D$ is less than $1/2$. Thus more than $|C_1|/2$ vertices in $C_1$
have no neighbour in $D$. On the other hand, the expectation of the number of neighbours of $v$ in $D$ is at least
$|D|y/c_1$; and so there exists $v\in C_1$ with more than $2|D|y/c_1$ neighbours in $D$. Also there exists $v'\in C_1$
with no neighbours in $D$. It follows that
$$|N_B(v)\cup N_B(v')|\ge y|B|+2|D|y/c_1> (y+2(x+y-b_2)y/c_1)|B|.$$
Some vertex in $A_0$ is not a second neighbour of either of $v, v'$, and so
$$|N_B(v)\cup N_B(v')|< (b_1-(x+y-b_2))|B|.$$
Consequently
$y+2(x+y+b_1-1)y/c_1\le 1-x-y$. Now $c_1\le (1-x-2y+yb_1)/(1-x-y)$ since 
$$1-b_1=b_2\le (1-x-y)c_2/y= (1-x-y)(1-c_1)/y$$. So
$$2y(x+y+b_1-1)(1-x-y)/(1-x-2y+yb_1)\le 1-x-2y,$$
that is, 
$$b_1y(1-x)\le x^2(1+2y) +x( -2 +4y^2) +1 -2y +2y^3.$$
But $b_1\ge x$, contrary to the hypothesis. This proves \ref{bite}.~\bbox

\begin{thm}\label{23topsi12}
Let $x,y\in (0,1]$, such that $x+2y>1$, $x\ge 1/4$ and $y\ge 1/3$. Then $\psi(x,y)\ge 1/2$.
\end{thm}
\Proof The only lower bound constraints on $y$ are $y\ge 1-2x$ and $y\ge 1/3$, and these are both satisfied if $y=0.38$ since $x\ge 1/4$. Hence we may 
assume that $y\le 0.38$, by replacing $y$ by $\min(y,0.38)$. Consequently $y^2-3y+1> 0$, and so 
$$(1-y)^3<1-2y<x.$$
Let $G$ be $(x,y)$-biconstrained via $(A,B,C)$, and suppose that $|N^2_A(w)|<|A|/2$ for each $w\in C$.
Choose $w_1,w_2,w_3\in C$ uniformly at random. The expected number of vertices in  $B$ nonadjacent to all of $w_1,w_2,w_3$
is at most $(1-y)^3|B|<x|B|$; so we may choose $w_1,w_2,w_3$ such that 
fewer than $x|B|$ vertices in $B$ are nonadjacent to all of $w_1,w_2,w_3$. For $i=1,2,3$ let $A_i=N^2_A(w_i)$. Thus $A_1\cup A_2\cup A_3=A$.
In particular, one of $A_2,A_3$, say $A_2$, includes at least half of $A\setminus A_1$; and since $|A_1|<|A|/2$, it follows that
$|A_2\setminus A_1|>|A|/4$. Since $|A_2|<|A|/2$, it follows that $|A_1\cap A_2|<|A|/4<x|A|$; and so $N_B(w_1), N_B(w_2)$ are disjoint (because 
any common neighbour would have at least $x|A|$ neighbours in $A$, all belonging to $A_1\cap A_2$). Hence 
$$|N_B(w_1)\cup N_B(w_2)|\ge 2y|B|>(1-x)|B|,$$
and so $A_1\cup A_2=A$, contradicting that $|A_1|,|A_2|<|A|/2$. This proves \ref{23topsi12}.~\bbox

\begin{thm}\label{psi12extracurve}
Let $x,y\in (0,1]$, such that $x\le 13/27$ and $y\le 1/7$ and $3x+5y\le 2$. Then $\psi(x,y)\le 13/27$. If in addition $y\le 1/8$,
then $\psi(y,x)< 1/2$.
\end{thm}
\Proof
We claim that, for both statements of the theorem, we may assume that $3x+5y=2$. 
By increasing $x$, we may assume that either $x=13/27$ or $3x+5y=2$; 
and if $x=13/27$ then $y\le 1/9$, since $3x+5y\le 2$,
and by increasing $y$ we may assume that $3x+5y=2$. This proves our claim. 
Since $3x+5y=2$ and $x\le 13/27$, it follows that $y\ge 1/9$; and since $y\le 1/7$ it follows that $x\ge 3/7$.

We return to the graph of figure \ref{fig:exactcount}. Let $A,B,C$ be the three rows of vertices, in order where $A$ is
the top row.
We need to adjust the vertex weights. Define $p=1/2-x/2-y$ and $r=(x-y)/2$.
With the vertices in the same order as the figure,
take vertex weights as follows:
\begin{eqnarray*}
5/27, 5/27,& 1/9, 1/9, 1/9,& 4/27, 4/27\\
p,\ \ \  \ \   p,\ \ & y,\ \ \  y,\ \ \  y,& \ \  r,\ \ \ \ \  r\\
1/7,\ \  1/7,\ & 1/7, 1/7, 1/7,& \ 1/7, \ 1/7
\end{eqnarray*}
One can check (it takes some time and we omit the details) that this defines an $(x,y)$-biconstrained weighted graph
showing that $\psi(x,y)\le 13/27$. For the second statement, take the same graph and same vertex weighting, except
replace the third row (of all one-sevenths) in the table above, by
\begin{eqnarray*}
p',\ \ \  p',\  & y,\  y,\  y,\ & r',\ \ \ r'
\end{eqnarray*}
where $p'=1/2-3y$, and $r'=3y/2$. This weighted graph is $(y,x)$-biconstrained via $(C,B,A)$, and shows
that $\psi(y,x)<1/2$.
(Again, we leave the reader to check that this works.)
This proves \ref{psi12extracurve}.~\bbox

Now the mono-constrained case: for which pairs $(x,y)$ is $\phi(x,y)\ge 1/2$? Now we have symmetry between $x$ and $y$,
and we found some examples of pairs $(x,y)$ with $\phi(x,y)<1/2$ on a computer searching randomly. 
(Conjecture~\ref{monosymconj}
says that all points above both the lines $x+2y=1$ and $2x+y=1$ should be good, and indeed, all the maximal examples 
the computer found lie
in the wedges between the lines.) 

The next result strengthens \ref{limk1} when $k=2$:

\begin{thm}\label{improvedphi12}
Let $x,y\in (0,1]$, such that $2x^2y \ge (1-x-y)^2$. Then $\phi(x,y)\ge 1/2$.
\end{thm}
\Proof
Suppose that $\phi(x,y)=1/2-\varepsilon$ where $\varepsilon>0$. Then by \ref{permute} and \ref{rotate} we have
$\phi(x,1/2+\varepsilon)\le 1-y$. Let $y'=1/2+\varepsilon$ and $z=1-y$. Since $y'>1/2$ and $\phi(x,y')\le z$, the 
second statement of \ref{hompe4} implies that $4x^2y'(1-z)\le (z-x)^2$, and so $2x^2y< (1-x-y)^2$, a contradiction. This 
proves~\ref{improvedphi12}.~\bbox
In particular, \ref{improvedphi12} implies that $\phi(x,x)\ge 1/2$ if $x\ge 0.352202$, which is stronger than \ref{limk1} when $k=2$.

\begin{thm}\label{phi12bettercurve}
Let $x,y\in (0,1]$ with $y< 1/2$ and $y< (1-x)^2/(2-4x+6x^2)$; then $\phi(x,y)<1/2$.
\end{thm}
\Proof By increasing $y$ slightly if necessary, we may assume that $s=(1/y-2)^{1/2}$ is rational.
Thus $s\le 1/y-1$, $s^2\le 1/y-2$ and
$1+2/s\le 1/x$, since $y\le 1/2$ and $y\le (1-x)^2/(2-4x+6x^2)$. 
Choose an integer $N\ge 2$ such that $sN$ is an integer.

Choose a graph $G_1$ that is $(s,1-s)$-constrained via a tripartition $(A_1,B_1,C_1)$,
such that $|A_1|=|B_1|=|C_1|=N$
and $N^2_{A_i}(v)\ne A_1$ for each $v\in C_1$. (It is easy to see that such a graph exists, for instance, one of the
graphs used in \ref{cayleybound}.) Let $G_2$ be isomorphic to $G_1$, and let $(A_2,B_2,C_3)$ be the corresponding 
tripartition. Take the disjoint union of $G_1$ and $G_2$, and add edges to make every vertex in $B_1$
adjacent to every vertex in $C_2$. Add three more vertices $a,b,c$,
where $a$ is adjacent to $b$, $b$ is adjacent to every vertex in $C_1$, and $c$ is adjacent to every vertex in $B_2$,
forming $G$.
We define a weighting $w$ of $G$ as follows. Let $p=1/(2N)$ and $q=1/2-1/(4N^2)$.
Define $w$ by:
\begin{eqnarray*}
w(a)&=& 1-p-q\\
w(v)&=&p/N \text{ for each }v\in A_1\\
w(v)&=& q/N \text{ for each }v\in A_2\\
w(b)&=&1-2x/s\\
w(v)&=& x/(Ns) \text{ for each }v\in B_1\cup B_2\\
w(c)&=& 1-(1+s)y\\
w(v)&=& y/N \text{ for each }v\in C_1\\
w(v)&=& sy/N \text{ for each }v\in C_2
\end{eqnarray*}
Define $A=A_1\cup A_2\cup\{a\}$ and define $B,C$ similarly.
Then the weighted graph $(G,w)$ is $(x,y)$-constrained via $(A,B,C)$, and proves that $\phi(x,y)<1/2$.
(To see the latter, observe that, for instance, if $v\in C_1$ then 
$$w(N^2_A(v))\le 1-p-q+ p(1-1/N)=1-q-p/N=1 -(1/2-1/(4N^2))-1/(2N^2)<1/2,$$
from the choice of $G_1$).
This proves \ref{phi12bettercurve}.~\bbox

\section{The $2/3$ level}
When is $\phi(x,y)\ge 2/3$; or the same question for $\psi$? In this section we
say what we know about these.

\begin{figure}[ht]
\centering

\begin{tikzpicture}[scale=.68,auto=left]

\begin{scope}[shift ={(-12,0)}]
\draw[->] (0,0) -- (11,0) node[anchor=north]{$x$};
\draw[->] (0,0) -- (0,11) node[anchor=east] {$y$};
\node at (8.4,7.7) {\begin{scriptsize}$(1/2,1/2)$\end{scriptsize}};
%\node at (15*12/25+1.1,15/7) {\begin{scriptsize}$(12/25,1/7)$\end{scriptsize}};
\node at (15*3/5+1,15/5+.1) {\begin{scriptsize}$(3/5,1/5)$\end{scriptsize}};
%\node at (15*8/17+1, 15*1/6+.1) {\begin{scriptsize}$(8/17,1/6)$\end{scriptsize}};
\node at (4.2,8.1) {$?$};
\node at (7.85, 3.95) {$?$};
\node at (9.5, 2.25) {$?$};
%Start of insert
\draw

(15*65/98, 15*0/1) -- (15*65/98, 15*1/98) -- 
(15*63/95, 15*1/98) -- (15*63/95, 15*1/95) -- 
(15*61/92, 15*1/95) -- (15*61/92, 15*1/92) -- 
(15*59/89, 15*1/92) -- (15*59/89, 15*1/89) -- 
(15*57/86, 15*1/89) -- (15*57/86, 15*1/86) -- 
(15*55/83, 15*1/86) -- (15*55/83, 15*1/83) -- 
(15*53/80, 15*1/83) -- (15*53/80, 15*1/80) -- 
(15*51/77, 15*1/80) -- (15*51/77, 15*1/77) -- 
(15*49/74, 15*1/77) -- (15*49/74, 15*1/74) -- 
(15*47/71, 15*1/74) -- (15*47/71, 15*1/71) -- 
(15*45/68, 15*1/71) -- (15*45/68, 15*1/68) -- 
(15*43/65, 15*1/68) -- (15*43/65, 15*1/65) -- 
(15*41/62, 15*1/65) -- (15*41/62, 15*1/62) -- 
(15*39/59, 15*1/62) -- (15*39/59, 15*1/59) -- 
(15*37/56, 15*1/59) -- (15*37/56, 15*1/56) -- 
(15*35/53, 15*1/56) -- (15*35/53, 15*1/53) -- 
(15*33/50, 15*1/53) -- (15*33/50, 15*1/50) -- 
(15*31/47, 15*1/50) -- (15*31/47, 15*2/47) -- 
(15*29/44, 15*2/47) -- (15*29/44, 15*1/22) -- 
(15*25/38, 15*1/22) -- (15*25/38, 15*1/19) -- 
(15*23/35, 15*1/19) -- (15*23/35, 15*2/35) -- 
(15*21/32, 15*2/35) -- (15*21/32, 15*1/16) -- 
(15*19/29, 15*1/16) -- (15*19/29, 15*2/29) -- 
(15*17/26, 15*2/29) -- (15*17/26, 15*1/13) -- 
(15*15/23, 15*1/13) -- (15*15/23, 15*2/23) -- 
(15*11/17, 15*2/23) -- (15*11/17, 15*2/17);

%end of insert

\draw (15*11/17, 15*2/17) -- (15*3/5,15*2/15) -- (15*3/5,15*1/5) -- 
%(15*3/5, 15*1/8) -- (15*3/5, 15*1/5) -- 
%(15*17/30, 15*1/5);
%\draw (15*17/30, 15*1/5) -- (15*5/9, 15*2/9) --
(15*1/2, 15*1/4) -- (15*1/2, 15*1/2)
-- (15*1/4, 15*1/2) -- (15*1/5, 15*3/5) ;
\draw
%(15*2/17,15*32/51) -- 
(15*1/8, 15*5/8) %this is \ref{psi23extra}
--(15*2/17, 15*11/17) -- (15*2/23, 15*11/17) 
-- (15*2/23, 15*15/23) -- (15*3/35, 15*15/23) 
-- (15*3/35, 15*23/35) -- (15*1/19, 15*23/35) 
-- (15*1/19, 15*25/38) -- (15*2/41, 15*25/38) 
-- (15*2/41, 15*27/41) -- (15*1/41,15*27/41);

\node at (9,9) { \large{$\psi\ge 2/3$}};
\node at (5.3,5.3) { \large{$\psi< 2/3$}};
\node at (5,7.2) {\begin{scriptsize}$(1/3,1/2)$\end{scriptsize}};
\node at (8.5,5) {\begin{scriptsize}$(1/2,1/3)$\end{scriptsize}};
\node at (2.8,7.2) {\begin{scriptsize}$(1/4,1/2)$\end{scriptsize}};
\node at (3.9,9.1) {\begin{scriptsize}$(1/5,3/5)$\end{scriptsize}};
\draw[thick, dotted] (15/6, 10) -- (0,10);
\draw[domain=5/9:.58579, smooth, variable=\y, thick,dotted] plot ({15*(2-3*\y)}, {15*\y});
\draw[domain=.58579:.6, smooth, variable=\y, thick, dotted] plot ({15*2*(1-\y)*(1-\y)/(2-\y)}, {15*\y});
\draw[domain=.62:2/3, smooth, variable=\y, thick, dotted] plot ({15*2*(1-\y)*(1-\y)/(2-\y)}, {15*\y});
%\draw[domain=.305:1/3, smooth, variable=\x, thick, dotted] plot ({15*\x}, {15*(1-\x)*(1-\x)/(1-2*\x*\x)});
%\draw[domain=1-1/sqrt(6):2/3,smooth,variable=\y,thick,dotted] plot ({15*2*(1-\y)*(1-\y)},{15*\y});
\draw[thick, dotted] (5,{15*5/9}) -- (5,7.5);
\draw[very thick, dotted](5,7.5) -- (7.5,7.5) -- (7.5,5);
\draw[domain=1/2-1/sqrt(12):1/3,smooth,variable=\y,thick,dotted] plot ({15*(2-\y/(1-\y))/3},{15*\y});
\draw[very thick, dotted] ({15/sqrt(3)}, {15*(1/2-1/sqrt(12))})-- (15*3/5,15*1/5);
\draw[thick, dotted] ({15*3/5}, {15*1/5})-- (10,2.5) -- (10,0);

%\draw[domain=4/7:11/17, smooth, variable=\x, very thick] plot ({15*\x}, {15*(1-\x)/3});
%\draw[domain=5/8:11/17, smooth, variable=\y, very thick] plot ({15*(1-\y)/3}, {15*\y});

%\draw[domain=0:1/12,smooth,variable=\x, thick] plot ({15*\x},{15*(3-5*\x)/4});
%\draw[domain=0:1/12,smooth,variable=\x, thick] plot ({15*\x},{15*(1-8*\x/3-20*\x*\x)});
%\draw[domain=0:1/12,smooth,variable=\x, thick] plot ({15*\x},{15/(4/3+5*\x/2)});
%\draw[domain=0:1/12,smooth,variable=\x, thick] plot ({15*\x},{15*27/41});

\draw[domain=3/5:5/8,smooth,variable=\y] plot ({15*(2-3*\y)},{15*\y});
\draw[domain=27/41:2/3,smooth,variable=\y] plot ({15*(2-3*\y)},{15*\y});
\node at (-0.8, 10) {\begin{scriptsize}$y=2/3$\end{scriptsize}};
\node (r1) at (2, 10.8) {\begin{footnotesize}$\ref{maxbound}$\end{footnotesize}};
\draw[
    gray, ultra thin, decoration={markings,mark=at position 1 with {\arrow[black,scale=2]{>}}},
    postaction={decorate},
    ]
(r1) to (2,10);

\node (r2) at (6.5, 8.8) {\begin{footnotesize}$\ref{3+2}$\end{footnotesize}};
\draw[
    gray, ultra thin, decoration={markings,mark=at position 1 with {\arrow[black,scale=2]{>}}},
    postaction={decorate},
    ]
(r2) to (5,8);
\draw[
    gray, ultra thin, decoration={markings,mark=at position 1 with {\arrow[black,scale=2]{>}}},
    postaction={decorate},
    ]
    (r2) -- (6,7.5);
\node (r3) at (9,6) {\begin{footnotesize}$\ref{2+3}$\end{footnotesize}};
\draw[
    gray, ultra thin, decoration={markings,mark=at position 1 with {\arrow[black,scale=2]{>}}},
    postaction={decorate},
    ]
(r3) -- (7.5,6.5);

\node (r4) at (10.4,4) {\begin{footnotesize}$\ref{4-7+2-7}$\end{footnotesize}};
\draw[
    gray, ultra thin, decoration={markings,mark=at position 1 with {\arrow[black,scale=2]{>}}},
    postaction={decorate},
    ]
(r4) -- ({15*5/9},{15/4});

\node (r5) at (10.8,2.2) {\begin{footnotesize}$\ref{4-7+2-7}$\end{footnotesize}};
\draw[
    gray, ultra thin, decoration={markings,mark=at position 1 with {\arrow[black,scale=2]{>}}},
    postaction={decorate}, 
    ]
(r5) to [bend right=30] ({15*7/11},{15*2/11});

\node (r6) at (11,1.25) {\begin{footnotesize}$\ref{maxbound}$\end{footnotesize}};
\draw[
    gray, ultra thin, decoration={markings,mark=at position 1 with {\arrow[black,scale=2]{>}}},
    postaction={decorate},
    ]
(r6) -- ({15*2/3},{1.25});

\node (r7) at (5.7, 10.8) {\begin{footnotesize}$\ref{hompepsi}$\end{footnotesize}};
\draw[
    gray, ultra thin, decoration={markings,mark=at position 1 with {\arrow[black,scale=2]{>}}},
    postaction={decorate},
    ]
(r7) to ({2.9},{9.55});
\draw[
    gray, ultra thin, decoration={markings,mark=at position 1 with {\arrow[black,scale=2]{>}}},
    postaction={decorate},
    ]
(r7) to ({4.7},{8.4});

\node (s1) at (1.3,8.3) {\begin{footnotesize}$\ref{psi23reversecurve}$\end{footnotesize}};
\draw[
    gray, ultra thin, decoration={markings,mark=at position 1 with {\arrow[black,scale=2]{>}}},
    postaction={decorate},
    ]
(s1) -- ({15*2/9},{15*5/9});
\draw[
    gray, ultra thin, decoration={markings,mark=at position 1 with {\arrow[black,scale=2]{>}}},
    postaction={decorate},
    ]
(s1) -- ({15*1/6},{15*11/18});
\draw[
    gray, ultra thin, decoration={markings,mark=at position 1 with {\arrow[black,scale=2]{>}}},
    postaction={decorate},
    ]
(s1) to [bend left =10] ({15*1/180},{15*119/180});

\node (s2) at (6.5,2.7) {\begin{footnotesize}$\ref{phi23extracurve}$\end{footnotesize}};
\draw[
    gray, ultra thin, decoration={markings,mark=at position 1 with {\arrow[black,scale=2]{>}}},
    postaction={decorate},
    ]
(s2) -- ({15*11/21},{15*5/21});

\node (s3) at (8.5,1) {\begin{footnotesize}$\ref{psi23extra}$\end{footnotesize}};
\draw[
    gray, ultra thin, decoration={markings,mark=at position 1 with {\arrow[black,scale=2]{>}}},
    postaction={decorate},
    ]
(s3) -- ({15*19/31},{15*4/31});

\node (s4) at (0.8,9.5) {\begin{footnotesize}$\ref{psi23extra}$\end{footnotesize}};
\draw[
    gray, ultra thin, decoration={markings,mark=at position 1 with {\arrow[black,scale=2]{>}}},
    postaction={decorate},
    ]
(s4) -- ({15*4/33},{15*21/33});

%%%%%%%%%%%%%%%%%%%%%%%%%%%%%%%%%%%%%%%%%%%%%%%%%%%%%%%%%%%%%%%%%%%%%%%%%%%%%%%%%%%%%%%%%%%%%%%%%%%%%%%%%%%
\end{scope}
\draw[->] (0,0) -- (11,0) node[anchor=north]{$x$};
\draw[->] (0,0) -- (0,11) node[anchor=east] {$y$};
\node at (8.3,7.7) {\begin{scriptsize}$(1/2,1/2)$\end{scriptsize}};
\node at (15*5/9+.9,60/13+.2) {\begin{scriptsize}$(5/9,4/13)$\end{scriptsize}};
\node at (4.25,7.2) {\begin{scriptsize}$(1/3,1/2)$\end{scriptsize}};
\node at (-1,10) {\begin{scriptsize}$y=2/3$\end{scriptsize}};
\draw[domain=0:(3-sqrt(2))/7,smooth,variable=\x,thick, dotted] plot ({15*\x},{10});
\draw[domain=(3-sqrt(2))/7:.35855, smooth, variable=\x, thick, dotted] plot ({15*\x}, {15*(1-\x)*(1-\x)/(1-2*\x*\x)});
\draw[domain=.35855:.36701, smooth, variable=\x, thick,dotted] plot ({15*\x},{15*(3*\x-2)*(3*\x-2)/(12*\x*\x)});
\draw[very thick,dotted] ({15*.36701},{15*.5}) -- ({15*.5},{15*.5}) -- ({15*.5},{15*.36701});
\draw[domain=.35855:.36701, smooth, variable=\y, thick, dotted] plot ({15*(3*\y-2)*(3*\y-2)/(12*\y*\y)},{15*\y});

\draw[domain=.34:.35855, smooth, variable=\y, thick, dotted] plot ({15*(1-\y)*(1-\y)/(1-2*\y*\y)},{15*\y});
\draw[domain=(3-sqrt(2))/7:.305, smooth, variable=\y, thick, dotted] plot ({15*(1-\y)*(1-\y)/(1-2*\y*\y)},{15*\y});
\draw[domain=0:(3-sqrt(2))/7,smooth,variable=\y,thick, dotted] plot ({10},{15*\y});

\node at (5,5) {\large{$\phi<2/3$}};
\node at (9,9) { \large{$\phi\ge 2/3$}};
\node at (5.15,8) {$?$};
\node at (7.9, 5.1) {$?$};
\draw
(15*65/98, 15*1/35) -- (15*65/98, 15*1/13) -- 
(15*63/95, 15*1/13) -- (15*63/95, 15*1/12) -- 
(15*59/89, 15*1/12) -- (15*59/89, 15*1/11) -- 
(15*51/77, 15*1/11) -- (15*51/77, 15*1/10) -- 
(15*35/53, 15*1/10) -- (15*35/53, 15*1/9) -- 
(15*52/79, 15*1/9) -- (15*52/79, 15*8/71) -- 
(15*25/38, 15*8/71) -- (15*25/38, 15*4/35) -- 
(15*48/73, 15*4/35) -- (15*48/73, 15*8/69) -- 
(15*23/35, 15*8/69) -- (15*23/35, 15*2/17) -- 
(15*44/67, 15*2/17) -- (15*44/67, 15*8/67) -- 
(15*65/99, 15*8/67) -- (15*65/99, 15*3/25) -- 
(15*21/32, 15*3/25) -- (15*21/32, 15*4/33) -- 
(15*61/93, 15*4/33) -- (15*61/93, 15*6/49) -- 
(15*19/29, 15*6/49) -- (15*19/29, 15*1/8) -- 
(15*36/55, 15*1/8) -- (15*36/55, 15*8/63) -- 
(15*53/81, 15*8/63) -- (15*53/81, 15*6/47) -- 
(15*17/26, 15*6/47) -- (15*17/26, 15*4/31) -- 
(15*49/75, 15*4/31) -- (15*49/75, 15*3/23) -- 
(15*15/23, 15*3/23) -- (15*15/23, 15*1/7) -- 
(15*31/48, 15*1/7) -- (15*31/48, 15*12/83) -- 
(15*51/79, 15*12/83) -- (15*51/79, 15*10/69) -- 
(15*20/31, 15*10/69) -- (15*20/31, 15*8/55) -- 
(15*29/45, 15*8/55) -- (15*29/45, 15*6/41) -- 
(15*47/73, 15*6/41) -- (15*47/73, 15*5/34) -- 
(15*9/14, 15*5/34) -- (15*7/11, 15*1/6);

\draw
(15*5/9,15*3/10) -- (15*5/9,15*4/13)--(15*6/11,15*4/13) ;
\draw(15*1/2, 15*1/3) -- (15*1/2, 15*1/2) 
-- (15*1/3, 15*1/2) ;
\draw
(15*4/13,15*6/11)--(15*4/13, 15*5/9) -- (15*3/10, 15*5/9) ;
\draw
(15*1/6, 15*7/11) -- (15*5/34, 15*9/14) 
-- (15*5/34, 15*47/73) -- (15*6/41, 15*47/73) 
-- (15*6/41, 15*29/45) -- (15*8/55, 15*29/45) 
-- (15*8/55, 15*20/31) -- (15*10/69, 15*20/31) 
-- (15*10/69, 15*51/79) -- (15*12/83, 15*51/79) 
-- (15*12/83, 15*31/48) -- (15*1/7, 15*31/48) 
-- (15*1/7, 15*15/23) -- (15*3/23, 15*15/23) 
-- (15*3/23, 15*49/75) -- (15*4/31, 15*49/75) 
-- (15*4/31, 15*17/26) -- (15*6/47, 15*17/26) 
-- (15*6/47, 15*53/81) -- (15*8/63, 15*53/81) 
-- (15*8/63, 15*36/55) -- (15*1/8, 15*36/55) 
-- (15*1/8, 15*19/29) -- (15*6/49, 15*19/29) 
-- (15*6/49, 15*61/93) -- (15*4/33, 15*61/93) 
-- (15*4/33, 15*21/32) -- (15*3/25, 15*21/32) 
-- (15*3/25, 15*65/99) -- (15*8/67, 15*65/99) 
-- (15*8/67, 15*44/67) -- (15*2/17, 15*44/67) 
-- (15*2/17, 15*23/35) -- (15*8/69, 15*23/35) 
-- (15*8/69, 15*48/73) -- (15*4/35, 15*48/73) 
-- (15*4/35, 15*25/38) -- (15*8/71, 15*25/38) 
-- (15*8/71, 15*52/79) -- (15*1/9, 15*52/79) 
-- (15*1/9, 15*35/53) -- (15*1/10, 15*35/53) 
-- (15*1/10, 15*51/77) -- (15*1/11, 15*51/77) 
-- (15*1/11, 15*59/89) -- (15*1/12, 15*59/89) 
-- (15*1/12, 15*63/95) -- (15*1/13, 15*63/95) 
-- (15*1/13, 15*65/98) -- (15*1/35, 15*65/98);
\draw[domain=1/2:6/11,smooth,variable=\x] plot ({15*\x},{15*(2-3*\x)/(5-7*\x)});
\draw[domain=5/9:7/11,smooth,variable=\x] plot ({15*\x},{15*(2-3*\x)/(5-7*\x)});
\draw[domain=65/98:2/3,smooth,variable=\x] plot ({15*\x},{15*(2-3*\x)/(5-7*\x)});

\draw[domain=1/2:6/11,smooth,variable=\y] plot ({15*(2-3*\y)/(5-7*\y)},{15*\y});
\draw[domain=5/9:7/11,smooth,variable=\y] plot ({15*(2-3*\y)/(5-7*\y)},{15*\y});
\draw[domain=65/98:2/3,smooth,variable=\y] plot ({15*(2-3*\y)/(5-7*\y)},{15*\y});

\node (t1) at (6.3,9.3) {\begin{footnotesize}$\ref{23weightedthm}$\end{footnotesize}};
\draw[
    gray, ultra thin, decoration={markings,mark=at position 1 with {\arrow[black,scale=2]{>}}},
    postaction={decorate},
    ]
%(t1) to ({15*1/3},{15*(1-(1/3)*(2-1/3)/(1+1/3))});
(t1) to ({15*1/3},{15*4/7});

\node (t2) at (2, 10.8) {\begin{footnotesize}$\ref{maxbound}$\end{footnotesize}};
\draw[
    gray, ultra thin, decoration={markings,mark=at position 1 with {\arrow[black,scale=2]{>}}},
    postaction={decorate},
    ]
(t2) to (2,10);

\node (t3) at (8.7, 6.7) {\begin{footnotesize}$\ref{hompe4}$\end{footnotesize}};
\draw[
    gray, ultra thin, decoration={markings,mark=at position 1 with {\arrow[black,scale=2]{>}}},
    postaction={decorate},
    ]
(t3) to (7.5,7);
\draw[
    gray, ultra thin, decoration={markings,mark=at position 1 with {\arrow[black,scale=2]{>}}},
    postaction={decorate},
    ]
(t3) to (8,5.43);

\node (u1) at (7,3.05) {\begin{footnotesize}$\ref{phi23curve}$\end{footnotesize}};
\draw[
    gray, ultra thin, decoration={markings,mark=at position 1 with {\arrow[black,scale=2]{>}}},
    postaction={decorate},
    ]
(u1) to ({15*11/21},{15*9/28});
\draw[
    gray, ultra thin, decoration={markings,mark=at position 1 with {\arrow[black,scale=2]{>}}},
    postaction={decorate},
    ]
(u1) to ({15*3/5},{15*1/4});
\draw[
    gray, ultra thin, decoration={markings,mark=at position 1 with {\arrow[black,scale=2]{>}}},
    postaction={decorate},
    ]
(u1) to [bend right=10] ({15*119/180},{15*1/180});

\end{tikzpicture}

\caption{When $\psi(x,y)<2/3$ and when $\phi(x,y)<2/3$.} \label{fig:phi23}
\end{figure}
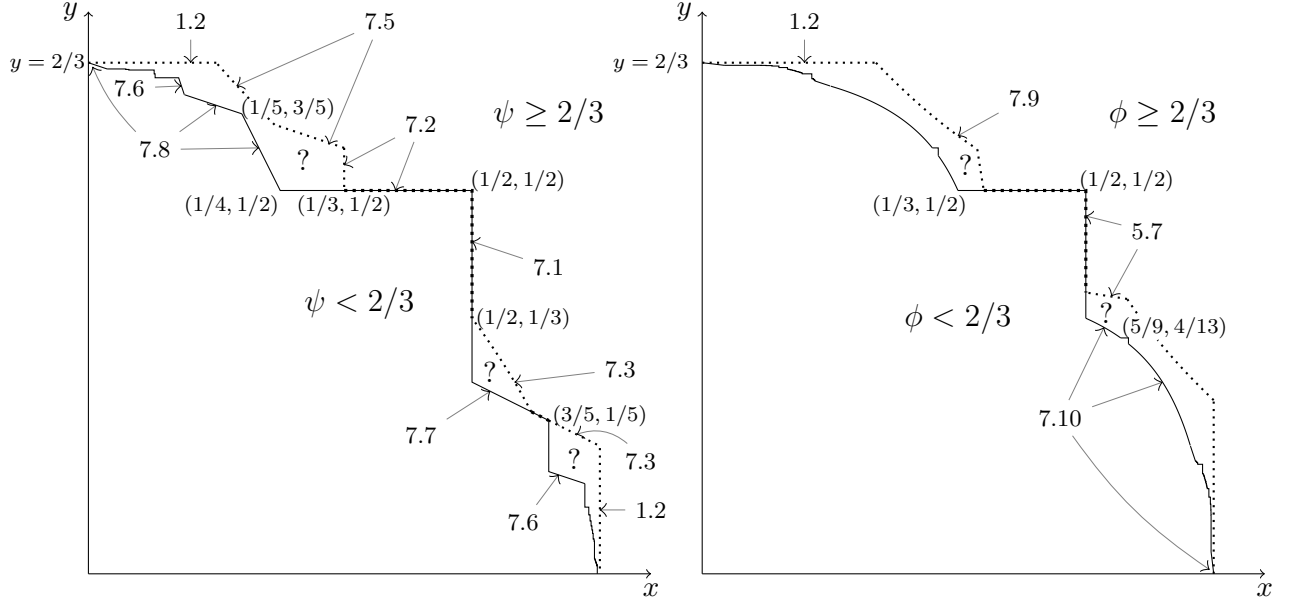

\begin{thm}\label{2+3}
If $x>1/2$ then $\psi(x,1/3)\ge 2/3$.
\end{thm}
\Proof
Let $G$ be $(x,1/3)$-biconstrained via $(A,B,C)$, and suppose for a contradiction that
$|N^2_A(v)|<2|A|/3$ for all $v\in C$. By averaging, there exists $v_0\in A$ such that
$|N^2_C(v_0)|<2|C|/3$. Let $B_0 = N(v_0)$ and $C_0=N^2_C(v_0)$. Hence $|B_0|\ge x|B|$, and $|C_0|<2|C|/3$,
and there are no edges between $B_0$ and $C\setminus C_0$, and every vertex in $C_0$ has a neighbour in $B_0$.

Choose $v_1\in C_0$. Thus $N(v_1)\cap B_0\ne \emptyset$. Let $B_1 = N(v_1)$ and $A_1 = N^2_A(v_1)$.
So $|A_1|\ge x|A|$, and $|A_1|<2|A|/3$. Every vertex $v\in A\setminus A_1$ has a neighbour in $B_0$,
since $|B\setminus B_0|<|B|/2<|N(v)|$. Consequently every vertex in $A\setminus A_1$ has at least $|C|/3\ge |C_0|/2$
second neighbours in $C_0$, and by averaging it follows that some vertex $v_2\in C_0$
has at least $|A\setminus A_1|/2$ second neighbours in $A\setminus A_1$. Let $B_2 = N(v_2)$ and $A_2 = N^2_A(v_2)$.
Then $|A_2\setminus A_1|\ge |A\setminus A_1|/2\ge |A|/6$. If there exists $u\in B_1\cap B_2$, then since
$u$ has at least $x|A|$ neighbours in $A_1$, and they all belong to $A_2$, it follows that
$$|A_2|=|A_2\cap A_1|+|A_2\setminus A_1| \ge x|A|+|A|/6\ge 2|A|/3,$$
a contradiction. Consequently $B_1\cap B_2=\emptyset$.

In particular, $|B_1\cup B_2|\ge 2|B|/3$, and so every vertex in $A$ has a neighbour in $B_1\cup B_2$;
and so $A_1\cup A_2=A$. Since $|A_1|,|A_2|<2|A|/3$, it follows that $|A_1\cap A_2|<|A|/3$. For $i = 1,2$,
choose $b_i\in B_i\cap B_0$. Then $N(b_i)\cap A\subseteq A_i$ for $i = 1,2$, and so
$|N(b_1)\cap N(b_2)\cap A|<|A|/3$. Consequently $|(N(b_1)\cup N(b_2))\cap A|>2|A|/3$.
Since $b_1,b_2\in B_0$ and they each have at least $|C|/3$ neighbours in $C_0$, and $|C_0|<2|C|/3$,
it follows that they have a common neighbour $v\in C_0$. But then $N(b_1)\cup N(b_2)\cap A\subseteq N^2_A(v)$,
and so $|N^2_A(v)|\ge 2|A|/3$, a contradiction. This proves \ref{2+3}.~\bbox

\begin{thm}\label{3+2}
If $y> 1/2$ then $\psi(1/3,y)\ge 2/3$.
\end{thm}
\Proof
Let $G$ be $(1/3,y)$-biconstrained via $(A,B,C)$.
There exists $v_1\in C$ with at least $y|B|$ neighbours in $B$ (in fact, every vertex in $C$ has this property).
Let $B_1 = N(v_1)$ and $A_1 = N^2_A(v_1)$. Thus $|B_1|\ge y|B|$. Since every vertex in
$B\setminus B_1$ has at least $y|C|$ neighbours in $C$, some vertex $v_2\in C$ has at least $y|B\setminus B_1|$
neighbours in $B\setminus B_1$. Let $B_2 = N(v_2)$ and $A_2 = N^2_A(v_2)$. Thus $|B_2\setminus B_1|\ge y|B\setminus B_1|$,
and so
$$|B_1\cup B_2|\ge |B_1|+ y|B\setminus B_1|=y|B|+ (1-y)|B_1|\ge y|B|+y(1-y)|B|= (2-y)y|B|>3|B|/4.$$
In particular, since every vertex in $A$ has at least $|B|/3$ neighbours in $B$, it follows that $A_1\cup A_2=A$.
But $B_1\cap B_2\ne \emptyset$, since $|B_1|,|B_2|\ge y|B|>|B|/2$, and so there exists $b\in B_1\cap B_2$;
and since $b$ has at least $|A|/3$ neighbours in $A$, and they all belong to $A_1\cap A_2$, it follows that
$|A_1\cap A_2|\ge |A|/3$. Since $|A_1\cup A_3|=|A|$, it follows that $|A_1|+|A_2|\ge 4|A|/3$, and so
one of $|A_1|,|A_2|$ is at least $2|A|/3$. This proves \ref{3+2}.~\bbox

The last two results are closely related, via reformulation into triangular language.
First, we need some shorthand for results of the form ``if $x'>x$ then $(x',y,z)$ is triangular'';
let us say ``$(x^+,y,z)$ is triangular'' to mean ``$(x',y,z)$ is triangular for all $x'>x$'', and treat the
other two coordinates similarly. We will mix the two systems of notation, in expressions such as
``$(x^{+*},y^+,z)$ is triangular'', meaning ``$(x'^{*},y^+,z)$ is triangular for all $x'>x$''.

Thus, in triangular language, we have the following.
\begin{itemize}
\item $(1/2^{+*},1/3^*,1/3^+)$ is triangular: because \ref{2+3} says that $(x^*,1/3^*,1/3^+)$ is triangular when $x>1/2$.
\item $(1/2^+,1/3^+, 1/3^*)$ is triangular: because
the proof of \ref{3+2} did not use that
every vertex in $C$ has at least $x|B|$ neighbours in $B$, and so it proves that $(1/3^*,1/2^+,1/3^+)$ is triangular,
and rotating gives that $(1/2^+,1/3^+, 1/3^*)$ is triangular.
\item $(1/2^+,1/3^{+*},1/3^*)$ is triangular; this follows from \ref{bisym} with $k=2$ and rotating.
\item $(1/2^+,1/3^*,1/3^{+*})$ is triangular; this also follows from \ref{bisym} with $k=2$ and rotating.
\end{itemize}
These four statements are similar, but no two are equivalent, and
it would be good to find a common strengthening. Note, however, that $(1/2^{+*},1/3^*,1/3^*)$ is not triangular,
and indeed $(2/3^*,1/3^*,1/3^*)$ is not triangular. We have not been able to decide whether
$(1/2^+,1/3^+,1/3)$ and $(1/2^+,1/3,1/3^+)$ are triangular, or indeed whether $(1/2^{+*},1/3^+,1/3^+)$ is triangular.

Pursuing this further, what about $(1/2^+,1/3^+,x)$ when $x<1/3$ (perhaps with some sprinkling of asterisks)?
How small can $x$ be such that the triple remains triangular? We have examples that show that $(5/9, 5/14, 4/13)$
and $(4/7, 3/8, 2/7)$ are not triangular, and $(3/5^*,2/5, 1/4)$ is not triangular,
but as far as we know, $(1/2^{+*},1/3^{+*},1/5)$ might be triangular.

\begin{figure}[ht]
\centering

\begin{tikzpicture}[xscale=0.8,yscale=0.5,auto=left]
\tikzstyle{every node}=[inner sep=1.5pt, fill=black,circle,draw]

\def\r{0}
\node (a1) at (\r,0) {};
\node (a2) at (\r,2) {};
\node (a3) at (\r,4) {};
\node (a4) at (\r,6) {};
\node (a5) at (\r,8) {};
\def\r{3}
\node (b1) at (\r,0) {};
\node (b2) at (\r,2) {};
\node (b3) at (\r,4) {};
\node (b4) at (\r,6) {};
\node (b5) at (\r,8) {};
\def\r{6}
\node (c1) at (\r,0) {};
\node (c2) at (\r,2) {};
\node (c3) at (\r,4) {};
\node (c4) at (\r,6) {};
\node (c5) at (\r,8) {};
\def\r{-3}
\node (d1) at (\r,0) {};
\node (d2) at (\r,2) {};
\node (d3) at (\r,4) {};
\node (d4) at (\r,6) {};
\node (d5) at (\r,8) {};
\foreach \from/\to in {a1/b1,a2/b3,a3/b2,a2/b4, a3/b4,a4/b2,a4/b3,a5/b5,b1/c4,b1/c2,b1/c3,b2/c2,b2/c1,b3/c3,b3/c1,b4/c4,b4/c1,b5/c5, d1/a1,d1/a5,d2/a2,d2/a5,d3/a3,d3/a5,d4/a4,d4/a5,d5/a1,d5/a2,d5/a3,d5/a4}
\draw [-] (\from) -- (\to);
\tikzstyle{every node}=[]
\draw (a1) node [above]           {\scriptsize$\frac17$};
\draw (a2) node [above]           {\scriptsize$\frac17$};
\draw (a3) node [above]           {\scriptsize$\frac17$};
\draw (a4) node [above]           {\scriptsize$\frac17$};
\draw (a5) node [above]           {\scriptsize$\frac37$};
\draw (b1) node [above]           {\scriptsize$\frac27$};
\draw (b2) node [above]           {\scriptsize$\frac17$};
\draw (b3) node [above]           {\scriptsize$\frac17$};
\draw (b4) node [above]           {\scriptsize$\frac17$};
\draw (b5) node [above]           {\scriptsize$\frac27$};
\draw (c1) node [above]           {\scriptsize$\frac14$};
\draw (c2) node [above]           {\scriptsize$\frac18$};
\draw (c3) node [above]           {\scriptsize$\frac18$};
\draw (c4) node [above]           {\scriptsize$\frac18$};
\draw (c5) node [above]           {\scriptsize$\frac38$};
\end{tikzpicture}

\caption{$(4/7,2/7,3/8)$ is not triangular} \label{fig:4-7-2-7}
\end{figure}
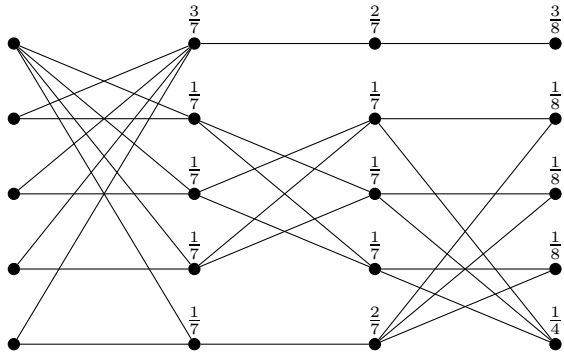

This extends to weighted graphs in the natural way.
For instance, the weighted graph of figure \ref{fig:4-7-2-7} (identify the vertices on the left with those on the right, in order)
 shows that $(4/7,2/7,3/8)$ is not triangular.

The graph of figure \ref{fig:4-7-2-7} shows that $\phi(4/7,2/7)\le 5/8$, and so 
we studied $\psi(4/7,2/7)$, and proved the following.
\begin{thm}\label{4-7+2-7}
If $x,y\in (0,1]$ such that $\max(x,y)>1/2$, $x\ge 1/3$, $x+2y> 1$, and $3x+y/(1-y)> 2$, then
$\psi(x,y)\ge 2/3$.
\end{thm}
\Proof
Let $G$ be $(x,y)$-biconstrained, via $(A,B,C)$, and suppose for a contradiction that $|N_A^2(v)|< 2|A|/3$
for each $v\in C$. By \ref{3+2}, $y\le 1/2$ since $x\ge 1/3$; and so $x>1/2$ since $\max(x,y)>1/2$.
Hence $y<1/3$ by \ref{2+3}. Also $x<2/3$, by \ref{maxbound}.
\\
\\
(1) {\em For all $v_1,v_2\in C$, if $|N(v_1)\cup N(v_2)|>(1-x)|B|$ then $N(v_1)\cap N(v_2)=\emptyset$,
$N^2_A(v_1)\cup N^2_A(v_2)=A$, and $|N^2_A(v_1)\cap N^2_A(v_2)|<|A|/3$.}
\\
\\
Every vertex in $A$ has a neighbour in $N(u)\cup N(v)$, and so $N_A^2(u)\cup N_A^2(v)=A$.
Since $|N_A^2(u)|< 2|A|/3$ and $|N_A^2(v)|< 2|A|/3$ it follows that $|N_A^2(u)\cap N_A^2(v)|< |A|/3$,
and so there is no vertex in $N(u)\cap N(v)$ (since any such vertex would have at least $x|A|$ neighbours in $A$,
all belonging to $N_A^2(u)\cap N_A^2(v)$). This proves (1).
\\
\\
(2) {\em There exist $v_1,v_2\in C$ with $N(v_1)\cap N(v_2)=\emptyset$.}
\\
\\
Choose $v_1\in C$. Since every vertex in $A\setminus N^2_A(v_1)$ has at least $x|B|$ second neighbours in
$B\setminus N(v_1)$,
some vertex $u_2\in B\setminus N(v_1)$ has at least $(x/(3(1-y)))|A|$ neighbours in $A\setminus N^2_A(v_1)$.
Let $v_2\in C$ be adjacent to $u_2$. If $v_1,v_2$ have a common neighbour $u_1$, then since $N_A(u_1)\subseteq N^2_A(v_2)$,
it follows that $|N^2_A(v_2)|\ge (x/(3(1-y))+x)|A|$, and so $x/(3(1-y))+x<2/3$, that is,
$(4-3y)x<2-2y<(4-3y)/2$, and so $x<1/2$, a contradiction.
This proves (2).
\\
\\
(3) {\em If $v_1,v_2,v_3\in C$ and $N(v_1)\cap N(v_2)=\emptyset$ then $N(v_3)$ is disjoint from exactly one of
$N(v_1), N(v_2)$.}
\\
\\
If $N(v_3)$ is disjoint from both $N(v_1), N(v_2)$, then every two of $N(v_1), N(v_2), N(v_3)$ have union of 
cardinality more than $(1-x)|B|$, and so every vertex in $A$ belongs to at least two of $N^2_A(v_i)\;(i=1,2,3)$.
Consequently one of $N^2_A(v_i)\;(i=1,2,3)$ has cardinality at least $a|A|/3$, a contradiction.
Now suppose that $N(u_3)$ has nonempty intersection with both $N(v_1), N(v_2)$. Thus $|N^2_A(v_i)\cap N^2_A(v_3)|\ge x|A|$
for $i = 1,2$, and since $|N^2_A(v_1)\cap N^2_A(v_2)|<|A|/3$, it follows that
$|N^2_A(v_3)|\ge (2x-1/3)|A|\ge 2|A|/3$ since $x\ge 1/2$, a contradiction. This proves (3).

\bigskip
Let $H$ be the bipartite graph $G[B\cup C]$. From (2) and (3), $H$ has exactly two components $H_1$ and $H_2$ say.
Let $C_i=V(H_i)\cup C$ and $B_i = V(H_i)\cap B$ for $i = 1,2$. Then from (3), every two vertices
in $C_i$ have a common neighbour in $B_i$, for $i = 1,2$.
Let $c_i=|C_i|/|C|$, for $i = 1,2$. Thus $c_1+c_2=1$.
We may assume that $b_1\ge 1/2$.
Choose $v_1\in C_1$. Since $|A\setminus N^2_A(v_1)|>|A|/3$, and every vertex in $A\setminus N^2_A(v_1)$ has at least $y|C|$
second neighbours in $C_1$, some vertex $v_2\in C_1$ has at least $(y/(3c_1))|A|$ second neighbours in
$A\setminus N^2_A(v_1)$. But since $v_1,v_2\in C_1$, they have a common neighbour in $B_1$; therefore
$|N^2_A(v_2)|\ge (y/(3c_1)+x)|A|$, and so $y/(3c_1)+x<2/3$. Now $c_1\le 1-y$, so $3x+y/(1-y)<2$, contrary to the
hypothesis.
This proves \ref{4-7+2-7}.~\bbox

If $A\subseteq V(G)$ and $X\subseteq V(G)\setminus A$, $N_A(X)$ denotes the set of vertices in $A$ with a neighbour in $X$.
The next result is a useful lemma which says, roughly speaking, the larger $X$ is, the larger $N_A(X)$ is. In this section we only use it for $k=1$, but we will
use it with $k=2$ in the section discussing when $\psi(x,y)\ge 3/4$.

\begin{thm}\label{hompe1-y}
Let $x,y,z\in (0,1]$, and suppose that $G$ is $(x,y)$-biconstrained via $(A,B,C)$, and $|N_A^2(w)| < z|A|$ for all
$w \in C$.
Then for all integers $k \ge 1$, if $B' \subseteq B$ with
$$\frac{|B'|}{|B|} > (k-1)(1-y)+\max(1-y,1-x/(1-y)),$$
then $|N_A(B')| > (x+k(1-z))|A|$.
\end{thm}

\Proof
We proceed by induction on $k$, and so we assume that either $k=1$ or the result holds for $k-1$.
Let $A' = N_A(B')$.  Every vertex in $B \setminus B'$ has at least $y|C|$ neighbours in $C$, so there exists $v \in C$
with at least $y|B \setminus B'|$ neighbours in $B \setminus B'$.
By hypothesis,
$|B'|/|B| > 1-x/(1-y)$,
that is, $x|B|+y|B\setminus B'|> |B\setminus B'|$.
But every vertex in $A\setminus A'$ has at least $x|B|$ neighbours in $B\setminus B'$, and therefore has one such
neighbour adjacent to $v$; and so $A\setminus A'\subseteq N^2_A(v)$.
Let $B''=B'\cap N_B(v)$, and $A''=N_A(B'')$.
\\
\\
(1) {\em $|A''|\ge (x+(k-1)(1-z))|A|$.}
\\
\\
Since
$|N_B(v)|\ge y|B|$, it follows that
$|B''|\ge y|B|-(|B|-|B'|)=|B'|-(1-y)|B|$. If $k=1$, then $|B'|/|B|>1-y$ by hypothesis,
and therefore $B''\ne \emptyset$, and so $|A''|\ge x|A|$ as claimed. If $k\ge 2$, then since
$|B'|/|B|> (k-1)(1-y)+\max(1-y,1-x/(1-y))$, it follows that
$$\frac{|B''|}{|B|}> (k-2)(1-y)+ \max(1-y,1-x/(1-y)),$$
and so $|A''|\ge (x+(k-1)(1-z))|A|$ from the inductive hypothesis. This proves (1).

\bigskip

Since $A''\subseteq A'$, and $A''\cup (A\setminus A')\subseteq N^2_A(v)$,
it follows that
$$z|A|\ge |N^2_{A}(v)|\ge (x+(k-1)(1-z))|A|+|A\setminus A'|$$
and so $|A'|\ge  (x+k(1-z))|A|$. This proves \ref{hompe1-y}.~\bbox

\begin{thm}\label{hompepsi}
Let $x,y\in (0,1]$, such that $y>1/2$, $x+3y> 2$ and $x>2(1-y)^2/(2-y)$. Then $\psi(x,y)\ge 2/3$.
\end{thm}
\Proof
Let $G$ be $(x,y)$-biconstrained, via $(A,B,C)$, and suppose for a contradiction that $|N_A^2(v)|< 2|A|/3$
for each $v\in C$. 
\\
\\
(1) {\em If $B'\subseteq B$ with $|B'|/|B|> \max(1-y,1-x/(1-y))$, then there are at least $(x+1/3)|A|$ vertices in $A$
with a neighbour in $B'$. In particular, if $|B'|\ge (x+2y-1)|B|$ then the same conclusion holds.}
\\
\\
The first statement follows from \ref{hompe1-y} with $k=1$.
For the second, $x+2y-1>1-y$ since $x+3y\ge 2$;
and $x+2y-1>1-x/(1-y)$ since $x>2(1-y)^2/(2-y)$. Consequently $x+2y-1 > \max(1-y,1-x/(1-y))$. 
This proves the second statement and so proves (1).

\bigskip

Say $w_1,w_2\in C$ are {\em close} if $|N_B(w_1)\cup N_B(w_2)|\le (1-x)|B|$. 
\\
\\
(2) {\em There exists $w_1\in C$ such that the set of vertices in $C$ that are close to $w_1$ has cardinality at least $|C|/2$.}
\\
\\
This is trivial if every two vertices in $C$ are close; so we assume there exist $w_1,w_2\in C$ that are not close. Consequently
every vertex in $A$ has a neighbour in $N_B(w_1)\cup N_B(w_2)$, and so $N_A^2(w_1)\cup N^2_A(w_2)=A$. If there exists $w\in C$
that is not close to either of $w_1,w_2$ then similarly $N_A^2(w_1)\cup N^2_A(w)=A$ and $N_A^2(w_2)\cup N^2_A(w)=A$; 
and so every vertex in $A$
belongs to at least two of $N^2_A(w_1),N^2_A(w_2), N^2_A(w)$, and therefore one of these three sets has cardinality at least $2|A|/3$,
a contradiction. Thus, exchanging $w_1,w_2$ if necessary, we may assume that at least half of all vertices in $C$ are close to
$w_1$. This proves (2).

\bigskip
Let $C_1$ be the set of vertices in $C$ that are close to $w_1$; thus $|C_1|\ge |C|/2$. Let $B_1=N_B(w_1)$ and $A_1=N^2_A(w_1)$.
Since $|B_1|\ge y|B|$ and every vertex in $A\setminus A_1$ has at least $x|B|$ neighbours in $B\setminus B_1$, there exists 
$v_2\in B\setminus B_1$ with at least $\frac{x}{1-y}|A\setminus A_1|$ neighbours in $A\setminus A_1$. Since $y>1/2$ and $|C_1|\ge |C|/2$,
$v_2$ has a neighbour $w_2\in C_1$. Since $w_2$ is close to $w_1$ it follows that $|N_B(w_1)\cup N_B(w_2)|\le (1-x)|B|$, and so
$|N_B(w_1)\cap N_B(w_2)|\ge  (x+2y-1)|B|$. From (1), there are at least $(x+1/3)|A|$ vertices in $A$
with a neighbour in $N_B(w_1)\cap N_B(w_2)$. These vertices all belong to $A_1$, and so 
$$|N^2(w_2)|\ge \frac{x}{1-y}|A\setminus A_1| + \left(x+\frac13\right)|A|\ge \left(\frac{x}{3(1-y)}+x+\frac13\right)|A|.$$
Consequently $\frac{x}{3(1-y)}+x+1/3<2/3$, that is, $\frac{x}{1-y}+3x<1$. By hypothesis, $x>2(1-y)^2/(2-y)$, and so
substitution for $x$ yields 
$$\frac{2(1-y)^2/(2-y)}{1-y}+6(1-y)^2/(2-y)<1,$$ 
which simplifies to $(3y-2)(2y-3)<0$, contrary 
to \ref{maxbound}. This proves \ref{hompepsi}.~\bbox

\begin{thm}\label{psi23extra}
Let $x,y\in (0,1]$. If either
\begin{itemize}
\item $x\le 11/17$ and $y\le 1/7$ and $x+3y\le 1$; or
%$4/7\le x\le 11/17$ and $x+3y\le 1$; or
\item $x<1/8$ and $y\le 11/17$ and $3x+y\le 1$
%$5/8< y\le 11/17$ and $3x+y\le 1$
\end{itemize}
then $\psi(x,y)<2/3$.
\end{thm}
\Proof
Take the graph consisting of seven disjoint copies of a three-vertex path, numbered $a_i,b_i,c_i$ in order 
$(1\le i\le 7)$.
For $1\le i\le 3$ and $4\le j\le 7$, make $a_i$ adjacent to $b_j$ and make $a_j$ adjacent to $b_i$, forming $G$.
Let $A=\{a_i\;:1\le i\le 7\}$ and define $B,C$ similarly.

For the first statement, we may assume by increasing $x,y$ that $x+3y= 1$. It follows that $x\ge 4/7$ 
(because $y\le 1/7$) and similarly $y\ge 2/17$, and $4x=4-12y\ge 1+9y$.
For $1\le i\le 3$, let $w(a_i) = 3/17$, $w(b_i) = (4x-1)/9$, and $w(c_i)=1/7$.
For $4\le i\le 7$, let $w(a_i) = 2/17$, $w(b_i)=(1-x)/3$, and $w(c_i)=1/7$. 
Then this weighted graph is $(x,y)$-biconstrained via $(A,B,C)$ and shows that $\psi(x,y)<2/3$. This proves the first statement.

For the second statement, we may assume that $3x+y=1$, and so $x\ge 2/17$ and $y>5/8$, and 
so $3y=3-9x\le 1+8x$.
Let us take the same graph and redefine $w$, as follows. 
%Choose $p$ with $1/6<p<5/27$ and $x\le p\le (1-4x)/3$. 
For $1\le i\le 3$, let $w(a_i) = (1-y)/2$ and $w(b_i) = w(c_i)= (1-4x)/3$.
For $4\le i\le 7$, let $w(a_i) = (3y-1)/8$ and $w(b_i)=w(c_i)=x$. 
Then this weighted graph is $(x,y)$-biconstrained via $(C,B,A)$ and shows that $\psi(x,y)<2/3$. This proves the second statement,
and hence proves \ref{psi23extra}.~\bbox

\begin{thm}\label{phi23extracurve}
Let $x',y',z'\in (0,1]$ such that $\psi(x',y')\le z'<1/2$; and let $x,y\in (0,1]$ satisfy
$x\le \frac{1}{2-x'}$, $x< 1-\frac{1-x'}{3(1-z')}$, $y\le  \frac{y'}{1+y'}$ and $x+\frac{1-x'}{y'}y\le 1$.
Then $\psi(x,y)< 2/3$. Consequently:
\begin{itemize}
\item $\psi(x,y)<2/3$ if $x\le 3/5$ and $y\le 1/4$ and $x+2y\le 1$;
\item $\psi(x,y)<2/3$ if $x\le 5/8$ and $y\le 1/6$ and $x+3y\le 1$.
\end{itemize}
\end{thm}
\Proof
The first statement follows from \ref{psilb1} taking $z$ slightly less than $2/3$.
The two statements in bullets follow by setting $x'=y'=z'=1/3$, and
then $x'=z'=2/5$ and $y'=1/5$. This proves \ref{phi23extracurve}.~\bbox

\begin{thm}\label{psi23reversecurve}
Let $x',y',z'\in (0,1]$ such that $\psi(x',y')\le z'<1/2$; and let $x,y\in (0,1]$ satisfy
$x< \frac{2x'}{3}$, $y\le \frac{1}{2-y'}$, and $\frac{1-y'}{x'}x+y\le 1$.
Then $\psi(x,y)< 2/3$. Consequently:
\begin{itemize}
\item $\psi(x,y)<2/3$ if $y\le 3/5$ and $2x+y\le 1$; and
\item $\psi(x,y)<2/3$ if $y\ge 3/5$ and $x+3y\le 2$, and $x+3y< 2$ if $x$ or $y$ is irrational.
\end{itemize}
\end{thm}
\Proof
The first statement follows from \ref{psilb2}, taking $z$ slightly less than $2/3$.
To prove the first bullet, let $2x+y\le 1$, and so $x\le 1/2$. If also $y\le 1/2$ then $\psi(x,y)\le 1/2$
by \ref{intk}; so we may assume that $y>1/2$.
We claim there is an integer $k\ge 1 $ with
$$\frac{3-3y}{6y-2}< k\le \frac{1-y}{4y-2}.$$
%$$\frac{1}{3y-1}-\frac12\le k\le \frac{1}{8y-4}-\frac14.$$
To see this, if $y> 5/9$ we can take $k=1$ (because we are given that $y\le 3/5$), and if $y\le 5/9$, then
$$\frac{1-y}{4y-2} - \frac{3-3y}{6y-2}\ge 1,$$
and so again $k$ exists. Thus 
$y\le \frac{2k+1}{4k+1}$, and
$x\le \frac{1-y}{2} < \frac{2k}{6k+3}.$
Let $x'=z'=k/(2k+1)$, and $y'=1/(2k+1)$. Then the claim follows from the first statement.

For the second bullet, let $x+3y\le 2$ with $y\ge 3/5$, with $x+3y<2$ if $x$ or $y$ is irrational.
Consequently, we may assume that $x,y$ are rational, by increasing them slightly if necessary. 
Let $x'=2/y-3$, and $y'=2-1/y$; it follows that $x'+2y'\le 1$ and
$2x'+y'\le 1$, and $x',y'$ are rational,
and so $\psi(x',y')<1/2$ by \ref{psi12curve}. The result follows from the first statement.
This proves \ref{psi23reversecurve}.~\bbox

For the mono-constrained question, we have:
\begin{thm}\label{23weightedthm}
For $x,y\in (0,1]$, if $y\le 1/2$ and $x> (1-y)^2/(1-2y^2)$ then $\phi(x,y)\ge 2/3$.
\end{thm}
\Proof 
Let $G$ be $(x,y)$-constrained via $(A,B,C)$. 
If $x+y>1$ the result follows from \ref{trivial}, so we may assume that $x+y\le 1$. Since
$x> (1-y)^2/(1-2y^2)$, we may also assume that 
\begin{itemize}
\item $x,y$ are rational; and
\item every vertex in $A$ has strictly more than $x|B|$ neighbours in $B$
\end{itemize}
by reducing $x$ and $y$ a little if necessary while retaining the
property that $x> (1-y)^2/(1-2y^2)$. 

Let $p=(1-x-y)/(1-2y)$. Thus $p$ is rational, so we may assume (by multiplying vertices) that $p|B|$ is an integer. 
Also $p\le y$, since $x> (1-y)^2/(1-2y^2)$.
Let $s=(x-(1-y)^2)/(y(1-y))$. It follows that $0\le s\le 1$, since $x> (1-y)^2/(1-2y^2)$ and $x+y\le 1$.

Choose $v_1\in C$ with at
least $y|B|$ neighbours in $B$, and let
$B_1\subseteq N(v_1)$ with $|B_1|=y|B|$.
Choose $v_2\in C$ such that 
$sb_0+b_2\ge y(sy+(1-y))$,
where $b_0|B|=|N(v_2)\cap B_1|$ and $b_2|B|=|N(v_2)\setminus B_1|$. 
(Such a vertex exists by averaging.)
We claim that $b_0+b_2\ge  p$; for from the definition of $s$,
$$b_0+b_2\ge sb_0+b_2\ge y(sy+(1-y))=y((x-(1-y)^2)/(1-y)+1-y) = xy/(1-y),$$
and $p=(1-x-y)/(1-2y)\le xy/(1-y)$ since $x> (1-y)^2/(1-2y^2)$.

Also we claim that $b_2\ge  1-x-y$; for from the definition of $s$,
$$sy+1-x-y=  y(sy+1-y)\le sb_0+b_2\le sy+b_2.$$
Consequently $|N(v_1)\cup N(v_2)|\ge (1-x)|B|$, and there exist $P_1\subseteq N(v_1)$ and $P_2\subseteq N(v_2)$, 
both of cardinality $p|B|$. Choose $v_3\in C$ with at least $y(1-2p)|B|$ neighbours in $B\setminus (P_1\cup P_2)$.
Then for $i = 1,2$,
$$|P_i\cup N(v_3)|\ge (y(1-2p) +p)|B|\ge (1-x)|B|.$$
Since every vertex in $A$ has strictly more than $x|B|$ neighbours in $B$, it follows that every vertex in $A$ belongs to at least
two of the sets $N^2_A(v_i)\;(i = 1,2,3)$; and so one of these sets has cardinality at least $2|A|/3$. This proves \ref{23weightedthm}.~\bbox

\begin{thm}\label{phi23curve}
For all $x,y\in (0,1]$ with $y\le 1/2$,
if $\frac{x}{1-x}+ \frac{y}{1-2y}\le 2$, then $\phi(x,y)<2/3$.
\end{thm}
\Proof We may assume that $x>1/2$, or else the result is true since $\phi(1/2,1/2)=1/2$.
The claim follows from \ref{philb} taking $z$ slightly less than $2/3$.
This proves \ref{phi23curve}.~\bbox

\section{The $1/3$ level}

Next we do the same for $\psi(x,y)\ge 1/3$ and $\phi(x,y)\ge 1/3$. The figure summarizes our results.
\begin{figure}[H]
\centering

\begin{tikzpicture}[scale=1.2,auto=left]

\begin{scope}[shift ={(-6.5,0)}]
\draw[->] (0,0) -- (5.5,0) node[anchor=north]{$x$};
\draw[->] (0,0) -- (0,5.5) node[anchor=east] {$y$};
\node at (4.25,3.9) {\begin{scriptsize}$(1/4,1/4)$\end{scriptsize}};
\node at (2.25,4.8) {\begin{scriptsize}$(2/19,6/19)$\end{scriptsize}};
\node at (2.5,2.5) {\large{$\psi<1/3$}};
\node at (4.5,4.5) { \large{$\psi\ge 1/3$}};
\node at (-.7,5) {\begin{scriptsize}$y=1/3$\end{scriptsize}};
\node at (2.2,4.5) {$?$};
\node at (30/7+.2,2.6) {$?$};
\draw[thick, dotted] (0,5)-- (15/6,15/3);
\draw[domain=0.169:1/6,smooth,variable=\x, thick, dotted] plot ({15*\x},{15*(1-4*\x)/(2-6*\x)});
\draw[domain=0.177:sqrt(10)/6-1/3,smooth,variable=\x, thick, dotted] plot ({15*\x},{15*(1-4*\x)/(2-6*\x)});
\draw[domain=sqrt(10)/6-1/3:1/4,smooth,variable=\x,very thick, dotted] plot ({15*\x},{15*(1-\x)/3});
\draw[domain=1/4:2/9,smooth,variable=\y,very thick, dotted] plot ({15*(1-\y)/3},{15*\y});
\draw[domain=7/27:19/63,smooth,variable=\x,thick, dotted] plot ({15*\x},{15*2/9});
\draw[domain=2/9:1/5,smooth,thick,dotted,variable=\y] plot ({5*(1-\y/(3*(1-\y)))},{15*\y});
\draw[domain=11/36:1/3,smooth,thick,dotted,variable=\x] plot ({15*\x},{3});
\draw[domain=1/5:0, smooth, variable=\y, thick, dotted] plot ({5},{15*\y});

\draw[domain=2/19:1/4,smooth,variable=\y] plot ({15*(1-\y)/3}, {15*\y});
\draw[domain=0:1/46,smooth,variable=\y] plot ({15*(1-\y)/3}, {15*\y});
\draw[domain=2/19:1/4,smooth,variable=\x] plot ({15*\x},{15*(1-\x)/3});
\draw[domain=0:1/46,smooth,variable=\x] plot ({15*\x},{15*(1-\x)/3});

%$y>1/5$ and $x+4y>1$ and $2x+2y>1$ and $1/(2(1-y))\ge 1/x - 3$ and $3x+y/(3(1-y))\ge 1$

%insert start

\draw(15/3,0) -- (15*22/67, 15*1/67)--(15*22/67,15*2/67) --
(15*22/67, 15*1/70) -- (15*22/67, 15*2/67) -- 
(15*18/55, 15*2/67) -- (15*18/55, 15*2/55) -- 
(15*11/34, 15*2/55) -- (15*11/34, 15*1/17) -- 
(15*10/31, 15*1/17) -- (15*10/31, 15*2/31) -- 
(15*6/19, 15*2/31) -- (15*6/19, 15*2/19) -- (15*17/57, 15*2/19);

\draw (15*2/19,15*17/57) -- (15*2/19, 15*6/19) -- (15*2/31, 15*6/19)
-- (15*2/31, 15*10/31) -- (15*1/17, 15*10/31) 
-- (15*1/17, 15*11/34) -- (15*2/55, 15*11/34) 
-- (15*2/55, 15*18/55) -- (15*2/67, 15*18/55) 
-- (15*2/67, 15*22/67) -- (15*1/67, 15*22/67) -- (0,15/3);
%insert end
%\draw
%(15*15/46,15*1/46) -- (15*15/46, 15*1/23) -- 
%(15*14/43, 15*1/23) -- (15*14/43, 15*2/43) -- 
%(15*13/40, 15*2/43) -- (15*13/40, 15*1/20) -- 
%(15*12/37, 15*1/20) -- (15*12/37, 15*2/37) -- 
%(15*11/34, 15*2/37) -- (15*11/34, 15*1/17) -- 
%(15*10/31, 15*1/17) -- (15*10/31, 15*2/31) -- 
%(15*6/19, 15*2/31) -- (15*6/19, 15*2/19) --  (15*17/57, 15*2/19);

%\draw (15*2/19,15*17/57) -- (15*2/19, 15*6/19) -- (15*2/31, 15*6/19) 
%-- (15*2/31, 15*10/31) -- (15*1/17, 15*10/31) 
%-- (15*1/17, 15*11/34) -- (15*2/37, 15*11/34) 
%-- (15*2/37, 15*12/37) -- (15*1/20, 15*12/37) 
%-- (15*1/20, 15*13/40) -- (15*2/43, 15*13/40) 
%-- (15*2/43, 15*14/43) -- (15*1/23, 15*14/43) 
%-- (15*1/23, 15*15/46) -- (14*1/46, 15*15/46);

\node (r1) at (15/8, 5.6) {\begin{footnotesize}$\ref{maxbound}$\end{footnotesize}};
\draw[
    gray, ultra thin, decoration={markings,mark=at position 1 with {\arrow[black,scale=2]{>}}},
    postaction={decorate},
    ]
(r1) to (15/8,5);

\node (r2) at (5.6, 15/6) {\begin{footnotesize}$\ref{maxbound}$\end{footnotesize}};
\draw[
    gray, ultra thin, decoration={markings,mark=at position 1 with {\arrow[black,scale=2]{>}}},
    postaction={decorate},
    ]
(r2) to (5,15/6);

\node (r3) at (3.5, 5) {\begin{footnotesize}$\ref{bisym}$\end{footnotesize}};
\draw[
    gray, ultra thin, decoration={markings,mark=at position 1 with {\arrow[black,scale=2]{>}}},
    postaction={decorate},
    ]
(r3) to ({15*.18},{15*(1-4*.18)/(2-6*.18)});
\draw[
    gray, ultra thin, decoration={markings,mark=at position 1 with {\arrow[black,scale=2]{>}}},
    postaction={decorate},
    ]
(r3) to ({15*.22},{15*.26});

\node (r4) at (5.2, 3.8) {\begin{footnotesize}$\ref{semibi}$\end{footnotesize}};
\draw[
    gray, ultra thin, decoration={markings,mark=at position 1 with {\arrow[black,scale=2]{>}}},
    postaction={decorate},
    ]
(r4) to ({15*.255},{15*.235});
\draw[
    gray, ultra thin, decoration={markings,mark=at position 1 with {\arrow[black,scale=2]{>}}},
    postaction={decorate},
    ]
(r4) to [bend right = 20] ({15*.29},{15*2/9});

\node (r5) at (5.4, 3.3) {\begin{footnotesize}$\ref{psi13good}$\end{footnotesize}};
\draw[
    gray, ultra thin, decoration={markings,mark=at position 1 with {\arrow[black,scale=2]{>}}},
    postaction={decorate},
    ]
(r5) to ({15*41/135},{15*4/19});
\draw[
    gray, ultra thin, decoration={markings,mark=at position 1 with {\arrow[black,scale=2]{>}}},
    postaction={decorate},
    ]
(r5) to [bend right = 20] ({15*23/72},{15/5});

\node (s1) at (15/16, 15/4) {\begin{footnotesize}$\ref{psi12curve}$\end{footnotesize}};
\draw[
    gray, ultra thin, decoration={markings,mark=at position 1 with {\arrow[black,scale=2]{>}}},
    postaction={decorate},
    ]
(s1) to [bend right=20] ({15/8},{15*7/24});
\draw[
    gray, ultra thin, decoration={markings,mark=at position 1 with {\arrow[black,scale=2]{>}}},
    postaction={decorate},
    ]
(s1) to [bend left=20] ({15/90},{15*89/270});

\node (s2) at (15/4,15/16) {\begin{footnotesize}$\ref{psi12curve}$\end{footnotesize}};
\draw[
    gray, ultra thin, decoration={markings,mark=at position 1 with {\arrow[black,scale=2]{>}}},
    postaction={decorate},
    ]
(s2) to [bend left=20] ({15*7/24},{15/8});
\draw[
    gray, ultra thin, decoration={markings,mark=at position 1 with {\arrow[black,scale=2]{>}}},
    postaction={decorate},
    ]
(s2) to [bend right=20] ({15*89/270},{15/90});

\end{scope}
%%%%%%%%%%%%%%%%%%%%%%%%%%%%%%%%%%%%%%%%%%%%%%%%%%%%%%%%%%%%%%%%%%%%%%%%%%%%%%%%%%%%%%%%%%%%%%%%%%%%%%%%%%%%%%%%%%%%%%%%%%%%

\draw[->] (0,0) -- (6,0) node[anchor=north]{$x$};
\draw[->] (0,0) -- (0,6) node[anchor=east] {$y$};
\draw 
(15*32/97, 15*.075111) -- (15*32/97, 15*1/10) -- 
(15*18/55, 15*1/10) -- (15*18/55, 15*1/9) -- 
(15*17/52, 15*1/9) -- (15*17/52, 15*1/8) -- 
(15*29/90, 15*1/8) -- (15*29/90, 15*12/95) -- 
(15*19/59, 15*12/95) -- (15*19/59, 15*8/63) -- 
(15*28/87, 15*8/63) -- (15*28/87, 15*6/47) -- 
(15*9/28, 15*6/47) -- (15*9/28, 15*4/31) -- 
(15*26/81, 15*4/31) -- (15*26/81, 15*3/23) -- 
(15*17/53, 15*3/23) -- (15*17/53, 15*8/61) -- 
(15*25/78, 15*8/61) -- (15*25/78, 15*12/91) -- 
(15*8/25, 15*12/91) -- (15*8/25, 15*1/7) -- 
(15*6/19, 15*1/7) -- (15*6/19, 15*2/13) -- 
(15*81/259, 15*2/13);

\draw 
(15*2/13, 15*81/259) 
-- (15*2/13, 15*6/19) -- (15*1/7, 15*6/19) 
-- (15*1/7, 15*8/25) -- (15*12/91, 15*8/25) 
-- (15*12/91, 15*25/78) -- (15*8/61, 15*25/78) 
-- (15*8/61, 15*17/53) -- (15*3/23, 15*17/53) 
-- (15*3/23, 15*26/81) -- (15*4/31, 15*26/81) 
-- (15*4/31, 15*9/28) -- (15*6/47, 15*9/28) 
-- (15*6/47, 15*28/87) -- (15*8/63, 15*28/87) 
-- (15*8/63, 15*19/59) -- (15*12/95, 15*19/59) 
-- (15*12/95, 15*29/90) -- (15*1/8, 15*29/90) 
-- (15*1/8, 15*17/52) -- (15*1/9, 15*17/52) 
-- (15*1/9, 15*18/55) -- (15*1/10, 15*18/55) 
-- (15*1/10, 15*32/97)-- (15*.075111, 15*32/97);
% -- (15*1/35, 15*32/97) ;

\draw[thick,dotted] (0,15/3)--(15*6/25,15/3);
\draw[thick,dotted] (15/3,15*6/25) -- (15/3,0);

\tikzstyle{every node}=[]
\node at (5.1,5.1) { \large{$\phi\ge 1/3$}};
\node at (2.7,2.7) { \large{$\phi< 1/3$}};
\node at (15*2/13+.65, 15*6/19+.05) {\begin{scriptsize}$(2/13,6/19)$\end{scriptsize}};
%\draw[domain=.302809:1/3,smooth,variable=\y,thick, dotted] plot ({15*((-sqrt(2/3)*\y+sqrt(2*\y*\y/3-(1/(1-\y)+\y)*(2*\y/3-1)))/(1/(1-\y)+\y)) *((-sqrt(2/3)*\y+sqrt(2*\y*\y/3-(1/(1-\y)+\y)*(2*\y/3-1)))/(1/(1-\y)+\y)) },{15*\y});
%\draw[domain=.302809:1/3,smooth,variable=\y,thick, dotted] plot ({15*\y},{15*((-sqrt(2/3)*\y+sqrt(2*\y*\y/3-(1/(1-\y)+\y)*(2*\y/3-1)))/(1/(1-\y)+\y)) *((-sqrt(2/3)*\y+sqrt(2*\y*\y/3-(1/(1-\y)+\y)*(2*\y/3-1)))/(1/(1-\y)+\y)) });
%\draw[domain=2/11:5/17,smooth,variable=\x] plot ({15*\x},{15*(1-3*\x)/(3-8*\x)});
%\draw[domain=0:1/35,smooth,variable=\x] plot ({15*\x},{15*(1-3*\x)/(3-8*\x)});
%\draw[domain=32/97:1/3,smooth,variable=\x] plot ({15*\x},{15*(1-3*\x)/(3-8*\x)});
%\draw[domain=0:1/3,smooth,variable=\x] plot ({15*\x}, {15*(1-\x)*(1-\x)/3});

%\draw[domain=.29914:1/3,smooth,variable=\x,thick, dotted] plot ({15*\x}, {15*(1-\x)*(1-\x)*(1-\x)/(1+2*\x-5*\x*\x)});
%\draw[domain=.29914:1/3,smooth,variable=\y, thick,dotted] plot ({15*(1-\y)*(1-\y)*(1-\y)/(1+2*\y-5*\y*\y)},{15*\y});

%\draw[domain=1/4:1/3,smooth,variable=\x,] plot ({15*\x}, {15*\x/(2-3*\x)});
%\draw[domain=0.05:0.19142,smooth,variable=\y] plot ({15*\y/(1-2/sqrt((2*\y*\y-\y)/(1-\y)))},{15*\y});

\draw[domain=1/4:.28231,smooth,variable=\x,thick, dotted] plot ({15*\x}, {15*(1-\x)*(1-\x)/((2*\x+1)*(1+2*\x-5*\x*\x))});
\draw[domain=.31994:1/3,smooth,variable=\x,thick, dotted] plot ({15*\x}, {15*(1-\x)*(1-\x)/((2*\x+1)*(1+2*\x-5*\x*\x))});
%\draw[domain=0:1/3,smooth,variable=\y, thick] plot ({15*\y/(1-2*\y/sqrt(1-\y))},{15*\y});
\draw[domain=6/19:.31994,smooth,variable=\y, thick, dotted] plot ({15/4},{15*\y});

\draw[domain=1/4:.28231,smooth,variable=\y,thick,dotted] plot ({15*(1-\y)*(1-\y)/((2*\y+1)*(1+2*\y-5*\y*\y))},{15*\y});
%\draw[domain=0:.16170,smooth,variable=\x, thick] plot ({15*\x},{15*\x/(1-2*\x/sqrt(1-\x))});
\draw[domain=.31992:1/3,smooth,variable=\y,thick,dotted] plot ({15*(1-\y)*(1-\y)/((2*\y+1)*(1+2*\y-5*\y*\y))},{15*\y});
\draw[domain=6/19:.31994,smooth,variable=\x, thick, dotted] plot ({15*\x},{15/4});

%\draw[domain=1/4:1/3,smooth,variable=\x,] plot ({15*\x}, {15*(1-\x)*(1-\x)/(3+6*\x-15*\x*\x)});
%\draw[domain=1/6:1/3,smooth,variable=\x] plot ({15*\x}, {15/(1+1/((1-\x)*(1-\x)))});

\draw[domain=0:.075111,smooth,variable=\x] plot ({15*\x},{15*(1-2*\x)*(1-2*\x)/(3-12*\x+16*\x*\x)});
\draw[domain=2/13:1/4,smooth,variable=\x] plot ({15*\x},{15*(1-2*\x)*(1-2*\x)/(3-12*\x+16*\x*\x)});
\draw[domain=0:.075111,smooth,variable=\y] plot ({15*(1-2*\y)*(1-2*\y)/(3-12*\y+16*\y*\y)},{15*\y});
\draw[domain=2/13:1/4,smooth,variable=\y] plot ({15*(1-2*\y)*(1-2*\y)/(3-12*\y+16*\y*\y)},{15*\y});

\node at (-.7, 5) {\begin{scriptsize}$y=1/3$\end{scriptsize}};
\node at (4,4) {$?$};

\node (t1) at (15/8, 5.6) {\begin{footnotesize}$\ref{maxbound}$\end{footnotesize}};
\draw[
    gray, ultra thin, decoration={markings,mark=at position 1 with {\arrow[black,scale=2]{>}}},
    postaction={decorate},
    ]
(t1) to (15/8,5);

\node (t3) at (5.1, 4.6) {\begin{footnotesize}$\ref{phi13monster}$\end{footnotesize}};
\draw[
    gray, ultra thin, decoration={markings,mark=at position 1 with {\arrow[black,scale=2]{>}}},
    postaction={decorate},
    ]
(t3) to (15*.267,15*.3);

\node (u1) at (2, 15/4) {\begin{footnotesize}$\ref{phi13bettercurve}$\end{footnotesize}};
\draw[
    gray, ultra thin, decoration={markings,mark=at position 1 with {\arrow[black,scale=2]{>}}},
    postaction={decorate},
    ]
(u1) to [bend right=20] ({15/5},{15*2/7});
\draw[
    gray, ultra thin, decoration={markings,mark=at position 1 with {\arrow[black,scale=2]{>}}},
    postaction={decorate},
    ]
(u1) to [bend left=20] ({15/13},{15*121/367});
\end{tikzpicture}

\caption{When $\psi(x,y)<1/3$ and when $\phi(x,y)<1/3$.} \label{fig:phi13}
\end{figure}
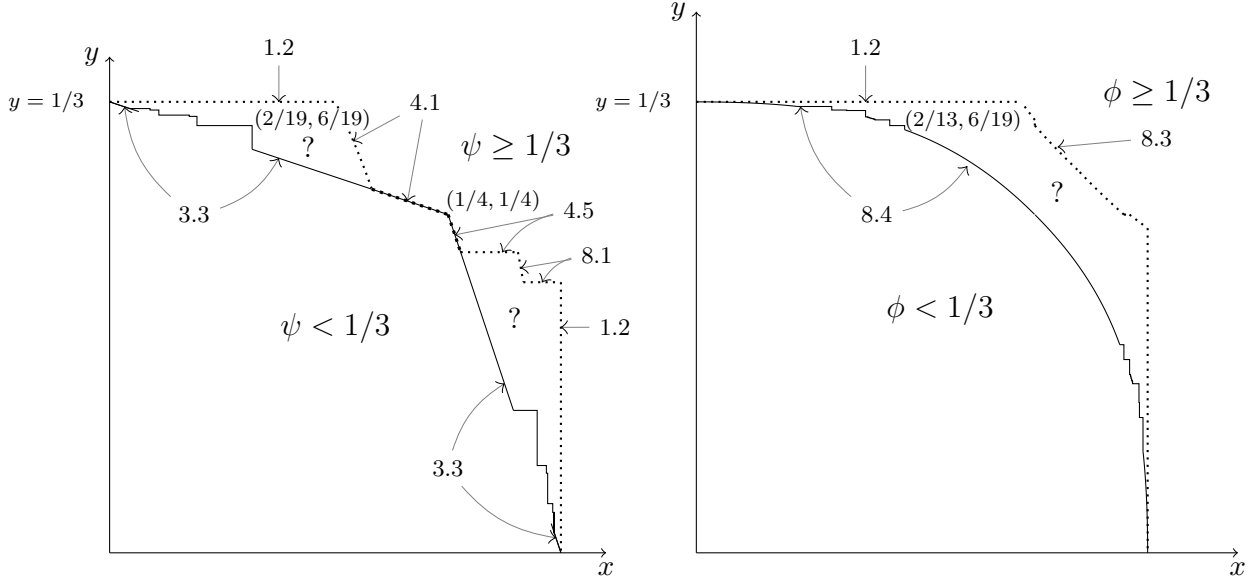

\begin{thm}\label{psi13good}
Let $x,y\in (0,1]$ with $y>\frac15$ and $3x+\frac{y}{3(1-y)}\ge 1$. Then $\psi(x,y)\ge \frac13$.
\end{thm}
\Proof
We may assume that $y\le 1/3$, by \ref{maxbound}, and so $y/(3(1-y))\le 1/6$. Consequently $x\ge 5/18$, and in particular
$x>2y/3$ (we will need this later). Also, since $1/5\le y\le 1/3$, it follows that $3y-y/(3(1-y))> 1/2$; and so 
$$\left(3x+\frac{y}{3(1-y)}\right) + \left(3y-\frac{y}{3(1-y)}\right)> \frac32,$$
and consequently $x+y>1/2$.
Let $G$ be $(x,y)$-biconstrained via $(A,B,C)$, and suppose that $|N^2_A(v)|<|A|/3$ for each $v\in C$. 
Let $H$ be the subgraph induced on $B\cup C$, and let $H_1\ll H_k$
be its components. Let $B_i = V(H_i)\cap B$ and $C_i = V(H_i)\cap C$, and
$b_i = |B_i|/|B|$, $c_i = |C_i|/|C|$, for $1\le i\le k$. Since $y>0$, $B_i,C_i$ are both nonempty
and so $b_i, c_i\ge y$ for $1\le i\le k$.
%y>1/5
For $1\le i\le k$, let $A_i$ be the set of vertices in $A$ with a neighbour in $B_i$, and let $A_i^*$ be the set of vertices in $A$
such that $N(v)\subseteq B_i$. 
\\
\\
(1) {\em $k\ge 2$.}
\\
\\
Suppose that $k=1$, and let $H'$ be the graph with vertex set $B$ in which $u,u'$ are adjacent if $u,u'$ have
a common neighbour in $H$. Then every stable set of $H'$ has cardinality at most $4$. By~\ref{radius}
there is a vertex $u_1\in B$ with $H'$-distance at most
four to every other
vertex in $B$; and so the $H$-distance from $u_1$ to each vertex in $B$ is at most eight. Let $v_1\in C$ be adjacent to $u_1$.
Let $A'=A\setminus N^2_A(v_1)$ and $B'=B\setminus N(v_1)$. Hence $|A'|>2|A|/3$.
Since every vertex in
$A'$ has at least $x|B|$ neighbours in $B'$, and $|B'|\le (1-y)|B|$, some vertex $u\in B'$ has at least
$$\frac{x|A'|}{1-y}\ge \frac{2x|A|}{3(1-y)}$$
neighbours in $A'$. Choose a path of $H$ between $u_1$ and $u$ of length at most eight, and let its vertices
be $u_1\d v_2\d u_2\c v_t\d u_t=u$
say, in order. Thus $t\le 5$,
and so there exists $i$ with $1\le i\le t-1$ such that there are at least $|N_{A'}(u)|/4$ vertices that belong to
$N_{A}(u_{i+1})\setminus N_{A}(u_{i})$. Since $|N_{A}(u_{i})|\ge x|A|$, it follows that
$$|N^2_A(v_{i+1})|\ge x|A|+|N_{A'}(u)|/4\ge \frac{x+2x}{12(1-y)}|A|\ge |A|/3,$$
%3x+x/(2(1-y))< 1
a contradiction, since $x\ge 2y/3$ and so $x+x/(6(1-y))\ge x+y/(9(1-y))\ge 1/3$ . This proves (1).
\\
\\
(2) {\em $b_i\le 1-x-y<1/2$ for $1\le i\le k$, and so $k\ge 3$.}
%x+y>1/2
\\
\\
Suppose that $b_1>1-x-y$ say.
Thus, if $u\in A\setminus A_1$, then $u\in N^2_A(v)$ for every $v\in C\setminus C_1$; and so $|A\setminus A_1|<|A|/3$, and
so $|A_1|>2|A|/3$.
Let $H'$ be the graph with vertex set $C_1$ in which $v,v'$ are adjacent if they have a common $H_1$-neighbour in $B_1$.
Thus $H'$ has stability number at most three (by (1)) and so has radius at most three, by~\ref{radius}. Choose
$v_1\in C_1$ such that every vertex in $C_1$ has $H_1$-distance at most six from $v_1$. Let $A' = A_1\setminus N^2_A(v_1)$;
thus $|A'|>|A|/3$. Since every vertex in
$A'$ has a neighbour in $B_1$ and hence has at least $y|C|$ second neighbours in $C_1$, there exists $v\in C_1$
such that
$$|N^2_{A'}(v)|\ge \frac{y}{|C_1|}|A'|\ge \frac{y}{3(1-y)}|A|,$$
since $|C_1|\le (1-y)|C|$. Choose a path of $H_1$ between $v_1,v$ of length at most six, with vertices
$v_1\d u_1\d v_2\c u_{t-1}\d v_t=v$
say where $t\le 4$. Then for some $i$ with $1\le i\le t-1$,
$$|N^2_{A'}(v_{i+1})\setminus N^2_{A'}(v_i)|\ge \frac{y}{9(1-y)}|A|,$$
and hence
$$|N^2_{A'}(v_{i+1})|\ge \left(x+\frac{y}{9(1-y)}\right)|A|$$
since all vertices of $N_A(u_i)$ belong to $N^2_{A}(v_{i+1})$ and do not belong to $N^2_{A'}(v_{i+1})\setminus N^2_{A'}(v_i)$.
But $3x+y/(3(1-y))\ge 1$, a contradiction. This proves (2).
%$3x+y/(3(1-y))\ge 1$

\bigskip

By (2), $k\ge 3$; and $k\le 4$ since $y>1/5$. We may assume that $|B_1|,|B_2|\ge |B_i|$ for $i\ge 3$; let $B_0=\bigcup_{3\le i\le k}B_i$,
and $C_0=\bigcup_{3\le i\le k}C_i$. Hence $|B_0|\le |B|/2$ since $k\le 4$. Let $b_0=|B_0|/|B|$ and $c_0=|C_0|/|C|$; let
$A_0$ be the set of vertices in $A$ with a neighbour in $B_0$, and let $A_0^*$ be the set of vertices in $A$
such that $N(v)\subseteq B_0$. For $0\le i<j\le 2$, choose $A_{ij}=A_{ji}\subseteq A_i\cap A_j$ such that the sets
$A_{12}, A_{13}, A_{23}, A_0^*, A_1^*, A_2^*$ are pairwise disjoint and have union $A$.
For $0\le i\le 2$ let $a_i=|A_i|/|A|$ and $a_i^* = |A_i^*|/|A|$, and for $0\le i,j\le 2$ with $i\ne j$ let $a_{ij}= |A_{ij}|/|A|$.
Since $b_1,b_2\le 1-x-y< x+y$ and $b_0\le 1/2< x+y$, we have $b_i<x+y$
for $i = 0,1,2$. Let $0\le i\le 2$, and choose $v\in C_i$ uniformly at random. Then $A_i^*\subseteq N^2_A(v)$ because $b_i<x+y$, and the
expected value of $|N^2_A(v)\cap (A_i\setminus A_i^*)|$ is at least $(y/c_i)|A_i\setminus A_i^*|$; so the expected value of
$|N^2_A(v)|$ is at least 
$$|A_i^*| + \frac{y}{c_i} |A_i|\setminus |A_i^*|= \left(a^*_i + \frac{y}{c_i} (a_i-a^*_i)\right)|A|.$$
Since $|N^2_A(v)|<|A|/3$, it follows that
$a^*_i + (y/c_i) (a_i-a^*_i)<1/3$. Now $A^*_0, A_{01}, A_{02}$ are pairwise disjoint subsets of $A_0$, so
$a_{01}+a_{02}\le a_0-a_0^*$; and hence
$$a^*_0 + (y/c_0) (a_{01}+a_{02})\le a^*_0 + (y/c_0) (a_0-a^*_0)<1/3.$$
Similarly we have $a^*_1 + (y/c_1) (a_{01}+a_{12}) <1/3$ and $a^*_2 + (y/c_2) (a_{02}+a_{12}) <1/3$; and by summing these three 
inequalities and using the equation 
$$a_1^*+a_2^*+a_3^*+a_{12}+a_{13}+a_{23}= 1$$
we obtain
$$a_{12}\left(\frac{y}{c_1}+\frac{y}{c_2}-1\right) + a_{13}\left(\frac{y}{c_1}+\frac{y}{c_3}-1\right)+a_{23}\left(\frac{y}{c_2}+\frac{y}{c_3}-1\right)<0.$$
Consequently there exist distinct $i,j\in \{0,1,2\}$ with $y/c_i+y/c_j-1<0$. But
$1/c_i+1/c_j\ge 4/(c_i+c_j)$, and $c_i+c_j\le 1-y$, 
and so $4y/(1-y)<1$, a contradiction.  
This proves \ref{psi13good}.~\bbox

For $\phi$, we need the following. 
\begin{thm}\label{leftbound}
Let $x,y\in (0,1]$. Let $G$ be $(x,2/3)$-constrained via $(A,B,C)$, such that every three vertices in $B$ have a
common neighbour in $C$, and every vertex $w\in C$
satisfies $|N^2_A(w)|< (1-y)|A|$. If $v_1\in B$ has $a|A|$ neighbours in $A$ then
$$a\le 1- (1+2x-5x^2)y/(1-x)^2.$$
\end{thm}
\Proof
Let $v_1\in B$ have $a|A|$ neighbours in $A$. Define $A_1=N_A(v_1)$ and choose $C_1\subseteq N_C(v_1)$ with $|C_1|=2|C|/3$
(we may assume this is an integer). Let $A_1'=A\setminus A_1$
and $C_1'=C\setminus C_1$.
\\
\\
(1) {\em If some vertex $v_2\in B$ has a set $A_2$ of $t|B|$ neighbours in $A_1'$, then
\begin{itemize}
\item every $v\in B$ has at most $(1-y-a-t)|A|$ neighbours in $A\setminus (A_1\cup A_2)$; and
\item the sum over all $v\in B$ of the number of neighbours of $v$ in $A\setminus (A_1\cup A_2)$
is $(1-a-t)x|A||B|$.
\end{itemize}}
\noindent The first claim follows since $v_1,v_2,v$ have a common neighbour in $C$. The second holds since every vertex in
$A\setminus (A_1\cup A_2)$ has $x|B|$ neighbours. The third follows. This proves (1).

\bigskip

The sum over $u\in A_1'$, of $|N^2_{C_1}(u)|$, is at most $\frac23 (1-y-a)|A||C|$; and for each $u$,
$|N^2_{C_1}(u)|\ge \max_{v\in N_B(u)}|N_{C_1}(v)|$. But the latter is at least
$$ \sum_{v\in N_B(u)}|N_{C_1}(v)|/(x|B|).$$
It follows that
$$\sum_{u\in A_1'}\sum_{v\in N_B(u)}|N_{C_1}(v)|   \le \frac23 x(1-y-a)|A||B||C|.$$
Consequently
$$\sum_{v\in B}|N_{A_1'}(v)||N_{C_1}(v)| \le \frac23 x(1-y-a)|A||B||C|.$$

Moreover, each vertex in $C_1'$ has at least $x|B|$ nonneighbours in $B$, and so there are at most $(1-x)/3|B||C|$
edges between $B$ and $C_1'$. Hence there are at least $(1+x)/3 |B||C|$ edges between $B$ and $C_1$.

For each $v\in B$, let $p(v)=|N_{A_1'}(v)|/|A|$.
Thus $\sum_{v\in B}p(v)=x(1-a)|B|$. By setting $q(v)=3|N_{C_1}(v)|/|C|-1$ we deduce:
for each $v\in B$ there exists $q(v)$ such that
\begin{itemize}
\item for each $v\in B$, $1/3\le q(v)/3+1/3\le 2/3$, that is, $0\le q(v)\le 1$;
\item $\sum{v\in B}(q(v)/3+1/3)\ge (1+x)/3 |B|$, that is, $\sum{v\in B}q(v)\ge x |B|$;
\item $\sum_{v\in B}p(v)(q(v)/3+1/3) \le \frac23 x(1-y-a)|B|$, that is, $\sum_{v\in B}p(v)q(v) \le (x-2xy-xa)|B|$.
\end{itemize}

Let $Q\subseteq B$ be the $x|B|$ vertices in $B$ (we may assume this is an integer) with $p(v)$ smallest.
Then the expression in the last bullet above is minimized by setting $q(v)=1$ for $v\in Q$, and $q(v)=0$ for $v\in B\setminus Q$.
Consequently $\sum_{v\in Q}p(v) \le (x-xa- 2xy)|B|$.

Choose $v_2\in B\setminus Q$ with $|N_{A_1'}(v_2)|$ maximum; $A_2$ say, where $|A_2|=t|A|$. By (1),
every $v\in B$ has at most $(1-y-a-t)|A|$ neighbours in $A_1'\setminus A_2$, and the sum over all $v\in B$ of the number
of neighbours of $v$ in $A_1'\setminus A_2$
is $(1-a-t)x|A||B|$. So the number of edges between $A_1'\setminus A_2$ and $B\setminus Q$ is at most
$(1-y-a-t)(1-x)|A||B|$; and the number between $A_1'\setminus A_2$ and $Q$ is at most $(x-xa- 2xy)|A||B|$, since
$\sum_{v\in Q}p(v) \le (x-xa- 2xy)|B|$. Hence
the number between $A_1'\setminus A_2$ and $B$ is at most $((1-y-a-t)(1-x) +x-xa- 2xy)|A||B|$, and since this number
equals $(1-a-t)x|A||B|$, it follows that
$$(1-y-a-t)(1-x) +x-xa- 2xy\ge (1-a-t)x,$$
that is,
$$1-y-a-t-x+xa- xy+2tx \ge 0.$$
Consequently $t\le (1-x-y-a+xa- xy)/(1-2x)$.

Now $t(|B|-|Q|)|A|+ (x-xa- 2xy)|A||B| \ge x(1-a)|A||B|$, so
$$t(1-x) \ge 2xy.$$
Hence
$$(1-x-y-a+xa- xy)/(1-2x)\ge t\ge 2xy/(1-x),$$
that is,
$$(1-x)(1-x-y-a+xa- xy)\ge (1-2x)2xy.$$
Consequently
$$(1-x)^2(1-a)\ge (1+2x-5x^2)y.$$
This proves \ref{leftbound}.~\bbox

\begin{thm}\label{phi13monster}
Let $x,y\in (0,1]$ with $(1-x)^2> (2x+1)(1+2x-5x^2)y$ 
and $x\ge 1/4$ and $4x^2y^2\ge (1-y)(x-y)^2$. Then $\phi(x,y)\ge 1/3$.
Consequently if $x>0.28231$ then $\phi(x,x)\ge 1/3$.
\end{thm}
\Proof
Let $\phi(x,y)=z$ and suppose that $z<1/3$. Then there is a graph $G$ that is $(x,1-z)$-constrained via $A,B,C$, such that
$|N^2_A(w)|\le (1-y)|A|$ for each $w\in C$, by \ref{rotate}
and \ref{permute}. As in \ref{hompe4}, there exists $w\in C$ such that there are at least $x(1-z)|A|\cdot|B|$ edges between 
$N_B(w)$ and $N^2_A(w)$. Define $B_1=N_B(w)$ and $B_2=B\setminus B_1$; and let $B_2=t|B|$. For each $u\in A$, 
let $u$ have $d(u)|B|$ neighbours in $B_1$. 
Let $A_1=N^2_A(w)$ and $A_2=A\setminus A_1$; and let $|A_2|=s|A|$. Thus $d(u)=0$ for each $u\in A_2$.
\\
\\
(1) {\em $t\ge x$, and $s\ge y$, and $0\le d(u)\le \min(x,1-t)$ for each $u\in A$. Also
$$x(1-y)|A|\ge \sum_{u\in A_1}d(v)\ge x(1-z)|A|.$$}
\noindent We may assume that every vertex in $A$ has degree exactly $x|B|$;
so $0\le d(u)\le \min(x,1-t)$ for each $u\in A$. Since $|A_1|\le (1-y)|A|$, it follows that $s\ge y$. In particular, $A_2\ne\emptyset$, and so some vertex in $A_2$
has $x|B|$ neighbours in $B_2$, and so $t\ge x$. Since there are at least $x(1-z)|A|\cdot|B|$ edges between 
$N_B(w)$ and $N^2_A(w)$, it follows that $\sum_{u\in A_1}d(v)\ge x(1-z)|A|$. Since $|A_1|\le (1-y)|A|$
and every vertex in $A_1$ has degree exactly $x|B|$, it follows that the number of edges between $A_1$ and $B_1$
is at most $x(1-y)|A|\cdot |B|$, and so $x(1-y)|A|\ge \sum_{u\in A_1}d(v)$.  
This proves (1).
\\
\\
(2) {\em $1-t\ge \frac{2x(1-x)^2/3}{(1-x)^2-(1+2x-5x^2)y}$. Consequently $x<1-3t/2$ and so $x<1-t$.}
\\
\\
There are at least $2x|A|\cdot |B|/3$ edges between $B_1$ and $A$, since $z<1/3$. But each vertex in $B_1$ has at most
$$\left(1-\frac{(1+2x-5x^2)y}{(1-x)^2}\right)|A|$$
neighbours in $A$, by \ref{leftbound}, and the first claim follows. 
To show that $x<1-3t/2$, suppose not; then
$$2x/3+1/3>1-t  \ge \frac{2x(1-x)^2/3}{(1-x)^2-(1+2x-5x^2)y}$$
and so 
$$(2x+1)((1-x)^2-(1+2x-5x^2)y)> 2x(1-x)^2,$$
that is, 
$$(1-x)^2> (2x+1)(1+2x-5x^2)y$$
contrary to the hypothesis.
This proves (2). 

\bigskip

Let us choose $v_1,v_1'\in A_1$ uniformly and independently at random, and choose $v_2\in A_2$ uniformly at random. Then 
for $u\in A$, the probability that all of $v_1,v_1', v_2$ are nonadjacent to $u$ is
$$\frac{t-x+d(v)}{t}\left(\frac{1-t-d(v)}{1-t}\right)^2.$$
Since $1-y>2/3$ and so $v_1,v_1',v_2$ have a common neighbour in $C$, say $w'$, and $|N^2_A(w')|\le (1-y)|A|$, it follows
that 
$$\sum_{u\in A}\left(1-\frac{t-x+d(v)}{t}\left(\frac{1-t-d(v)}{1-t}\right)^2\right)\le (1-y)|A|,$$
that is,
$$\sum_{u\in A}\frac{t+d(u)-x}{t}\left(\frac{t+d(u)-1}{1-t}\right)^2\ge y|A|.$$
This can be rewritten as:
$$\sum_{u\in A}f(d(v))\ge t(1-t)^2(1-y)|A|,$$
where $f(r)$ is the polynomial $(r+t-x)(r+t-1)^2$.
We therefore need to investigate the maximum value of $\sum_{u\in A}f(d(v))$ (which we call ``the objective function'')
over all choices of the numbers $d(u)(u\in A)$ satisfying the various constraints, and verify that this maximum is
less than $t(1-t)^2(1-y)|A|$.

The derivative of $f(r)$ is zero when $3r^2+2(3t-x-2)r+(3t-2x-1)(t-1)=0$, which has roots 
$r=1-t$ and $r=(2x+1)/3-t$. Let us define $r_0=(2x+1)/3-t$. Since $r_0<1-t$, the function $f(r)$ increases for $r<r_0$
and for $r>1-t$, and decreases for $r_0<r<1-t$. 

The second derivative of $f(r)$ is zero when $3r+3t-x-2=0$, that is, when $r=r_1$ where $r_1=2/3+x/3-t$. 
By (2), $x<r_1$, and we are only concerned $f(r)$ for $r$ in with the range $0\le r\le x$; so in particular 
all such $r$ are less than $r_1$. The function $f(r)$ is concave through the range $0\le r\le r_1$, since its second derivative
is at most zero.

Let us choose real numbers $d(v)(v\in A)$ satisfying the constraints 
\begin{itemize}
\item $0\le d(u)\le x$  for each $u\in A$;
\item $d(u)=0$ for at least $y|A|$ vertices $u\in A$;
\item $x(1-y)|A|\ge \sum_{u\in A_1}d(v)\ge 2x|A|/3$
\end{itemize}
to maximize the function $\sum_{u\in A}f(d(v))$. From the concavity of $f$, it follows that 
there exists $r^*$ with $0<r^*\le x$ such that $d(v)\in \{0,r^*\}$ for all $v$ (because if there were $u,v$ with $d(u),d(v)$
distinct and nonzero, replacing them both by $(d(u)+d(v))/2$ would still satisfy the constraints and increase the 
objective function). Similarly, if there were more than $y|A|$ vertices $v$ with $d(v)=0$, then choose some one of them, $v$
say, and choose some $u$ with $d(u)>0$; then again replacing them both by $(d(u)+d(v))/2$ would still satisfy the constraints and increase the
objective function. We deduce that there are exactly $y|A|$ vertices $v$ with $d(v)=0$.

Now the problem breaks into three cases, depending which of the constraints $x(1-y)|A|\ge \sum_{u\in A_1}d(v)\ge 2x|A|/3$
hold with equality. 

Suppose first that neither holds with equality. Then from the optimality of the objective function, it follows that 
$r^*=r_0$, and since 
$$x(1-y)|A|\ge \sum_{u\in A_1}d(v)\ge x(1-z)|A|,$$ it follows that
$$x(1-y)\ge (1-y)r_0\ge x(1-z),$$
that is, 
$$x\ge (2x+1)/3-t \ge 2x/(3-3y).$$ Thus, if $t$ satisfies 
$$(1-x)/3\le t\le (2x+1)/3-2x/(3-3y)$$ then there is a possible 
optimal solution where the objective function has value
$$y|A|f(0)+(1-y)|A|f(r_0).$$ 
Now $f(0)=(t-x)(t-1)^2$, and 
$$f(r_0)=(r_0+t-x)(r_0+t-1)^2=((1-x)/3)((2x-2)/3)^2=4(1-x)^3/27.$$
We must therefore check that for $t$ in the given range, 
$$y|A|(t-x)(1-t)^2+4(1-y)|A|(1-x)^3/27 < t(1-t)^2(1-y)|A|,$$
This simplifies to:
$$4(1-x)^3(1-y)<27(1-t)^2(t-2ty+xy).$$
Now the function $27(1-t)^2(t-2ty+xy)$ has no local minimum at $t$ with $t<1$, and so is minimized at one of 
the ends of the range. Since $t\ge x$, we might as well replace the lower extreme of the range by $t\ge x$ (because 
it makes the arithmetic easier); so to check the lower extreme, we need to check that
$$4(1-x)^3(1-y)<27x(1-x)^2(1-y),$$
that is,
$4(1-x)<27x,$ which is true by hypothesis.

For the upper extreme, $t\le (2x+1)/3-2x/(3-3y)<1/3$; so it suffices to check that
$$4(1-x)^3(1-y)<27(1-t)^2(t-2ty+xy)$$
when $t=1/3$, that is, to check
$(1-x)^3(1-y)<(1-2y+3xy)$. But $(1-x)^3<1/2$ by hypothesis, and $(1-2y+3xy)/(1-y)\ge (1-2y)/(1-y)\ge 1/2$ since $y<1/3$.
This finished the first of the three cases.

Now let us assume that $x(1-y)|A|= \sum_{u\in A_1}d(v)$. It follows from the optimality of the objective function that
$r^*\le r_0$. Moreover, since $x(1-y)|A|= \sum_{u\in A_1}d(v)$, it follows that $x(1-y)|A|= (1-y)|A|r^*$, so
$r^*=x$. This is only possible if $x\le r_0$, that is, $t<(1-x)/3$; and this is impossible since $t\ge x\ge 1/4$.
This finishes the second case.

Finally, we assume that $\sum_{u\in A_1}d(v) = 2x|A|/3$. It follows from the optimality of the objective function that
$r^*\ge r_0$. Moreover, since $\sum_{u\in A_1}d(v)=2x|A|/3$, it follows that $(1-y)|A|r^*=2x|A|/3$, that is,
$r^*=2x/(3-3y)$. 
We must check that
$$y(t-x)(t-1)^2+(1-y)(r^*+t-x)(r^*+t-1)^2< t(1-t)^2(1-y).$$
This is cubic in $t$, and, collecting the various powers of $t$, it becomes:
$$yt^3 +t^2(-xy-2y + (1-y)(r^*-x) + 2(1-y)(r^*-1)+ 2(1-y))$$
$$+t(y+2xy + (1-y)(r^*-1)^2+2(1-y)(r^*-x)(r^*-1) -(1-y))$$
$$+(-xy + (1-y)(r^*-x)(r^*-1)^2)<0.$$
This simplifies to:
$$yt^3 +(x-2y)t^2 + t(y -2x/3 +4x^2y/(3-3y)) -xy+x(3y-1)(2x/3-1+y)^2/(3(1-y)^2)<0. $$
The derivative of the left side with respect to $t$ is
$$3yt^2 +2(x-2y)t + y -2x/3 +4x^2y/(3-3y),$$ 
which can be rewritten as
$$3y(t+(x-2y)/(3y))^2  -(x-y)^2/(3y) +4x^2y/(3-3y).$$
Since by hypothesis, $-(x-y)^2/(3y) +4x^2y/(3-3y)\ge 0$, the derivative is nonnegative, at every value of $t$.
Thus we only need verify the inequality for the maximum value of $t$ that lies in the range.

By (2), $t\le 2(1-x)/3$; so it is enough to verify that
$$y(t-x)(t-1)^2+(1-y)(r^*+t-x)(r^*+t-1)^2< t(1-t)^2(1-y)$$
holds when $t=2(1-x)/3$. Thus we need to check that
$$y(2(1-x)/3)^3 +(x-2y)(2(1-x)/3)^2 + (2(1-x)/3)(y -2x/3 +4x^2y/(3-3y))$$
$$ -xy+x(3y-1)(2x/3-1+y)^2/(3(1-y)^2)<0. $$
According to WolframAlpha, this is true for all $x,y$ with $1/4<x<1/3$ and $0<y\le x$. This proves \ref{phi13monster}.~\bbox

\begin{thm}\label{phi13bettercurve}
Let $x,y\in (0,1]$ with $x\le \frac14$ and $y< \frac13$ and $y< \frac{(1-2x)^2}{3-12x+16x^2}$; then $\phi(x,y)<\frac13$.
\end{thm}
\Proof
Apply \ref{phiaddpath} to \ref{phi12bettercurve}. This proves \ref{phi13bettercurve}.~\bbox

\section{The 3/4 Level}
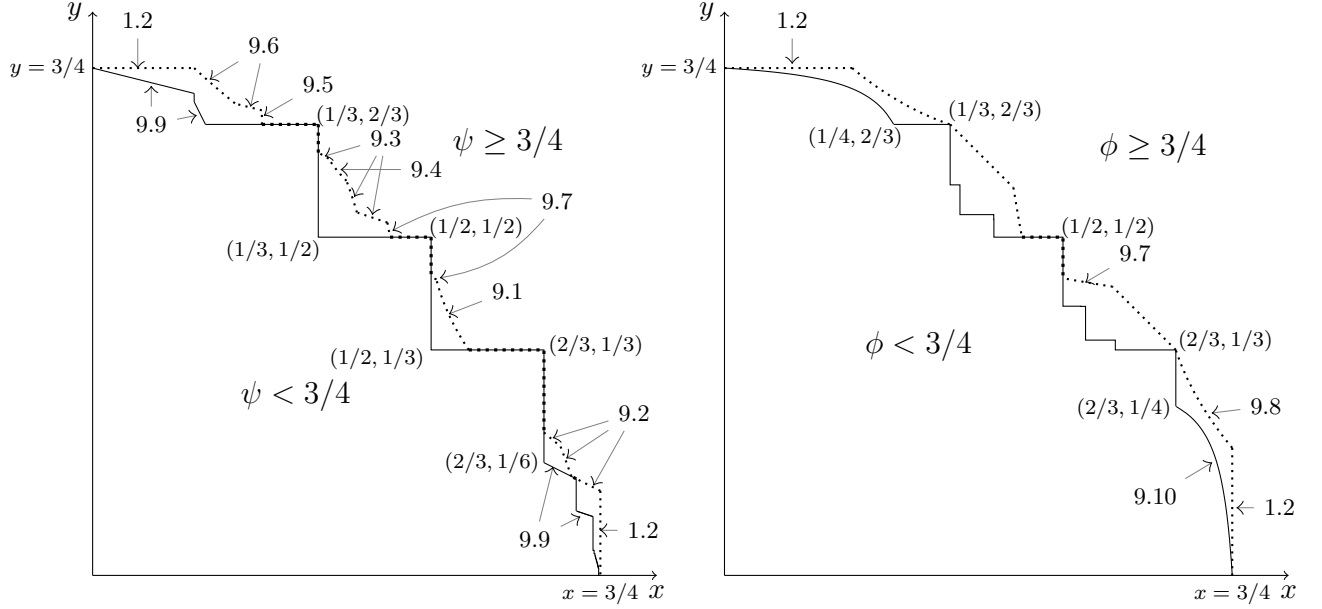
\begin{figure}[ht]
\centering

\begin{tikzpicture}[scale=.6,auto=left]

\begin{scope}[shift ={(-14,0)}]
\draw[->] (0,0) -- (12.5,0) node[anchor=north]{$x$};
\draw[->] (0,0) -- (0,12.5) node[anchor=east] {$y$};

\node at (9.2,9.6) { \large{$\psi\ge 3/4$}};
\node at (4.5,4) { \large{$\psi< 3/4$}};
\node at (8.5,7.7) {\begin{scriptsize}$(1/2,1/2)$\end{scriptsize}};
\node at (4,7.2) {\begin{scriptsize}$(1/3,1/2)$\end{scriptsize}};
\node at (6.3,4.8) {\begin{scriptsize}$(1/2,1/3)$\end{scriptsize}};
\node at (11.15, 5.1) {\begin{scriptsize}$(2/3,1/3)$\end{scriptsize}};
\node at (6, 10.2) {\begin{scriptsize}$(1/3,2/3)$\end{scriptsize}};
\node at (8.85,2.5) {\begin{scriptsize}$(2/3,1/6)$\end{scriptsize}};
\node at (-1, 15*3/4) {\begin{scriptsize}$y=3/4$\end{scriptsize}};
\node at (15*3/4, -.4) {\begin{scriptsize}$x=3/4$\end{scriptsize}};

\draw (15*3/4, 15*0/1) --
(15*74/99, 15*0/1) -- (15*74/99, 15*1/99) --
(15*71/95, 15*1/99) -- (15*71/95, 15*1/95) --
(15*68/91, 15*1/95) -- (15*68/91, 15*1/91) --
(15*65/87, 15*1/91) -- (15*65/87, 15*1/87) --
(15*62/83, 15*1/87) -- (15*62/83, 15*1/83) --
(15*59/79, 15*1/83) -- (15*59/79, 15*1/79) --
(15*56/75, 15*1/79) -- (15*56/75, 15*1/75) --
(15*53/71, 15*1/75) -- (15*53/71, 15*1/71) --
(15*50/67, 15*1/71) -- (15*50/67, 15*1/67) --
(15*47/63, 15*1/67) -- (15*47/63, 15*1/63) --
(15*44/59, 15*1/63) -- (15*44/59, 15*1/59) --
(15*41/55, 15*1/59) -- (15*41/55, 15*1/55) --
(15*38/51, 15*1/55) -- (15*38/51, 15*1/51) --
(15*35/47, 15*1/51) -- (15*35/47, 15*1/47) --
(15*32/43, 15*1/47) -- (15*32/43, 15*1/43) --
(15*29/39, 15*1/43) -- (15*29/39, 15*1/39) --
(15*26/35, 15*1/39) -- (15*26/35, 15*1/35) --
(15*23/31, 15*1/35) -- (15*23/31, 15*1/31) --
(15*20/27, 15*1/31) -- (15*20/27, 15*1/27) --
(15*17/23, 15*1/27) -- (15*17/23, 15*2/23);
\draw[domain=2/23:2/21,smooth,variable=\y] plot ({15*(1-3*\y)},{15*\y});%3.9
\draw (15*5/7,15*2/21) -- (15*5/7, 15*1/7);
\draw[domain=1/7:1/6,smooth,variable=\y] plot ({15*(1-2*\y)},{15*\y});%3.9
\draw 
(15*3/20, 15*57/80)--(15*3/20, 15*7/10);
\draw(15*3/20, 15*7/10) -- (15*1/6, 15*2/3);%3.9
\draw (15*1/6, 15*2/3) --
(15*1/3, 15*2/3) -- (15*1/3,15*1/2) -- (15*1/2, 15*1/2)
-- (15*1/2, 15*1/3) -- (15*2/3, 15*1/3) -- (15*2/3, 15*1/6);
\draw[domain=57/80:3/4,smooth,variable=\y] plot ({15*(3-4*\y)},{15*\y});%3.9
\draw[domain=2/23:2/21,smooth,variable=\y] plot ({15*(1-3*\y)},{15*\y});%3.9

\draw[domain=.4355:.4393, smooth, variable=\x, thick, dotted] plot ({15*\x}, {15*(4*\x-3)*(4*\x-3)/(16*\x*\x)});%3/4from5.7
\draw[domain=.5215:.5365, smooth, variable=\y, thick, dotted] plot ({15*(2-3*\y)}, {15*(\y)});%3.3
\draw[domain=.5365:7/12, smooth, variable=\y, thick, dotted] plot ({15*(5-6*\y)/(11-12*\y)}, {15*\y});%3.3
\draw[domain=7/12:.6162, smooth, variable=\y, thick, dotted] plot ({15*(6*\y*\y-8*\y+3)}, {15*\y});%3.4
\draw[domain=.6162:.625, smooth, variable=\y, thick, dotted] plot({15*(3-4*\y)/(4-4*\y)}, {15*\y});%3.3
\draw[very thick, dotted] (15*1/3, 15*.625) -- (15*1/3, 15*2/3);

\draw[very thick, dotted] (15*1/2,15*1/2) -- (15*1/2, 15*.4393);
\draw[domain=.4379:.4393, smooth, variable=\y, thick, dotted] plot ({15*((4*\y-3)*(4*\y-3)/(16*\y*\y))}, {15*(\y)});%3/4from5.7

\draw[domain=1/3:.4379, smooth, variable=\y, thick, dotted] plot ({15*(1-2*\y+2*\y*\y)}, {15*(\y)});%3.1
\draw[very thick, dotted] (15*5/9, 15*1/3) -- (15*2/3, 15*1/3) -- (15*2/3, 15*.2113);
\draw[domain=.2113:.191, smooth, variable=\y, thick, dotted] plot ({15*(1-2*\y+2*\y*\y)}, {15*(\y)});%3.2
\draw[domain=.1464:.191, smooth, variable=\y, thick, dotted] plot ({15*(.25*(3-(\y)/(1-\y)))}, {15*(\y)});%3.2
%\draw[domain=.1962:.2113, smooth, variable=\y, thick, dotted] plot ({15*(1-\y-\y*(2/3-\y)/(1-\y))}, {15*(\y)});%3.2
\draw[domain=1/7:.1464, smooth, variable=\y, very thick, dotted] plot ({15*(1-2*\y)}, {15*(\y)});%3.2
\draw[domain=.125:1/7, smooth, variable=\y, thick, dotted] plot ({15*(1-2*\y)}, {15*(\y)});%3.2
\draw[thick, dotted] (15*.75, 15*.125) -- (15*.75, 0);
\draw[very thick, dotted] (15*.25, 15*2/3) -- (15*1/3, 15*2/3);
\draw[domain=.15:.2111,smooth,variable=\x,thick,dotted] plot ({15*(\x)},{15*((1/6)*((-1)*\x-sqrt(\x*\x+12*\x)+6))});%3.6
\draw[thick, dotted] (15*.25, 15*2/3) -- (15*.25, 15*.6875);
\draw[thick, dotted] (0, 15*.75) -- (15*.15, 15*.75);
\draw[very thick, dotted] (15*.4393, 15*.5) -- (15*.5, 15*.5);
\draw[domain=.6875:.6972,smooth,variable=\y,thick,dotted] plot ({15*(3-4*\y)},{15*\y});%3.6

\node (r1) at (1, 12.3) {\begin{footnotesize}$\ref{maxbound}$\end{footnotesize}};
\draw[
    gray, ultra thin, decoration={markings,mark=at position 1 with {\arrow[black,scale=2]{>}}},
    postaction={decorate},
    ]
(r1) to (1,15*3/4);

\node (r2) at (3.8, 15*3/4+.5) {\begin{footnotesize}$\ref{3/4ub5}$\end{footnotesize}};
\draw[ 
    gray, ultra thin, decoration={markings,mark=at position 1 with {\arrow[black,scale=2]{>}}}, 
    postaction={decorate}, ] 
    (r2) to (2.6,11.0);
\draw[
    gray, ultra thin, decoration={markings,mark=at position 1 with {\arrow[black,scale=2]{>}}},
    postaction={decorate},
    ]
    (r2) -- (3.5,10.4);

\node (r1) at (5, 10.9) {\begin{footnotesize}$\ref{3/4ub4}$\end{footnotesize}};
\draw[
    gray, ultra thin, decoration={markings,mark=at position 1 with {\arrow[black,scale=2]{>}}},
    postaction={decorate},
    ]
(r1) to (3.8,10.2);

\node (r1) at (6.5, 9.7) {\begin{footnotesize}$\ref{3/4ub3}$\end{footnotesize}};
%\draw[ gray, ultra thin, decoration={markings,mark=at position 1 with {\arrow[black,scale=2]{>}}}, postaction={decorate}, ] (r1) to (5.1,9.75);

\draw[
    gray, ultra thin, decoration={markings,mark=at position 1 with {\arrow[black,scale=2]{>}}},
    postaction={decorate},
    ]
(r1) to (5.15,9.3);

\draw[
    gray, ultra thin, decoration={markings,mark=at position 1 with {\arrow[black,scale=2]{>}}},
    postaction={decorate},
    ]
(r1) to (6.2,8.0);
\draw[
    gray, ultra thin, decoration={markings,mark=at position 1 with {\arrow[black,scale=2]{>}}},
    postaction={decorate},
    ]
(r1) to (5.8,8.4);

\node (r1) at (7.4, 9.0) {\begin{footnotesize}$\ref{edgecounting}$\end{footnotesize}};
\draw[
    gray, ultra thin, decoration={markings,mark=at position 1 with {\arrow[black,scale=2]{>}}},
    postaction={decorate},
    ]
(r1) to (5.5,9);

\node (r1) at (10.3, 8.3) {\begin{footnotesize}$\ref{3/4from5.7}$\end{footnotesize}};
%\draw[ gray, ultra thin, decoration={markings,mark=at position 1 with {\arrow[black,scale=2]{>}}}, postaction={decorate}, ] (r1) to (7,7.35);
%\draw[ gray, ultra thin, decoration={markings,mark=at position 1 with {\arrow[black,scale=2]{>}}}, postaction={decorate}, ] (r1) to (7.35,7);

\draw[
    gray, ultra thin, decoration={markings,mark=at position 1 with {\arrow[black,scale=2]{>}}},
    postaction={decorate},
    ]
(r1) to [bend right=15] (6.65,7.65);

\draw[
    gray, ultra thin, decoration={markings,mark=at position 1 with {\arrow[black,scale=2]{>}}},
    postaction={decorate},
    ]
(r1) [bend left = 20] to (7.7,6.6);

\node (r1) at (9.2, 6.3) {\begin{footnotesize}$\ref{3/4ub1}$\end{footnotesize}};
\draw[
    gray, ultra thin, decoration={markings,mark=at position 1 with {\arrow[black,scale=2]{>}}},
    postaction={decorate},
    ]
(r1) to (7.9,5.8);
%\draw[ gray, ultra thin, decoration={markings,mark=at position 1 with {\arrow[black,scale=2]{>}}}, postaction={decorate}, ] (r1) to (9,5.3);

\node (r1) at (12, 3.6) {\begin{footnotesize}$\ref{3/4ub2}$\end{footnotesize}};
%\draw[ gray, ultra thin, decoration={markings,mark=at position 1 with {\arrow[black,scale=2]{>}}}, postaction={decorate}, ] (r1) to (10.15,3.7);
\draw[
    gray, ultra thin, decoration={markings,mark=at position 1 with {\arrow[black,scale=2]{>}}},
    postaction={decorate},
    ]
(r1) to (10.2,3.05);
\draw[
    gray, ultra thin, decoration={markings,mark=at position 1 with {\arrow[black,scale=2]{>}}},
    postaction={decorate},
    ]
(r1) to (10.5,2.6);
\draw[ 
    gray, ultra thin, decoration={markings,mark=at position 1 with {\arrow[black,scale=2]{>}}}, 
    postaction={decorate}, 
    ] 
    (r1) to (11.1,2);

\node (r1) at (12.2, 1) {\begin{footnotesize}$\ref{maxbound}$\end{footnotesize}};
\draw[
    gray, ultra thin, decoration={markings,mark=at position 1 with {\arrow[black,scale=2]{>}}},
    postaction={decorate},
    ]
(r1) to (15*3/4,1);

\node (r1) at (1.3, 10) {\begin{footnotesize}$\ref{3/4psilb}$\end{footnotesize}};
\draw[
    gray, ultra thin, decoration={markings,mark=at position 1 with {\arrow[black,scale=2]{>}}},
    postaction={decorate},
    ]
(r1) to (1.3,10.85);
\draw[
    gray, ultra thin, decoration={markings,mark=at position 1 with {\arrow[black,scale=2]{>}}},
    postaction={decorate},
    ]
(r1) to (2.3,10.2);

\node (r1) at (9.8, .8) {\begin{footnotesize}$\ref{3/4psilb}$\end{footnotesize}};
\draw[
    gray, ultra thin, decoration={markings,mark=at position 1 with {\arrow[black,scale=2]{>}}},
    postaction={decorate},
    ]
(r1) to (10.8,1.3);
\draw[
    gray, ultra thin, decoration={markings,mark=at position 1 with {\arrow[black,scale=2]{>}}},
    postaction={decorate},
    ]
(r1) to (10.2,2.4);

%%%%%%%%%%%%%%%%%%%%%%%%%%%%%%%%%%%%%%%%%%%%%%%%%%%%%%%%%%%%%%%%%%%%%%%%%%%%%%%%%%%%%%%%%%%%%%%%%%%%%%%%%%%
\end{scope}

\draw[->] (0,0) -- (12.5,0) node[anchor=north]{$x$};
\draw[->] (0,0) -- (0,12.5) node[anchor=east] {$y$};

\node at (9.5,9.5) { \large{$\phi \ge 3/4$}};
\node at (4.3,5.1) { \large{$\phi < 3/4$}};

\node at (8.5,7.7) {\begin{scriptsize}$(1/2,1/2)$\end{scriptsize}};
\node at (11.1, 5.2)
{\begin{scriptsize}$(2/3,1/3)$\end{scriptsize}};
\node at (6, 10.3)
{\begin{scriptsize}$(1/3,2/3)$\end{scriptsize}};
\node at (2.9, 9.7)
{\begin{scriptsize}$(1/4,2/3)$\end{scriptsize}};
\node at (8.85,3.7)
{\begin{scriptsize}$(2/3,1/4)$\end{scriptsize}};

\node at (-1, 15*3/4) {\begin{scriptsize}$y=3/4$\end{scriptsize}};

\node at (15*3/4, -.4) {\begin{scriptsize}$x=3/4$\end{scriptsize}};

\draw (15*1/4,15*2/3) -- (15*1/3,15*2/3)--(15*1/3,15*56/97) -- (15*8/23, 15*56/97) -- (15*8/23, 15*8/15) -- (15*35/88, 15*8/15) -- (15*35/88, 15*10/19) -- (15*39/98, 15*10/19) -- (15*39/98, 15*1/2) -- (15*1/2,15*1/2) -- (15*1/2, 15*39/98) -- (15*10/19, 15*39/98) --(15*10/19, 15*35/88) --(15*8/15, 15*35/88) -- (15*8/15, 15*8/23) -- (15*56/97, 15*8/23) -- (15*56/97, 15*1/3) -- (15*2/3,15*1/3) -- (15*2/3,15*1/4);

\draw[domain=0:1/4,smooth,variable=\y] plot ({15*((3-10*\y)/(4-13*\y))},{15*\y});

\draw[domain=0:1/4,smooth,variable=\x] plot ({15*\x},{15*((3-10*\x)/(4-13*\x))});

\draw[domain=.4268:.4393, smooth, variable=\x, thick, dotted] plot ({15*(\x)}, {15*((4*\x-3)*(4*\x-3)/(16*\x*\x))});
\draw[domain=1/3:.4268,smooth,variable=\x, thick, dotted] plot
({15*\x},{15*(1-\x)});

\draw[domain=.4268:.4393,smooth,variable=\y, thick, dotted] plot
({15*((4*\y-3)*(4*\y-3)/(16*\y*\y))},{15*\y});
\draw[domain=1/3:.4268,smooth,variable=\y, thick, dotted] plot
({15*(1-\y)},{15*(\y)});

\draw[domain=.1886:1/3,smooth,variable=\y,thick,dotted] plot({15*((1-\y)/(1+\y-3*\y*\y))},{15*(\y)});

\draw[domain=.1886:1/3,smooth,variable=\x,thick,dotted] plot({15*\x},{15*((1-\x)/(1+\x-3*\x*\x))});

\draw[thick, dotted] (15*3/4, 15*.1886) -- (15*3/4, 0);
\draw[thick, dotted] (0, 15*3/4) -- (15*.1886, 15*3/4);
\draw[very thick, dotted] (15*.4393, 15*.5) -- (15*.5, 15*.5) -- (15*.5, 15*.4393);

\node (r1) at (1.5, 12.3) {\begin{footnotesize}$\ref{maxbound}$\end{footnotesize}};
\draw[
    gray, ultra thin, decoration={markings,mark=at position 1 with {\arrow[black,scale=2]{>}}},
    postaction={decorate},
    ]
(r1) to (1.5,15*3/4+.1);

\node (r1) at (12.3, 1.5) {\begin{footnotesize}$\ref{maxbound}$\end{footnotesize}};
\draw[
    gray, ultra thin, decoration={markings,mark=at position 1 with {\arrow[black,scale=2]{>}}},
    postaction={decorate},
    ]
(r1) to (15*3/4+.1,1.5);

\node (r1) at (9.1, 7.15) {\begin{footnotesize}$\ref{3/4from5.7}$\end{footnotesize}};
\draw[
    gray, ultra thin, decoration={markings,mark=at position 1 with {\arrow[black,scale=2]{>}}},
    postaction={decorate},
    ]
(r1) to (8,6.5);
%\draw[ gray, ultra thin, decoration={markings,mark=at position 1 with {\arrow[black,scale=2]{>}}}, postaction={decorate}, ] (r1) to (7.65,7);

\node (r1) at (12, 3.75) {\begin{footnotesize}$\ref{3/4phiub}$\end{footnotesize}};
\draw[
    gray, ultra thin, decoration={markings,mark=at position 1 with {\arrow[black,scale=2]{>}}},
    postaction={decorate},
    ]
(r1) to (10.8,3.6);

\node (r1) at (9.55, 1.75) {\begin{footnotesize}$\ref{3/4philb}$\end{footnotesize}};
\draw[
    gray, ultra thin, decoration={markings,mark=at position 1 with {\arrow[black,scale=2]{>}}},
    postaction={decorate},
    ]
(r1) to (10.8,2.5);

%%%%%%%%%%%%%%%%%%%%%%%%%%%%%%%%%%%%%%%%%%%%%%%%%%%%%%%%%%%%%%%%%%%%%%%%%%%%%%%%%%%%%%%%%%%%%%%%%%%%%%%%%%%
% \end{scope}

\end{tikzpicture}

\caption{When $\phi(x,y)<3/4$ and when $\psi(x,y)<3/4$.} \label{fig:psiphi34}
\end{figure}

In this section we investigate when $\psi(x,y) \ge 3/4$ and $\phi(x,y) \ge 3/4$. The results are shown in figure \ref{fig:psiphi34}.

\begin{thm}\label{3/4ub1}
Let $x,y\in (0,1]$, such that $y > 1/3$, $x \ge 1/2$, and $2y-2y^2 > 1-x$. Then $\psi(x,y) \ge 3/4$.
\end{thm}
\Proof
Let $G$ be a graph that is $(x,y)$-biconstrained via $(A,B,C)$. 
We can assume that $x,y$ are rational, and by multiplying vertices if necessary that $y|B| \in \mathbb{Z}$. Let $v_1 \in C$, and 
let $B_1 \subseteq N(v_1)$ be such that $|B_1| = y|B|$. Choose $v_2 \in C$ with at least $y|B \setminus B_1|=y(1-y)|B|$ neighbours 
in $B \setminus B_1$, and choose $B_2\subseteq N(v_2)$ with $|B_2|=y|B|$. Choose $v_3\in C$ with at least $y(1-2y)$ neighbours in
$B\setminus (B_1\cup B_2)$. Thus $|N(v_1)\cup N(v_2)|\ge y+y(1-y)>1-x$, and for $i=1,2$, $|N(v_i)\cup N(v_3)|\ge y+y(1-2y)>1-x$.

For $1\le i\le 3$, let $A_i=N_A^2(v_i)$. Since $y>1/3$, it follows that there exist $i,j$ with $1\le i<j\le 3$ 
such that $N(v_i)\cap N(v_j)\ne \emptyset$, and so $|A_i\cap A_j|\ge x|A|\ge |A|/2$. 
But $A_i\cup A_j=A$ since $|N(v_i)\cup N(v_j)|>1-x$, and
so $|A_i|+|A_j|\ge 3|A|/2$, and therefore one of $|A_i|,|A_j|\ge 3|A|/4$. This proves \ref{3/4ub1}.~\bbox

\begin{thm}\label{3/4ub2}
Let $x,y\in (0,1]$, such that $x > 2/3$, $x+2y>1$, and either $4(1-x)(1-y)\le 1$ or $x>1-2y+2y^2$.
Then $\psi(x,y) \ge 3/4$.
\end{thm}
\Proof
If $x+y> 1$ the result follows from (theorem 4.1 of the paper); so we may assume that $y\le 1-x<1/3$.
Let $G$ be $(x,y)$-biconstrained via $(A,B,C)$, and suppose that $|N_A^2(v)|<3|A|/4$ for each $v\in C$.
Let $H$ be the graph with $V(H) = V(C)$, where distinct $u,v$ are adjacent in $H$ if and only if $u$ and $v$ have distance two in $G$
(that is, in $G$ they have a common neighbour in $B$).
\\
\\
(1) {\em For all $w_1,w_2,w_3\in C$, if $N(w_1) \cap N(w_2) \ne\emptyset$ and $N(w_2) \cap N(w_3) \ne \emptyset$, then
$N(w_1) \cap N(w_3) \ne \emptyset$. Consequently each component of $H$ is a complete graph.}
\\
\\
Suppose that $w_1,w_2,w_3\in C$, and $N(w_1)\cap N(w_3)=\emptyset$, and $v_i \in N(w_i)\cap N(w_2)$ for $i = 1,3$.
Let $A_i = N_A^2(w_i)$ for $i = 1,2,3$. Since $x+2y>1$ we have $A_1 \cup A_3 = A$. Consequently $|A_1 \cap A_3| < |A|/2$, and since
$N_A(v_1) \cap N_A(v_3) \subseteq A_1 \cap A_3$, it follows that $|N_A(v_1) \cap N_A(v_3)|<|A|/2$. Thus
$$|N_A^2(w_2)| \ge |N_A(v_1) \cup N_A(v_3)| > (2x-1/2)|A| \ge 3|A|/4,$$
a contradiction. This proves (1).

\bigskip
Let $\alpha$ be the size of the largest stable set in $H$, that is, the number of components of $H$. 
Let the vertex sets of the components of $H$ be $C_1\ll C_{\alpha}$, and for $1\le i\le \alpha$ let $B_i$ be the set of 
vertices in $B$ with a
neighbour in $C_i$. The sets $B_1\ll B_{\alpha}$ have union $B$, and from the definition of $H$, they are pairwise disjoint.
For $1\le i\le \alpha$, let $w_i\in C_i$ and let $A_i=N^2_A(w_i)$. 
\\
\\
(2) {\em For $1\le i<j\le \alpha$, $A_i\cup A_j=A$, and so $|A_i\cap A_j|<|A|/2$. Consequently $\alpha\le 3$.}
\\
\\
Since $w_i, w_j$ have no common neighbour in $B$, it follows that $|N(w_i)\cup N(w_j)|\ge 2y|B|>(1-x)|B|$, and so
$A_i\cup A_j=A$. Since $|A_i|,|A_j|<3|A|/4$, it follows that $|A_i\cap A_j|<|A|/2$. This proves the first assertion.
Suppose that $\alpha\ge 4$. By the first assertion, every vertex in $A$ belongs to at least three
of $A_1\ll A_4$. Consequently some $A_i$ has cardinality at least $3|A|/4$, a contradiction. This proves (2).
\\
\\
(3) {\em $\alpha\ne 1$.}
\\
\\
Suppose that $\alpha=1$.  Every vertex in $A\setminus A_1$ has at least $x|B|$ neighbours in $B\setminus N(w_1)$, so we may
choose $v \in B$ with at least 
$$x|A \setminus A_1|/(1-y) > x|A|/(4-4y)\ge |A|/12$$ 
neighbours in $A \setminus A_1$. Let $w \in C$ be a neighbour of $v$. Since
$w_1,w$ have a common neighbour, it follows that
$$|N_A^2(w)| > (x+1/12)|A|\ge 3|A|/4,$$
a contradiction. This proves (3).
\\
\\
(4) {\em For $1\le i\le \alpha$, if $|B_i|>(1-x)|B|$ then $|C_i|>\frac{y}{3-4x}|C|$.}
\\
\\
Suppose that  $|B_i|>(1-x)|B|$ say, and let $|C_i|=c|C|$. Since $x > 2/3$, every vertex $u\in A \setminus A_i$ has a neighbour in 
$B_i$, and so $|N^2_C(u)\cap C_i|\ge y|C|= (y/c)|C_i|$. Hence there exists $w \in C_i$ such that 
$$|N^2_A(w)\cap (A \setminus A_i)|\ge \frac{y}{c}|A\setminus A_i|>\frac{y}{4c}|A|.$$
But $w,w_i$ have a common neighbour, and so 
$|N^2_A(w)\cap A_i|\ge x|A|$, and therefore
$$\frac34 |A|> |N^2_A(w)|\ge (x+\frac{y}{4c})|A|.$$
Consequently $3/4> x+y/(4c)$, and so $c>y/(3-4x)$. This proves (4).
\\
\\
(5) {\em $x>1-2y+2y^2$ and $\alpha=2$.}
\\
\\
Since $\alpha\le 3$, we may assume without loss of generality that $|B_1| \ge |B|/3> (1-x)|B|$. Since each $|C_i|\ge y|C|$,
it follows that $|C_1|\le (1-(\alpha-1)y)|C|$. By (4), $1-(\alpha-1)y >y/(3-4x)$, and since $\alpha\ge 2$, it follows that 
$1-y >y/(3-4x)$, that is, $4(1-x)(1-y)>1$. From the hypothesis it follows that $x>1-2y+2y^2$. This proves the first claim.
Suppose that $\alpha>2$; then
$$1-2y>\frac{y}{3-4x}>\frac{y}{3-4(1-2y+2y^2)}$$
which simplifies to $(1-y)(1-4y)^2<0$, a contradiction. This proves (5).

\bigskip

We may assume without loss of generality that $|C_1|\le |C|/2$. By (4), 
$|B_1|\le (1-x)|B|<|B|/2$, and so $|B_2|\ge |B|/2$.

Every vertex in $B_2\setminus N(w_2)$ is adjacent to at least a fraction $y/(1-y)$ of the vertices
of $C_2$, and hence there exists $w\in C_2$ with 
$$|N(w)\setminus N(w_2)|\ge \frac{y}{1-y}|B_2\setminus N(w_2)|.$$
Thus
$$|N(w)\cup N(w_2)|\ge |N(w_2)| + \frac{y}{1-y}|B_2\setminus N(w_2)| = \frac{y}{1-y}|B_2| + \frac{1-2y}{1-y}|N(w_2)|.$$
Since $|B_2|\ge x|B|$ and $|N(w_2)|\ge y|B|$, it follows that 
$$|N(w)\cup N(w_2)|\ge \left(\frac{xy}{1-y}+\frac{y(1-2y)}{1-y}\right)|B|>(1-x)|B|,$$
(because $x> 1-2y+2y^2$ by (5)).
Thus, $N_A^2(w) \cup N_A^2(w_2) = A$, and since $w$ and $w_2$ have a common neighbour in $B_2$ it follows that
$$|N_A^2(w)| + |N_A^2(w_2)| \ge (x+1)|A| \ge 3|A|/2$$
and so one of $|N^2_A(w)|,|N^2_A(w_2)|\ge 3|A|/4$, a contradiction. 
This proves \ref{3/4ub2}.~\bbox

\begin{thm}\label{3/4ub3}
Let $x,y\in (0,1]$, with $x>1/3$, $y > 1/2$, $x+3y>2$, and either $x\ge (5-6y)/(11-12y)$ or $x \ge (3-4y)/(4-4y)$. 
Then $\psi(x,y) \ge 3/4$.
\end{thm}
\Proof Let $G$ be $(x,y)$-biconstrained via $(A,B,C)$, and suppose that $N^2_A(v)<3|A|/4$
for each $v\in C$. By (theorem 4.1 of the paper) it follows that $x+y\le 1$. 
Let $H$ be the graph with $V(H) = V(C)$, in which distinct $u,v$
are adjacent if and only if $|N(u) \cup N(v)| \le (1-x)|B|$. 
\\
\\
(1) {\em For all $u,v\in C$, if $u,v$ are nonadjacent in $H$ then $N^2_A(u)\cup N^2_A(v)=A$. If $u,v$
are adjacent in $H$ then $|N^2_A(u)\cap N^2_A(v)|>(x+1/4)|A|$.}
\\
\\
If $u,v$ are nonadjacent in $H$, then $|N(u) \cup N(v)| > (1-x)|B|$, and so every vertex in $A$ has a  neighbour 
in $N(u) \cup N(v)$, that is, $N^2_A(u)\cup N^2_A(v)=A$. Now we assume that $u,v$ are adjacent in $H$. Consequently
$$|N(v_1)\cap N(v_2)| > 2y-(1-x)>1-y$$ 
by hypothesis. Moreover, $x>1/3$ and $y>1/2$ imply that
$$|N(u)\cap N(v)| > 2y+x-1>1/3>1-x/(1-y).$$
Thus, \ref{hompe1-y} implies that $|N^2_A(u)\cap N^2_A(v)|>(x+1/4)|A|$. This proves (1).

\bigskip
If there exist $w_1\ll w_4\in C$, pairwise nonadjacent in $H$, then by (1) each pair of the sets 
of the sets $N^2_A(w_i)\;(1\le i\le 4)$ has union $A$, and so each vertex in $A$ belongs to at least three of the four sets; and 
so one of the four sets has cardinality at least $3|A|/4$, a contradiction. Thus we may choose $w_1,w_2,w_3\in C$
such that every other vertex in $C$ is adjacent in $H$ to at least one of $w_1,w_2,w_3$. Choose a partition $C=C_1\cup C_2\cup C_3$
such that for $1\le i\le 3$, every vertex in $C_i$ is equal or adjacent in $H$ to $w_i$.
\\
\\
(2) {\em $|C_i|\le (1-y)|C|$ for $1\le i\le 3$.}
\\
\\
Suppose that $|C_1|>(1-y)|C|$. Define $B_1=N(w)$ and $A_1=N^2_A(w)$. Choose $v \in B \setminus N(w)$ with at least 
$x|A\setminus A_1|/(1-y) > x|A|/(4-4y)$ neighbours in $A \setminus A_1$. 
Since $|C_1|>(1-y)|C|$, there exists 
$w \in C_1$ adjacent to $v$. Then
$$|N_A^2(w)| > x|A|/(4-4y) + (x+1/4)|A| \ge 3|A|/4$$
since $x > 1/3$ and $y > 1/2$, a contradiction. This proves (2).
\\
\\
(3) {\em Every vertex in $B$ has neighbours in exactly two of $C_1,C_2,C_3$.}
\\
\\
Since each $|C_i|\le (1-y)|C|<y|C|$ by (2), it follows that every vertex in $B$ has neighbours in at least two of $C_1,C_2,C_3$.
Suppose that $v\in B$ has a neighbour $w_i'\in C_i$ for $i = 1,2,3$.
Let $A_i = N_A^2(w_i')$ for $i = 1,2,3$. For $1\le i<j\le 3$, $A_i\cup A_j=A$ by (1), and so every vertex of $A$ belongs to
at least two of $A_1,A_2,A_3$, and $N_A(u)$ is a subset of all three of $A_1,A_2,A_3$. Consequently
$$|A_1|+|A_2|+|A_3|\ge 2|A|+|N_A(w)|\ge (2+x)|A|\ge 9|A|/4$$
and so some $|A_i|\ge 3|A|/4$, a contradiction. This proves (3).

\bigskip

From (3) we may partition $B$ into $B_1,B_2,B_3$ such that
every vertex
in $B_1$ has neighbours in $C_2$ and in $C_3$ but not in $C_1$, and similarly for $B_2,B_3$. 
Without loss of generality,
we may assume that
$|B_1| \le 1/3$. Let $A_1=N^2_A(w_1)$.
\\
\\
(4) {\em $|N^2_A(w)\setminus A_1|< (2-4x)|A\setminus A_1|$ for each $w\in C_1$.}
\\
\\
By (1) $|N^2_A(w_1)\cap N^2_A(w)|\ge (x+1/4)|A|$, and since $|N^2_A(w)|<3|A|/4$,
it follows that 
$$|N^2_A(w)\setminus A_1|< (1/2-x)|A|\le (2-4x)|A\setminus A_1|.$$ 
This proves (4).

\bigskip

This has two consequences. The first is that $x < (3-4y)/(4-4y)$. To see this, by (4) we may choose $u\in A\setminus A_1$ such that $|N^2_C(u)\cap C_1|<(2-4x)|C_1|$.
Since $|B_1|\le |B|/3$, $u$ has a neighbour $v\in B_2\cup B_3$, and we may assume that $v\in B_2$ from the symmetry.
So at least $y|C|$ neighbours of $v$ belong to $C_1\cup C_3$, and therefore at least
$(y-(1-y))|C|$ neighbours belong to $C_1$, since $|C_2|\le (1-y)|C|$.
So $(2y-1)|C|<(2-4x)|C_1|\le (2-4x)(1-y)|C|$, and hence $x < (3-4y)/(4-4y)$ as claimed. 

The second consequence is that $x< (5-6y)/(11-12y)$. To see this, let $S= B\setminus (B_1\cup N(w_1))$. By (4)
and since each vertex in $B_2\cup B_3$ has a neighbour in $C_1$, it follows that each vertex $v\in B_2\cup B_3$
has fewer than $(2-4x)|A\setminus A_1|$ neighbours in $A\setminus A_1$. Since 
$S\subseteq B_2\cup B_3$, it follows that some vertex $u\in A\setminus A_1$
has fewer than $(2-4x)|S|$ neighbours in $S$. But $u$ has no neighbours in
$N(w_1)$, and only at most $r|B|$ neighbours in $B_1$; and since it has at least $x|B|$ neighbours in total, we deduce that
$$x|B|<(2-4x)|S|+r|B|\le  (2-4x)(1-r-y)|B|+r|B|$$ 
(since $B_1\cap N(w)=\emptyset$ and $|N(w_1)|\ge y|B|$ and therefore $|S|\le (1-r-y)|B|$). Consequently
$$x<  (2-4x)(1-r-y)+r= (2-4x)(1-y)+r(4x-1)\le (2-4x)(1-y)+(4x-1)/3$$
and so $x< (5-6y)/(11-12y)$.

We have shown then that $x < (3-4y)/(4-4y)$ and $x<(5-6y)/(11-12y)$; but this contradicts the hypothesis.
This proves \ref{3/4ub3}.~\bbox

\begin{thm}\label{edgecounting}
Let $x,y\in (0,1]$ with $1/2 < y < 2/3$ and $x>6y^2-8y+3$. Then $\psi(x,y) \ge 3/4$.
\end{thm}
\Proof 
We may assume that $y \in \mathbb{Q}$, by decreasing $y$ if necessary. Let $G$ be $(x,y)$-biconstrained via $(A,B,C)$, and suppose that
$|N^2_A(w)|<3|A|/4$ for each $w\in C$. By multiplying vertices if necessary, we may assume that $y|B| \in \mathbb{Z}$.
Since $x>6y^2-8y+3=6(y-2/3)^2+1/3$, it follows that $x>1/3$.
Let $H$ be the graph with $V(H) = C$ in which distinct $u,v\in C$ are adjacent  if and only if $|N(u) \cup N(v)| \le (1-x)|B|$. 
It follows that if $u,v$ are nonadjacent in $H$, then $N^2_A(u)\cup N^2_A(v)=A$.
As in the proof of \ref{3/4ub3}, there do not $w_1\ll w_4\in C$, pairwise nonadjacent in $H$; and so we may choose $w_1,w_2,w_3\in C$
and a partition $C=C_1\cup C_2\cup C_3$ such that for $1\le i\le 3$, every vertex in $C_i$ is equal or adjacent in $H$ to $w_i$.
Let $|C_i|=c_i|C|$, and choose $B_i\subseteq N_(w_i)$ with $|B_i|=y|B|$ for $1\le i\le 3$. 
Let $F$ be the set of all edges 
$vw$ of $G$ with $v\in B$ and $w\in C$, such that for $1\le i\le 3$, not both $v\in B_i$ and $w\in C_i$.
\\
\\
(1) {\em $\frac{|F|}{|B||C|}\le (1-x-y)< (2-3y)(2y-1)$.}
\\
\\
Let $w\in C$, with $w\in C_i$ say; then since $w,w_i$ are adjacent in $H$, it follows that
$w$ has at most $(1-x-y)|B|$ neighbours in $B\setminus B_i$. Thus $|F|\le (1-x-y)|B|\cdot |C|$. But $1-x-y<(2-3y)(2y-1)$
since $x>6y^2-8y+3$. This proves (1).

\bigskip

Let $p_1 = |B_1 \setminus (B_2 \cup B_3)|/|B|$, and define $p_2,p_3$ similarly. Let  
$q_1 = |(B_2\cup B_3) \setminus B_1|/|B|$, and define $q_2,q_3$ similarly. Let $p_0=|B\setminus (B_1\cup B_2\cup B_3)|/|B|$, and
$q_0=|B_1\cup B_2\cup B_3|/|B|$. Let $q=q_1+q_2+q_3$. Thus
\begin{eqnarray*}
p_0+p_1+p_2+p_3+q_0+q_1+q_2+q_3 &=& 1\\
p_1+q_0+q_2+q_3&=&y\\
p_2+q_0+q_3+q_1&=&y\\
p_3+q_0+q_1+q_2&=&y.
\end{eqnarray*}
By subtracting the last three of these from the first, we obtain
$$p_0-2q_0-(q_1+q_2+q_3)=1-3y,$$
and so $p_0=2q_0+q-3y+1$.

Every vertex in $B\setminus (B_1\cup B_2\cup B_3)$ is incident with at least $y|C|$ edges in $F$,
every vertex in $B_1 \setminus (B_2 \cup B_3)$ is incident with at least $(y-c_1)|C|$ edges in $F$, and every vertex
in $(B_2\cup B_3)\setminus B_1$ is incident with at least $\max(y-c_2-c_3,0) = \max(y+c_1-1,0)$ edges in $F$ (and similar
statements hold for $c_2,c_3$). Summing, we deduce that
$$\frac{|F|}{|B||C|}\ge p_0y+\sum_{1\le i\le 3}\left(p_i(y-c_i)+q_i\max(y+c_i-1,0)\right).$$
Since $p_i=y-q_0-q+q_i$ for $i = 1,2,3$, and $c_1+c_2+c_3=1$, it follows that
$$\sum_{1\le i\le 3}p_i(y-c_i)=\sum_{1\le i\le 3}(y-q_0-q+q_i)(y-c_i)=(y-q_0-q)(3y-1)+qy-\sum_{1\le i\le 3}q_ic_i.$$
Also, $p_0=2q_0+q-3y+1$, and so 
$$\frac{|F|}{|B||C|}\ge (2q_0+q-3y+1)y+(y-q_0-q)(3y-1)+qy-\sum_{1\le i\le 3}q_ic_i+\sum_{1\le i\le 3}q_i\max(y+c_1-1,0).$$
This simplifies to
$$\frac{|F|}{|B||C|}\ge (1-y)q_0+\sum_{i\in I}q_i(1-y-c_i)$$
where $I$ is the set of $i\in \{1,2,3\}$ such that $c_i<1-y$. From (1) we deduce that
$$(1-y)q_0+\sum_{i\in I}q_i(1-y-c_i)<(2-3y)(2y-1).$$
In particular it follows that $(1-y)q_0<(2-3y)(2y-1)\le (1-y)(2y-1)$, and so $q_0<2y-1$. Moreover, since $|B_2\cup B_3|\le |B|$,
it follows that $|B_2\cap B_3|\ge (2y-1)|B|$, and so $q_1\ge 2y-1-q_0$, and the same holds for $q_2,q_3$.
Consequently 
$$(1-y)q_0+\sum_{i\in I}(2y-1-q_0)(1-y-c_i)<(2-3y)(2y-1),$$
and so 
$$(1-y)q_0+\sum_{1\le i\le 3}(2y-1-q_0)(1-y-c_i)<(2-3y)(2y-1),$$
since $2y-1-q_0>0$. This simplifies to 
$(2y-1)q_0<0$, a contradiction. This proves \ref{edgecounting}.~\bbox

\begin{thm}\label{3/4ub4}
Let $x,y\in (0,1]$ with $x \ge 1/4$ and $y > 2/3$. Then $\psi(x,y) \ge 3/4$.
\end{thm}
\Proof
Suppose that $G$ is $(x,y)$-biconstrained via $(A,B,C)$, and $|N^2_A(w)|<3|A|/4$ for each $w\in C$.
Since $x \ge 1/4$ and $y > 2/3$, it follows that $x>2(1-y)^2$, that is, $y > (1-y)+1-x/(1-y)$ ; and also $y > 2-2y$.
Let $w\in C$. By \ref{hompe1-y}  with $k=2$ and $B'=N(w)$, it follows that 
$$|N_A^2(v)|>(x+1/2)|A| \ge 3|A|/4,$$
which is a contradiction. This proves \ref{3/4ub4}.~\bbox{}

\begin{thm}\label{3/4ub5}
Let $x,y\in (0,1]$ with $y> 2/3$, $x+4y>3$ and $x>3(1-y)^2/(2-y)$.
%$$2y+x-1>1-y+\max(1-y,1-x/(1-y)).$$
Then $\psi(x,y) \ge 3/4$.
\end{thm}
\Proof 
Suppose that $G$ is $(x,y)$-biconstrained via $(A,B,C)$, and $|N^2_A(w)|<3|A|/4$ for each $w\in C$. Consequently $y\le 3/4$,
and so $x\ge 3/20$ since $x>3(1-y)^2/(2-y)$.
From the hypotheses it follows
that 
$$2y+x-1>1-y+\max(1-y,1-x/(1-y)).$$
If $w_1,w_2 \in C$ with $|N(w_1) \cup N(w_2)| \le (1-x)|B|$ then $|N(w_1) \cap N(w_2)| \ge (2y+x-1)|B|$. Thus, \ref{hompe1-y} 
applied with $k=2$ tells us that, for all such $w_1,w_2 \in C$, more than $(x + 1/2)|A|$ vertices in $A$ have a neighbour in 
$N(v_1) \cap N(v_2)$. Let $H$ be the graph with vertex set $C$, in which $w_1,w_2$ are adjacent if $|N(w_1) \cup N(w_2)| \le (1-x)|B|$.
As in the proof of \ref{3/4ub3}, there is no stable set of size at least four in $H$.
It follows that there exist $w_1,w_2,w_3 \in C$ and a partition $C = C_1 \cup C_2 \cup C_3$ such that for $1\le i\le 3$, every
vertex in $C_i$ is equal to or adjacent in $H$ to $w_i$. We assume without loss of generality that $|C_1| \ge 1/3$. Since $y > 2/3$, 
every vertex in $B$ has a neighbour in $C_1$. Let $B_1= N(w_1)$ and $A_1 = N_A^2(w_1)$, and choose $v \in B \setminus B_1$ with 
more than $x|A|/(4-4y)$ neighbours in $A \setminus A_1$. 
Since $y > 2/3$, 
there exists $w\in C_1$ adjacent to $v$. Then
$$|N_A^2(w)| > (x+1/2+x/(4-4y))|A|\ge 3|A|/4$$
since $x\ge 3/20\ge 1/7$ and $y\ge 2/3$, a contradiction.
This proves \ref{3/4ub5}.~\bbox

\begin{thm}\label{3/4from5.7}
Let $x,y\in (0,1]$ with $y > 1/2$ and $x^2y \ge (3/4-x)^2$. Then $\phi(x,y) \ge 3/4$.
\end{thm}
\Proof
Apply \ref{hompe4} with $z = 3/4$.~\bbox{}

\begin{thm}\label{3/4phiub}
Let $x,y\in (0,1]$ with $y < 1/3$ and $x > \frac{1-y}{1+y-3y^2}$. Then $\phi(x,y) \ge 3/4$.
\end{thm}
\Proof
Suppose that $G$ is $(x,y)$-constrained via $(A,B,C)$, and $|N^2_A(w)|<3|A|/4$ for each $w\in C$. Consequently $x+y\le 1$.
By reducing $x$ or $y$
if necessary, we may assume that every vertex in $A$ has strictly more than $x|B|$ neighbours in $B$, and that 
$x,y$ are rational. Let 
$$p = \frac{1-x-y}{1-3y}.$$
By multiplying vertices, we may also assume that $y|B|$ and $p|B|$ are integers. Note that the hypotheses imply
that $p<xy$.
\\
\\
(1) {\em There exists $s\in [0,1]$ such that for all $b,c$, if $0\le a\le y$ and $sa+b\ge y(sy+1-y)$,
then $a+b\ge p$ and $b\ge 1-x-y$.}
\\
\\
We claim first that
$$\max\left(0,\frac{p-y(1-y)}{y^2}\right)\le \min\left(\frac{2y-y^2+x-1}{y-y^2}, 1\right).$$
To see this, we need to check that $0\le \frac{2y-y^2+x-1}{y-y^2}$, and $\frac{p-y(1-y)}{y^2}\le 1$, and 
$\frac{p-y(1-y)}{y^2}\le \frac{2y-y^2+x-1}{y-y^2}$. The first is true since
$$\frac{x}{1-y} > \frac{1}{1+y-3y^2}= 1-y+ \frac{4y^2-3y^3}{1+y-3y^2}\ge 1-y.$$
The second is true since $p\le xy\le y$. The third simplifies to $p/y\le x/(1-y)$, and this is true since $p\le xy$.
This proves the claim, and so there exists $s$ such that 
$$\max\left(0,\frac{p-y(1-y)}{y^2}\right)\le s\le \min\left(\frac{2y-y^2+x-1}{y-y^2}, 1\right).$$
We will show that $s$ satisfies (1). Suppose that $0\le a\le y$ and $sa+b\ge y(sy+1-y)$. Then
$$a+b\ge sa+b \ge  y(sy+1-y)\ge p$$
and 
$$sy+b\ge sa+b\ge y(sy+1-y)\ge sy+1-x-y$$
(and therefore $b\ge 1-x-y$). This proves (1).
\\
\\
(2) {\em There exists $t\in [0,1]$ such that for all $a,b$, if $0\le a\le 2p$ and $ta+b \ge y(1-2p(1-t))$ then $a+b\ge p$ and $p+b\ge 1-x$.}
\\
\\
We claim first that 
$$\max\left(0,\frac{2py+p-y}{2py}\right)\le \min\left(\frac{x+y+p-2py-1}{2p(1-y)},1\right).$$
To see this we must check that 
$0\le \frac{x+y+p-2py-1}{2p(1-y)}$, and $\frac{2py+p-y}{2py}\le 1$, and 
$\frac{2py+p-y}{2py}\le \frac{x+y+p-2py-1}{2p(1-y)}.$
The first is true since 
$$p-2py=\frac{(1-x-y)(1-2y)}{1-3y}\ge 1-x-y.$$ 
The second is true since $p\le xy\le y$; and the third simplifies to
$p\le xy$ and so is true. This proves the claim, and so there exists $t$ with 
$$\max\left(0,\frac{2py+p-y}{2py}\right)\le t\le \min\left(\frac{x+y+p-2py-1}{2p(1-y)},0\right).$$
We will show that $t$ satisfies (2).  Let $a,b$ satisfy $0\le a\le 2y$ and $ta+b \ge y(1-2p(1-t))$. Then 
$$a+b\ge ta+b \ge y(1-2p(1-t)) \ge p$$
and 
$$ 2tp + b \ge ta + b \ge y(1-2p(1-t)) \ge 2tp + 1 - x - p$$
(and so $b\ge 1-x-p$). This proves (2).

\bigskip

Choose $w_1 \in C$ with at least $y|B|$ neighbours in $B$.
\\
\\
(3) {\em There exists $w_2 \in C$ such that $|N(w_2)|\ge p|B|$ and $|N(w_1)\cup N(w_2)|\ge (1-x)|B|$.}
\\
\\
Choose $B_1 \subseteq N(w_1)$ with
$|B_1| = y|B|$. Choose $s$ as in (1).
Then
$$\sum_{w\in C}\left(s|N(w) \cap B_1|+ |N(w) \setminus B_1|\right)=\sum_{v\in B_1}s|N(v)\cap C|+\sum_{v\in B\setminus B_1}|N(v)\cap C|
\ge (sy^2+y(1-y))|B|\cdot|C|.$$
Consequently we may choose $w_2 \in C$ such that 
$$s\frac{|N(w_2) \cap B_1|}{|B|}+ \frac{|N(w_2) \setminus B_1|}{|B|} \ge y(sy+(1-y)).$$
Since 
$$0\le \frac{|N(w_2) \cap B_1|}{|B|}\le \frac{|B_1|}{|B|}=y,$$
the choice of $s$ implies that $|N(w_2)|\ge p|B|$ and $|N(w_2)\setminus B_1|\ge (1-x-y)|B|$, and so $|N(w_1)\cup N(w_2)|\ge (1-x)|B|$.
This proves (3).
\\
\\
(4) {\em There exists $w_3\in C$ such that $|N(w_3)|\ge p|B|$, and $|N(w_i)\cup N(w_3)|\ge (1-x)|B|$ for $i = 1,2$.}
\\
\\
Since $|N(w_1)|,|N(w_2)|\ge p|B|$ and $p|B|$ is an integer, we may choose $R\subseteq B$ with $|R|= 2p|B|$
such that $|N(w_1)\cap R|, |N(w_2)\cap R|\ge p|B|$. Choose $t$ as in (2). As in the proof of (3), there exists $w_3\in C$ with 
$$t\frac{|N(w_3) \cap R|}{|B|}+ \frac{|N(w_3) \setminus R|}{|B|} \ge y(2pt+(1-2p))=y(1-2y(1-t)).$$
From the choice of $t$, it follows that $|N(w_3)|/|B|\ge p$, and $|N(w_3)\setminus R|\ge 1-x-p$, and consequently
$|N(w_1)\cup N(w_3)|, |N(w_2)\cup N(w_3)|\ge (1-x)|B|$. This proves (4).

\bigskip

For $1\le i\le 3$, choose $B_i\subseteq N(w_i)$ with $|B_i|=p|B|$.
Since $|B_1\cup B_2\cup B_3|\le 3p$, we may choose $w_4 \in C$ with at least $y(1-3p)|B|$ neighbours in 
$B \setminus (B_1 \cup B_2 \cup B_3)$. Then for all $1 \le i \le 3$ we have:
$$|X_i \cup N(w_4)| \ge (p+y(1-3p))|B| \ge (1-x)|B|$$
by the definition of $p$. It follows that for $1\le i<j\le 4$, $|N(w_i)\cup N(w_j)|\ge (1-x)|B|$, and so (since every vertex in $A$
has strictly more than $x|B|$ neighbours in $B$) it follows that $N^2_A(w_i)\cup N^2_A(w_j)=A$. Thus every vertex in $A$ belongs to 
at least three of the four sets $N^2_A(w_i)\;(1\le i\le 4)$, and so one of them has cardinality at least $3|A|/4$, a contradiction.
This proves \ref{3/4phiub}.~\bbox

\begin{thm}\label{3/4psilb} Let $x,y\in (0,1]$. Then $\psi(x,y)<3/4$ if either:
\begin{itemize}
\item $x \le 1/6$ and $y\le 5/7$ and $2x+y \le 1$; or
\item $x \le 3/20$ and $x+4y\le 3$, and $x+4y<3$ if $x$ is irrational; or
\item $x \le 17/23$ and $y\le 1/8$ and $x+3y \le 1$; or
\item $x\le 5/7$ and $y\le 1/6$ and $x+2y\le 1$.
%x\le 13/19
\end{itemize}
\end{thm}
\Proof
If $x', y'$ with $2x'+y'\le 1$ and
$y'\le 3/5$, then $\psi(x',y')<2/3$ by the first bullet of \ref{psi23reversecurve}. Given $x,y$ as in the first bullet, 
the hypotheses imply that there is a choice of $x', y'$ with $2x'+y'\le 1$ and
$y'\le 3/5$, and which also satisfy the hypotheses of \ref{psilb2} 
with $z'=\psi(x',y')$ and $z$ slightly less than $3/4$ (checking this needs some lengthy calculation, which we omit); and 
so  the first statement follows from \ref{psilb2}.
The second statement follows similarly from \ref{psilb2} and the second bullet of \ref{psi23reversecurve}. The third statement 
follows 
from \ref{psilb1} and the first bullet of \ref{psi23extra}; and the fourth follows by applying \ref{psilb1} with 
$z=\max(2/7,x)$, taking $x'=3/5$, $y'=1/5$ and $z'=3/5$. This proves \ref{3/4psilb}.~\bbox

\begin{thm}\label{3/4philb}
If $x,y\in (0,1]$, with $y \le 1/3$ and 
$$\frac{x}{1-x}+\frac{y}{1-3y} \le 3,$$ 
then $\phi(x,y) < 3/4$.
\end{thm}
\Proof
Apply \ref{philb} with $z$ slightly less than $3/4$ to \ref{phi23curve}. This proves \ref{3/4philb}.~\bbox{}

\section{The 2/5 Level}
\begin{figure}[ht]
\centering

\begin{tikzpicture}[scale=.9,auto=left]

\begin{scope}[shift ={(-10,0)}]
\draw[->] (0,0) -- (7,0) node[anchor=north]{$x$};
\draw[->] (0,0) -- (0,7) node[anchor=east] {$y$};
\node at (5.6,5.15) {\begin{scriptsize}$(1/3,1/3)$\end{scriptsize}};
\node at (6.3,6.3) { \large{$\psi\ge 2/5$}};
\node at (3.2,3.2) { \large{$\psi< 2/5$}};
\node at (-.8,6) {\begin{scriptsize}$y=2/5$\end{scriptsize}};
\node at (6,-.33) {\begin{scriptsize}$x=2/5$\end{scriptsize}};

\draw
(15*2/5,15*0/1)--
(15*39/98, 15*0/1) -- (15*39/98, 15*1/98) -- 
(15*37/93, 15*1/98) -- (15*37/93, 15*4/93) -- 
(15*35/88, 15*4/93) -- (15*35/88, 15*1/22) -- 
(15*29/73, 15*1/22) -- (15*29/73, 15*4/73) -- 
(15*21/53, 15*4/73) -- (15*21/53, 15*3/53) -- 
(15*19/48, 15*3/53) -- (15*19/48, 15*1/16) -- 
(15*17/43, 15*1/16) -- (15*17/43, 15*3/43) -- 
(15*13/33, 15*3/43) -- (15*13/33, 15*1/11) -- 
(15*8/21, 15*1/11) -- (15*8/21, 15*2/21) -- 
(15*3/8, 15*2/21) -- (15*3/8, 15*1/8) -- 
(15*7/20, 15*1/8);
\draw (15*1/3,15*1/6) -- (15*1/3, 15*1/3) -- (15*1/6,15*1/3);
\draw
 (15*1/8, 15*7/20) 
-- (15*1/8, 15*3/8) -- (15*2/21, 15*3/8) 
-- (15*2/21, 15*8/21) -- (15*1/11, 15*8/21) 
-- (15*1/11, 15*13/33) -- (15*3/43, 15*13/33) 
-- (15*3/43, 15*17/43) -- (15*1/16, 15*17/43) 
-- (15*1/16, 15*19/48) -- (15*3/53, 15*19/48) 
-- (15*3/53, 15*21/53) -- (15*4/73, 15*21/53) 
-- (15*4/73, 15*29/73) -- (15*1/22, 15*29/73) 
-- (15*1/22, 15*35/88) -- (15*4/93, 15*35/88) 
-- (15*4/93, 15*37/93) -- (15*1/98, 15*37/93) 
-- (15*1/98, 15*39/98) -- (15*0/1, 15*39/98) 
--(15*0/1,15*2/5);

\draw[domain=1/8:1/6,smooth,variable=\y] plot({15*(2*(1-\y)/5)},{15*\y});
\draw[domain=1/8:1/6,smooth,variable=\x] plot({15*\x},{15*(2*(1-\x)/5)});

\draw[very thick, dotted] (15*1/3,15*1/3) -- (15*2/9,15*1/3);
\draw[domain=1/3:.3469,smooth,variable=\y,thick,dotted] plot({15*(2*\y*\y-3*\y+1)},{15*\y});
\draw[thick, dotted] (15*1/5,15*.3469) -- (15*1/5,15*2/5) -- (15*0,15*2/5);
\draw[very thick, dotted] (15*1/3,15*1/3) -- (15*1/3,15*1/4);
\draw[thick, dotted] (15*1/3,15*1/4) -- (15*.35,15*.25);
\draw[domain=.2123:.25,smooth,thick,dotted,variable=\y] plot({15*(.4 - \y/(10-20*\y))},{15*\y});
\draw[thick,dotted] (15*.3631,15*.2123) -- (15*.4,15*.2)--(15*.4,15*0);

\node (r1) at (1.7, 7.1) {\begin{footnotesize}$\ref{maxbound}$\end{footnotesize}};
\draw[
    gray, ultra thin, decoration={markings,mark=at position 1 with {\arrow[black,scale=2]{>}}},
    postaction={decorate},
    ]
(r1) to (1.7,6.05);

\node (r1) at (7.1, 1.5) {\begin{footnotesize}$\ref{maxbound}$\end{footnotesize}};
\draw[
    gray, ultra thin, decoration={markings,mark=at position 1 with {\arrow[black,scale=2]{>}}},
    postaction={decorate},
    ]
(r1) to (6.05,1.5);

\node (r1) at (4.3, 6) {\begin{footnotesize}$\ref{2/5psiub1}$\end{footnotesize}};
\draw[
    gray, ultra thin, decoration={markings,mark=at position 1 with {\arrow[black,scale=2]{>}}},
    postaction={decorate},
    ]
(r1) to (3.15,5.15);

\node (r1) at (5.8, 4.75) {\begin{footnotesize}$\ref{2/5psiub2}$\end{footnotesize}};
\draw[
    gray, ultra thin, decoration={markings,mark=at position 1 with {\arrow[black,scale=2]{>}}},
    postaction={decorate},
    ]
(r1) to (5.15,3.8);
\draw[
    gray, ultra thin, decoration={markings,mark=at position 1 with {\arrow[black,scale=2]{>}}},
    postaction={decorate},
    ]
(r1) to (5.4,3.5);
\draw[
    gray, ultra thin, decoration={markings,mark=at position 1 with {\arrow[black,scale=2]{>}}},
    postaction={decorate},
    ]
(r1) to (5.7,3.1);

\node (r1) at (1.7, 4.4) {\begin{footnotesize}$\ref{2/5psilb}$\end{footnotesize}};
\draw[
    gray, ultra thin, decoration={markings,mark=at position 1 with {\arrow[black,scale=2]{>}}},
    postaction={decorate},
    ]
(r1) to (2.1,5.1);

\node (r1) at (4.2, 1.7) {\begin{footnotesize}$\ref{2/5psilb}$\end{footnotesize}};
\draw[
    gray, ultra thin, decoration={markings,mark=at position 1 with {\arrow[black,scale=2]{>}}},
    postaction={decorate},
    ]
(r1) to (5.1,2.1);

%%%%%%%%%%%%%%%%%%%%%%%%%%%%%%%%%%%%%%%%%%%%%%%%%%%%%%%%%%%%%%%%%%%%%%%%%%%%%%%%%%%%%%%%%%%%%%%%%%%%%%%%%%%
\end{scope}
\draw[->] (0,0) -- (7,0) node[anchor=north]{$x$};
\draw[->] (0,0) -- (0,7) node[anchor=east] {$y$};
\node at (5.65,5.15) {\begin{scriptsize}$(1/3,1/3)$\end{scriptsize}};
\node at (6.3,6.3) { \large{$\phi\ge 2/5$}};
\node at (3.2,3.2) { \large{$\phi< 2/5$}};
\node at (-.8,6) {\begin{scriptsize}$y=2/5$\end{scriptsize}};
\node at (6,-.33) {\begin{scriptsize}$x=2/5$\end{scriptsize}};

\draw
(15*39/98, 15*.04348) -- (15*39/98, 15*1/11) -- 
(15*35/88, 15*1/11) -- (15*35/88, 15*1/9) -- 
(15*17/43, 15*1/9) -- (15*17/43, 15*1/8) -- 
(15*13/33, 15*1/8) -- (15*13/33, 15*2/15) -- 
(15*11/28, 15*2/15) -- (15*11/28, 15*1/7) -- 
(15*7/18, 15*1/7) -- (15*7/18, 15*2/13) -- 
(15*29/75,15*2/13);

\draw
(15*2/13, 15*29/75) 
-- (15*2/13, 15*7/18) -- (15*1/7, 15*7/18) 
-- (15*1/7, 15*11/28) -- (15*2/15, 15*11/28) 
-- (15*2/15, 15*13/33) -- (15*1/8, 15*13/33) 
-- (15*1/8, 15*17/43) -- (15*1/9, 15*17/43) 
-- (15*1/9, 15*35/88) -- (15*1/11, 15*35/88) 
-- (15*1/11, 15*39/98) -- (15*.04348, 15*39/98);

\draw
(15*9/25, 15*17/75) -- (15*9/25, 15*3/13) --
(15*5/14, 15*3/13) -- (15*5/14, 15*4/17) -- 
(15*11/31, 15*4/17) -- (15*11/31, 15*6/25) -- 
(15*6/17, 15*6/25) -- (15*6/17, 15*8/33) -- 
(15*7/20, 15*8/33) -- (15*7/20, 15*17/69) -- 
(15*9/26, 15*17/69) -- (15*9/26, 15*19/77) -- 
(15*.3387, 15*19/77);

\draw (15*1/4,15*1/3) -- (15*1/3,15*1/3) -- (15*1/3,15*1/4);

\draw (15*19/77,15*.3387) --
(15*19/77, 15*9/26) -- (15*17/69, 15*9/26) 
-- (15*17/69, 15*7/20) -- (15*8/33, 15*7/20) 
-- (15*8/33, 15*6/17) -- (15*6/25, 15*6/17) 
-- (15*6/25, 15*11/31) -- (15*4/17, 15*11/31) 
-- (15*4/17, 15*5/14) -- (15*3/13, 15*5/14) 
-- (15*3/13, 15*9/25) -- (15*17/75, 15*9/25);

\draw[domain=1/3:.3387,smooth,variable=\y] plot({15*(2-5*\y)/(7-17*\y)},{15*\y});
\draw[domain=9/25:29/75,smooth,variable=\y] plot({15*(2-5*\y)/(7-17*\y)},{15*\y});
\draw[domain=39/98:2/5,smooth,variable=\y] plot({15*(2-5*\y)/(7-17*\y)},{15*\y});
\draw[domain=2/13:17/75,smooth,variable=\x] plot({15*\x},{15*(2-7*\x)/(5-17*\x)});

\draw[domain=1/3:.3387,smooth,variable=\x] plot({15*\x},{15*(2-5*\x)/(7-17*\x)});
\draw[domain=9/25:29/75,smooth,variable=\x] plot({15*\x},{15*(2-5*\x)/(7-17*\x)});
\draw[domain=39/98:2/5,smooth,variable=\x] plot({15*\x},{15*(2-5*\x)/(7-17*\x)});
\draw[domain=2/13:17/75,smooth,variable=\y] plot({15*(2-7*\y)/(5-17*\y)},{15*\y});

\draw[very thick, dotted] (15*1/3,15*1/3)--(15*.297086,15*1/3);
\draw[thick, dotted] (15*.297086,15*1/3) -- (15*.297086,15*.3811);
\draw[very thick, dotted] (15*1/3,15*1/3)--(15*1/3,15*.297086);
\draw[thick, dotted] (15*1/3,15*.297086) -- (15*.3811,15*.297086);
\draw[domain=.3811:.4,smooth,thick,dotted,variable=\y] plot({15*.2*(6*\y*\y-2*sqrt(3*\y*\y*(3*\y*\y-5*\y+5))-5*\y+5)},{15*\y});
\draw[domain=.3811:.4,smooth,thick,dotted,variable=\x] plot({15*\x},{15*.2*(6*\x*\x-2*sqrt(3*\x*\x*(3*\x*\x-5*\x+5))-5*\x+5)});
\draw[thick, dotted] (15*.275,15*.4) -- (15*0,15*.4);
\draw[thick, dotted] (15*.4,15*.275)--(15*.4,15*0);

\node (r1) at (2.45, 7.1) {\begin{footnotesize}$\ref{maxbound}$\end{footnotesize}};
\draw[
    gray, ultra thin, decoration={markings,mark=at position 1 with {\arrow[black,scale=2]{>}}},
    postaction={decorate},
    ]
(r1) to (2.45,6.05);

\node (r1) at (7.1, 2.45) {\begin{footnotesize}$\ref{maxbound}$\end{footnotesize}};
\draw[
    gray, ultra thin, decoration={markings,mark=at position 1 with {\arrow[black,scale=2]{>}}},
    postaction={decorate},
    ]
(r1) to (6.05,2.45);

\node (r1) at (6.85, 4.8) {\begin{footnotesize}$\ref{2/5phiub1}$\end{footnotesize}};
\draw[
    gray, ultra thin, decoration={markings,mark=at position 1 with {\arrow[black,scale=2]{>}}},
    postaction={decorate},
    ]
(r1) to (5.9,4.3);

\node (r1) at (5.4, 5.75) {\begin{footnotesize}$\ref{2/5phiub2}$\end{footnotesize}};
\draw[
    gray, ultra thin, decoration={markings,mark=at position 1 with {\arrow[black,scale=2]{>}}},
    postaction={decorate},
    ]
(r1) to (4.5,5.5);

\node (r1) at (4.75, 2.5) {\begin{footnotesize}$\ref{2/5philb}$\end{footnotesize}};
\draw[
    gray, ultra thin, decoration={markings,mark=at position 1 with {\arrow[black,scale=2]{>}}},
    postaction={decorate},
    ]
(r1) to (5,3.7);
\draw[
    gray, ultra thin, decoration={markings,mark=at position 1 with {\arrow[black,scale=2]{>}}},
    postaction={decorate},
    ]
(r1) to (5.65,2.9);
\draw[
    gray, ultra thin, decoration={markings,mark=at position 1 with {\arrow[black,scale=2]{>}}},
    postaction={decorate},
    ]
(r1) to (5.9,.25);

\end{tikzpicture}

\caption{When $\psi(x,y)<2/5$ and when $\phi(x,y)<2/5$.} \label{fig:phi25}
\end{figure}
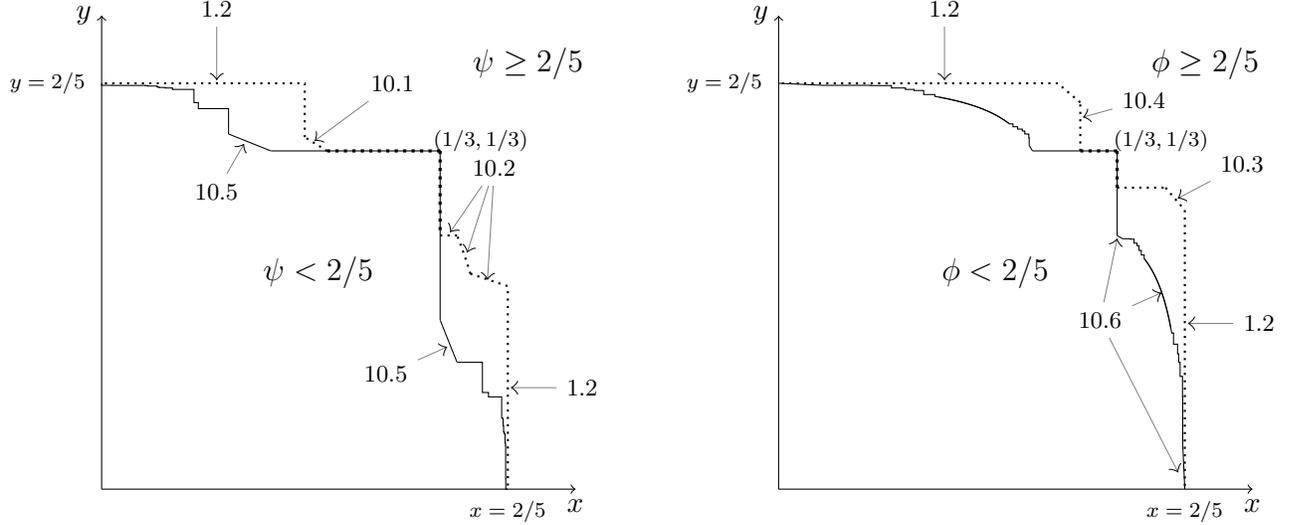

Next, we analyze when $\psi, \phi \ge 2/5$. The results are shown in figure \ref{fig:phi25}.

\begin{thm}\label{2/5psiub1}
Let $x,y\in (0,1]$ with $x \ge 1/5$, $y > 1/3$, and $3y-2y^2>1-x$; then $\psi(x,y) \ge 2/5$.
\end{thm}
\Proof
Suppose that $G$ is $(x,y)$-biconstrained via $(A,B,C)$, and $|N^2_A(w)|<(2/5)|A|$
for each $w\in C$. 
Let $w_1 \in C$, and let $A_1 = N_A^2(w_1)$. 
By averaging, there exists $w_2 \in C$ such that 
$$|A_2\setminus A_1|\ge y|A\setminus A_1|>(3y/5)|A|.$$ 
where $A_2=N^2_A(w_2)$.
Since
$|A_2|<2|A|/5$, it follows that 
$$|A_2\cap A_1|< (2/5-3y/5)|A|<x|A|,$$ 
and so $N(w_1)\cap N(w_2)=\emptyset$.
Let $B'=N(w_1)\cup N(w_2)$; thus $|B'|\ge 2y|B|$. By averaging, there exists $w_3\in C$ such that 
$|N(w_3)\setminus B'|\ge y|B\setminus B'|$,
and so 
$$|N(w_1)\cup N(w_2)\cup N(w_3)|\ge |B'|+y|B\setminus B'|=y|B|+(1-y)|B'|\ge (y+(1-y)(2y))|B|>(1-x)|B|.$$ 
Hence, setting $A_3 = N_A^2(w_3)$, it follows that $A_1\cup A_2\cup A_3=A$.

Since $y > 1/3$, some pair of $N(w_1), N(w_2), N(w_3)$ have nonempty intersection, and so some pair of $A_1,A_2,A_3$ have 
intersection of cardinality at least $x|A|\ge |A|/5$.
But then 
$$|A_1|+|A_2|+|A_3|\ge |A_1 \cup A_2 \cup A_3| +|A|/5=(6/5)|A|,$$ 
which is impossible since $|A_i|<(2/5)|A|$ for $1\le i \le 3$. 
This proves \ref{2/5psiub1}.~\bbox

\begin{thm}\label{2/5psiub2}
Let $x,y\in (0,1]$ with $x > 1/3$, $x+3y>1$, and either $y \ge 1/4$ or $x+y/(10(1-2y)) \ge 2/5$; then $\psi(x,y) \ge 2/5$.
\end{thm}
\Proof
Suppose that $G$ is $(x,y)$-biconstrained via $(A,B,C)$, and $|N^2_A(w)|<(2/5)|A|$ for each $w\in C$. Choose $w_1\ll w_{\alpha} \in C$
with $\alpha$
maximum such that $N(w_1)\ll N(w_{\alpha})$ are pairwise disjoint.
\\
\\
(1) {\em $\alpha=3$.}
\\
\\
Suppose that $\alpha\le 2$. Let $A'$ be the union of the sets $N^2_A(w_i)$ for $1\le i\le \alpha-1$, and let $B'$ be the union of the
sets $N(w_i)$  for $1\le i\le \alpha-1$. So 
$|A'|<(2\alpha/5)|A|\le (4/5)|A|$, and $|B'|\ge y|B|$. By averaging, there exists $v\in B\setminus B'$ such that 
$$|N(v)\cap (A\setminus A')|\ge (x/(1-y))|A\setminus A'|\ge (3x/(5(1-y)))|A|;$$
let $w\in C$ be adjacent to $v$. Since
$N(w)$ has nonempty intersection with $N(w_i)$ for some $i<\alpha$, it follows that $|N^2_A(w)\cap A'|\ge x|A|$.
Adding, we deduce that 
$$|N^2_A(w)|\ge  (3x/(5(1-y)))|A|+x|A|\ge (2/5)|A|,$$
a contradiction.
%3x/(5(1-y))+x\ge 2/5

Thus $\alpha\ge 3$; suppose that $\alpha\ge 4$.
Since $x + 3y > 1$ it follows that every vertex in $A$ belongs to at least two of the sets $N^2_A(w_i)\;(1\le i\le 4)$, and so one of these four sets has cardinality at least $|A|/2\ge (2/5)|A|$, a contradiction. This proves (1).
%x+3y>1
\\
\\
(2) {\em If $w\in C$ then $N(w)\cap N(w_i)$ is nonempty for exactly one value of $i\in \{1,2,3\}$.}
\\
\\
Since $\alpha=3$, it follows that $N(w)\cap N(w_i)$ is nonempty for at least one such $i$; suppose that
$N(w)$ has nonempty intersection with both $N(w_1), N(w_2)$ say. Let $A_i=N^2_A(w_i)$ for $i = 1,2,3$, and $A_0=N^2_A(w)$.
Since $A_1\cup A_2\cup A_3=A$, and each $|A_i|<(2/5)|A|$, there are fewer than $|A|/5$ vertices in $A$ that belong to more than 
one of $A_1,A_2,A_3$, and in particular $|A_1\cap A_2|<|A|/5$. But $|A_0\cap A_i|\ge x|A|$ for $i = 1,2$, and so 
$|A_0|\ge (2x-1/5)|A|\ge (2/5)|A|$, a contradiction. This proves (2).
%x\ge 3/10

\bigskip

From (2), we can partition $B=B_1\cup B_2\cup B_3$, and partition $C=C_1\cup C_2\cup C_3$, such that all six of these sets are 
nonempty,
and for all distinct $i,j\in \{1,2,3\}$ there is no edge between $B_i$ and $C_j$, and for all $i\in \{1,2,3\}$ and all $w,w'\in C_i$,
$N(w)\cap N(w')\ne \emptyset$. Let $|B_i|=b_i|B|$ and $C_i|=c_i|C|$
for $i = 1,2,3$. 
Without loss of generality we may assume that
$b_3\le 1/3<x$, and so
%x>1/3
every vertex in $A$ has a neighbour in $B_1\cup B_2$.
\\
\\
(3) {\em $x+y/(10(1-2y)) <2/5$ and so $y\ge 1/4$.}
\\
\\
Suppose that $x+y/(10(1-2y)) \ge 2/5$.
Without loss of generality, we may assume that at least $|A|/2$ vertices in $A$ have a neighbour in $B_1$. Choose $w\in C_1$.
Since 
$|N^2_A(w)|<(2/5)|A|$, there are at least $|A|/10$ vertices $u\in A \setminus N^2_A(w)$ that have a neighbour in $B_1$.
For each such $u$, $|N^2_C(u)\cap C_1|\ge y|C|$, and since $|C_1|\le (1-2y)|C|$ (because $|C_2|,|C_3|\ge y|C|$), it follows that
$|N^2_C(u)\cap C_1|\ge (y/(1-2y))|C_1|$. Consequently there exists $w'\in C_1$ such that $N^2_A(w')$ contains at least 
$(y/(10(1-2y)))|A|$ vertices in $A\setminus N^2_A(w)$. Since $w,w'$ have a common neighbour, it follows that $|N^2_A(w)\cap N^2_A(w')|\ge x|A|$, and so 
$$|N^2_A(w')|\ge (x+y/(10(1-2y)))|A|\ge (2/5)|A|,$$
a contradiction. Thus $x+y/(10(1-2y)) <2/5$, and so $y\ge 1/4$ from the hypothesis. This proves (3).

\bigskip

Since
$(b_1-y)+(b_2-y)+(b_3-y)=1-3y< x,$
it follows that for every vertex $u\in A$, there exists $i\in \{1,2,3\}$ such that $|N(u)\cap B_i|\ge (b_i-y)|B|$;
and consequently there is a partition $A=A_1\cup A_2\cup A_3$ such that 
for $i = 1,2,3$, every vertex in $A_i$ has more than $(b_i-y)|B|$ neighbours in $B_i$. It follows that
$A_i\subseteq N^2_A(w)$ for each $w\in C_i$. 
Let $|A_i|=a_i|A|$ for $i = 1,2,3$. 

For $i = 1,2$, let $D_i$ be the set of vertices in $A_3$ with a neighbour in $B_i$, and 
let $d_i = |D_i|/|A|$. For $i = 1,2$, if $u\in D_i$ then 
$$|N^2_C(u)\cap C_i|\ge y|C|\ge (y/(1-2y))|C_i|,$$
and so there exists $w\in C_i$ such that 
$$|N^2_A(w)\cap D_i|\ge  (y/(1-2y))|A_i|=(y/(1-2y))d_i|A|;$$ 
and since $A_i\subseteq N^2_A(w)$, it follows that 
$(y/(1-2y))d_i+a_i<2/5$. Since $d_1+d_2\ge a_3$, and $a_1+a_2=1-a_3$, summing for $i = 1,2$ yields
that
$$4/5>(y/(1-2y))(d_1+d_2)+(a_1+a_2)\ge(y/(1-2y))a_3+(1-a_3),$$
that is, $(1-3y)a_3/(1-2y)>1/5$; and since $a_3<2/5$, this implies that $y<1/4$, a contradiction.
This proves \ref{2/5psiub2}.~\bbox

\begin{thm}\label{2/5phiub1}
If $x,y\in (0,1]$ and $12x^2y \ge 5(1-x-y)^2$, then $\phi(x,y) \ge 2/5$.
\end{thm}
\Proof
Suppose not. Then $\phi(x,y) = 1-(3/5+\epsilon)$ for some $\epsilon > 0$, so by rotating we have $\phi(x,3/5+\epsilon) \le 1-y$. 
But \ref{hompe4} gives $\phi(x,3/5+\epsilon) > 1-y$, a contradiction. This proves \ref{2/5phiub1}.~\bbox

\begin{thm}\label{2/5phiub2}
If $x,y\in (0,1]$ and $x > 1/3$ and $y \ge (5-\sqrt{3})/11$, then $\phi(x,y) \ge 2/5$.
\end{thm}
\Proof Suppose not. Then $\phi(x,y) = 1 - (3/5+\epsilon)$ for some $\epsilon > 0$, so by rotating we have 
$\phi(3/5+\epsilon,y) \le 1-x < 2/3$. But $3/5 \ge (1-y)^2/(1-2y^2)$, since $y \ge (5-\sqrt{3})/11$; and so \ref{23weightedthm} 
gives that 
$\phi(3/5+\epsilon,y) \ge 2/3$, a contradiction. This proves \ref{2/5phiub2}.~\bbox

\begin{thm}\label{2/5psilb}
If $x,y\in (0,1]$, and either $5x/2 + y \le 1$ and $2y\le x$, or $x+5y/2 \le 1$ and $2x\le y$, then $\psi(x,y) < 2/5$.
\end{thm}
\Proof
Apply \ref{3.3gen} with $s/t = 2/5$. This proves \ref{2/5psilb}.~\bbox

\begin{thm}\label{2/5philb}
If $x,y\in (0,1]$ with $y \le 1/3$ and $x/(1-2x)+y/(1-3y) \le 2$, then $\phi(x,y) < 2/5$.
\end{thm}
\Proof
Apply \ref{phiaddpath} (with $z$ slightly less than $2/3$) to \ref{phi23curve}. This proves \ref{2/5philb}.~\bbox

\section{The 3/5 Level}
\begin{figure}[ht]
\centering

\begin{tikzpicture}[scale=.68,auto=left]

\begin{scope}[shift ={(-13,0)}]
\draw[->] (0,0) -- (10,0) node[anchor=north]{$x$};
\draw[->] (0,0) -- (0,10) node[anchor=east] {$y$};
\node at (8.4,7.7) {\begin{scriptsize}$(1/2,1/2)$\end{scriptsize}};
\node at (9,9) { \large{$\psi\ge 3/5$}};
\node at (5.3,5.3) { \large{$\psi< 3/5$}};
\node at (1.7,7.2) {\begin{scriptsize}$(1/6,1/2)$\end{scriptsize}};
\node at (6.45,2.5) {\begin{scriptsize}$(1/2,1/6)$\end{scriptsize}};
\node at (-0.8, 9) {\begin{scriptsize}$y=3/5$\end{scriptsize}};

%\node at (15*12/25+1.1,15/7) {\begin{scriptsize}$(12/25,1/7)$\end{scriptsize}};
%\node at (15*8/17+1, 15*1/6+.1) {\begin{scriptsize}$(8/17,1/6)$\end{scriptsize}};
%Start of insert

\draw (15*58/97, 15*0/1) -- (15*58/97, 15*1/97) -- 
(15*55/92, 15*1/97) -- (15*55/92, 15*3/92) -- 
(15*49/82, 15*3/92) -- (15*49/82, 15*2/41) -- 
(15*43/72, 15*2/41) -- (15*43/72, 15*1/18) -- 
(15*31/52, 15*1/18) -- (15*31/52, 15*3/52) -- 
(15*28/47, 15*3/52) -- (15*28/47, 15*3/47) -- 
(15*25/42, 15*3/47) -- (15*25/42, 15*1/14) -- 
(15*19/32, 15*1/14) -- (15*19/32, 15*3/32) -- 
(15*11/19, 15*3/32) -- (15*11/19, 15*2/19) -- 
(15*4/7, 15*2/19) -- (15*4/7, 15*1/7);

\draw (15*1/2, 15*1/6) -- (15*1/2, 15*1/2) --
(15*1/6, 15*1/2);

\draw (15*1/7, 15*4/7) -- (15*3/32, 15*93/160) ;
\draw (15*3/32, 15*93/160)
-- (15*3/32, 15*19/32) -- (15*1/14, 15*19/32) 
-- (15*1/14, 15*25/42) -- (15*3/47, 15*25/42) 
-- (15*3/47, 15*28/47) -- (15*3/52, 15*28/47) 
-- (15*3/52, 15*31/52) -- (15*1/18, 15*31/52) 
-- (15*1/18, 15*43/72) -- (15*2/41, 15*43/72) 
-- (15*2/41, 15*49/82) -- (15*3/92, 15*49/82) 
-- (15*3/92, 15*55/92) -- (15*1/97, 15*55/92) 
-- (15*1/97, 15*58/97) -- (15*0/1, 15*58/97);

\draw[domain=1/7:1/6,smooth,variable=\y] plot ({15*(1-3*\y)},{15*\y});
\draw[domain=1/7:1/6,smooth,variable=\x] plot ({15*(\x)},{15*(1-3*\x)});

% \draw[thick, dotted] (15*2/9, 10) -- (0,10);
% \draw[domain=.305:1/3, smooth, variable=\x, thick, dotted] plot ({15*\x}, {15*(1-\x)*(1-\x)/(1-2*\x*\x)});

\draw[very thick, dotted] (15*.2066, 15*1/2) -- (15*1/2, 15*1/2) -- (15*1/2, 15*1/4);
\draw[thick, dotted] (15*1/2, 15*1/4) -- (15*3/5, 15*1/5) -- (15*3/5, 0);

\draw[domain=.5:.5072,smooth,variable=\y,thick, dotted] plot({15*((5*\y*\y-6*\y+6)/6 - sqrt(25*\y*\y*\y*\y-60*\y*\y*\y+60*\y*\y)/6)},{15*\y});

% \draw[domain=.5:.5528, smooth, variable=\y, thick, dotted] plot ({15*(1-\y)*(1-\y)},{15*\y});
\draw[thick, dotted] (15*1/5, 15*.5072) -- (15*1/5, 15*3/5) -- (0, 15*3/5);

\node (r1) at (2.25, 10.1) {\begin{footnotesize}$\ref{maxbound}$\end{footnotesize}};
\draw[
    gray, ultra thin, decoration={markings,mark=at position 1 with {\arrow[black,scale=2]{>}}},
    postaction={decorate},
    ]
(r1) to (2.25,9.05);

\node (r1) at (10.1, 2.25) {\begin{footnotesize}$\ref{maxbound}$\end{footnotesize}};
\draw[
    gray, ultra thin, decoration={markings,mark=at position 1 with {\arrow[black,scale=2]{>}}},
    postaction={decorate},
    ]
(r1) to (9.05,2.25);

\node (r1) at (4.5, 8.5) {\begin{footnotesize}$\ref{3/5psiub1}$\end{footnotesize}};
\draw[
    gray, ultra thin, decoration={markings,mark=at position 1 with {\arrow[black,scale=2]{>}}},
    postaction={decorate},
    ]
(r1) to (3,8.5);
\draw[
    gray, ultra thin, decoration={markings,mark=at position 1 with {\arrow[black,scale=2]{>}}},
    postaction={decorate},
    ]
(r1) to (3.1,7.6);

\node (r1) at (8.5, 4.5) {\begin{footnotesize}$\ref{3/5psiub2}$\end{footnotesize}};
\draw[
    gray, ultra thin, decoration={markings,mark=at position 1 with {\arrow[black,scale=2]{>}}},
    postaction={decorate},
    ]
(r1) to (8.35,3.4);

\node (r1) at (1.15, 7.7) {\begin{footnotesize}$\ref{lostcurve35}$\end{footnotesize}};
\draw[
    gray, ultra thin, decoration={markings,mark=at position 1 with {\arrow[black,scale=2]{>}}},
    postaction={decorate},
    ]
(r1) to (2.25,8);
\draw[
    gray, ultra thin, decoration={markings,mark=at position 1 with {\arrow[black,scale=2]{>}}},
    postaction={decorate},
    ]
(r1) to (1.65,8.65);

\node (r1) at (7.8, 1.25) {\begin{footnotesize}$\ref{3/5psilb}$\end{footnotesize}};
\draw[
    gray, ultra thin, decoration={markings,mark=at position 1 with {\arrow[black,scale=2]{>}}},
    postaction={decorate},
    ]
(r1) to (8,2.25);

%%%%%%%%%%%%%%%%%%%%%%%%%%%%%%%%%%%%%%%%%%%%%%%%%%%%%%%%%%%%%%%%%%%%%%%%%%%%%%%%%%%%%%%%%%%%%%%%%%%%%%%%%%%
\end{scope}
\draw[->] (0,0) -- (11,0) node[anchor=north]{$x$};
\draw[->] (0,0) -- (0,11) node[anchor=east] {$y$};
\node at (8.3,7.7) {\begin{scriptsize}$(1/2,1/2)$\end{scriptsize}};
\node at (-1,9) {\begin{scriptsize}$y=3/5$\end{scriptsize}};
\node at (5,5) {\large{$\phi<3/5$}};
\node at (9,9) { \large{$\phi\ge 3/5$}};

%\draw[domain=1/2:.34, smooth, variable=\y, thick, dotted] plot ({15*(1-\y*(2-\y)/(1+\y))},{15*\y});
%\draw[domain=.305:(5-sqrt(13))/6, smooth, variable=\y, thick, dotted] plot ({15*(1-\y*(2-\y)/(1+\y))},{15*\y});
%\draw[domain=1/2:(5-sqrt(13))/6, smooth, variable=\y, thick, dotted] plot ({15*\y},{15*(1-\y*(2-\y)/(1+\y))});

%  \draw[domain=58/99:11/20, smooth, variable=\x] plot({15*\x},{15*(3-5*\x)/(11-18*\x)});

\draw[domain=19/77:1/4, smooth, variable=\y] plot({15*(3-11*\y)/(5-18*\y)}, {15*\y});
\draw[domain=2/13:22/97, smooth, variable=\y] plot({15*(3-11*\y)/(5-18*\y)}, {15*\y});
\draw[domain=19/77:1/4, smooth, variable=\x] plot({15*\x}, {15*(3-11*\x)/(5-18*\x)});
\draw[domain=2/13:22/97, smooth, variable=\x] plot({15*\x}, {15*(3-11*\x)/(5-18*\x)});

\draw
(15*58/97, 15*0/1) -- (15*58/97, 15*1/11) -- 
(15*55/92, 15*1/11) -- (15*55/92, 15*1/10) -- 
(15*52/87, 15*1/10) -- (15*52/87, 15*1/9) -- 
(15*25/42, 15*1/9) -- (15*25/42, 15*1/8) -- 
(15*19/32, 15*1/8) -- (15*19/32, 15*2/15)-- 
(15*16/27, 15*2/15) -- (15*16/27, 15*1/7) -- 
(15*10/17, 15*1/7) -- (15*10/17, 15*2/13) -- 
(15*58/99, 15*2/13);

\draw
(15*11/20, 15*22/97) -- (15*11/20, 15*3/13) -- 
(15*6/11, 15*3/13) -- (15*6/11, 15*4/17) --
(15*13/24, 15*4/17) -- (15*13/24, 15*6/25) -- 
(15*7/13, 15*6/25) -- (15*7/13, 15*8/33) -- 
(15*8/15, 15*8/33) -- (15*8/15, 15*17/69) -- 
(15*9/17, 15*17/69) -- (15*9/17, 15*18/73) -- 
(15*10/19, 15*18/73) -- (15*10/19, 15*19/77) --
(15*.5116, 15*19/77);

\draw
(15*1/2, 15*1/4) -- (15*1/2, 15*1/2) -- (15*1/4, 15*1/2);

\draw (15*19/77, 15*.5116) 
-- (15*19/77, 15*10/19) -- (15*18/73, 15*10/19) 
-- (15*18/73, 15*9/17) -- (15*17/69, 15*9/17) 
-- (15*17/69, 15*8/15) -- (15*8/33, 15*8/15) 
-- (15*8/33, 15*7/13) -- (15*6/25, 15*7/13) 
-- (15*6/25, 15*13/24) -- (15*4/17, 15*13/24) 
-- (15*4/17, 15*6/11) -- (15*3/13, 15*6/11) 
-- (15*3/13, 15*11/20) -- (15*22/97, 15*11/20);

 \draw (15*2/13, 15*58/99) 
-- (15*2/13, 15*10/17) -- (15*1/7, 15*10/17) 
-- (15*1/7, 15*16/27) -- (15*2/15, 15*16/27) 
-- (15*2/15, 15*19/32) -- (15*1/8, 15*19/32) 
-- (15*1/8, 15*25/42) -- (15*1/9, 15*25/42) 
-- (15*1/9, 15*52/87) -- (15*1/10, 15*52/87) 
-- (15*1/10, 15*55/92) -- (15*1/11, 15*55/92) 
-- (15*1/11, 15*58/97) -- (15*0/1, 15*58/97);
 
 \draw[very thick, dotted] (15*.3167, 15*1/2) -- (15*1/2, 15*1/2) -- (15*1/2, 15*.3167);

\draw[domain=.3031:.3167, smooth, thick, dotted, variable=\x] plot({15*\x}, {15*(5*\x-3)*(5*\x-3)/(40*\x*\x)});
\draw[thick, dotted] (0, 15*.6) -- (15*.3031, 15*.6);
\draw[domain=.3031:.3167, smooth, thick, dotted, variable=\y] plot({15*(5*\y-3)*(5*\y-3)/(40*\y*\y)}, {15*\y});
\draw[thick, dotted] (15*.6, 0) -- (15*.6, 15*.3031);

\node (r1) at (2.45, 10.1) {\begin{footnotesize}$\ref{maxbound}$\end{footnotesize}};
\draw[
    gray, ultra thin, decoration={markings,mark=at position 1 with {\arrow[black,scale=2]{>}}},
    postaction={decorate},
    ]
(r1) to (2.45,9.05);

\node (r1) at (6.8, 1.6) {\begin{footnotesize}$\ref{3/5philb}$\end{footnotesize}};
\draw[
    gray, ultra thin, decoration={markings,mark=at position 1 with {\arrow[black,scale=2]{>}}},
    postaction={decorate},
    ]
(r1) to (7.6,3.7);
\draw[
    gray, ultra thin, decoration={markings,mark=at position 1 with {\arrow[black,scale=2]{>}}},
    postaction={decorate},
    ]
(r1) to (8.7,2.6);
\draw[
    gray, ultra thin, decoration={markings,mark=at position 1 with {\arrow[black,scale=2]{>}}},
    postaction={decorate},
    ]
(r1) to (8.9,.5);

\node (r1) at (8.5, 6.25) {\begin{footnotesize}$\ref{3/5phiub1}$\end{footnotesize}};
\draw[
    gray, ultra thin, decoration={markings,mark=at position 1 with {\arrow[black,scale=2]{>}}},
    postaction={decorate},
    ]
(r1) to (8.25,4.7);

\end{tikzpicture}

\caption{When $\psi(x,y)<3/5$ and when $\phi(x,y)<3/5$.} \label{fig:phi35}
\end{figure}
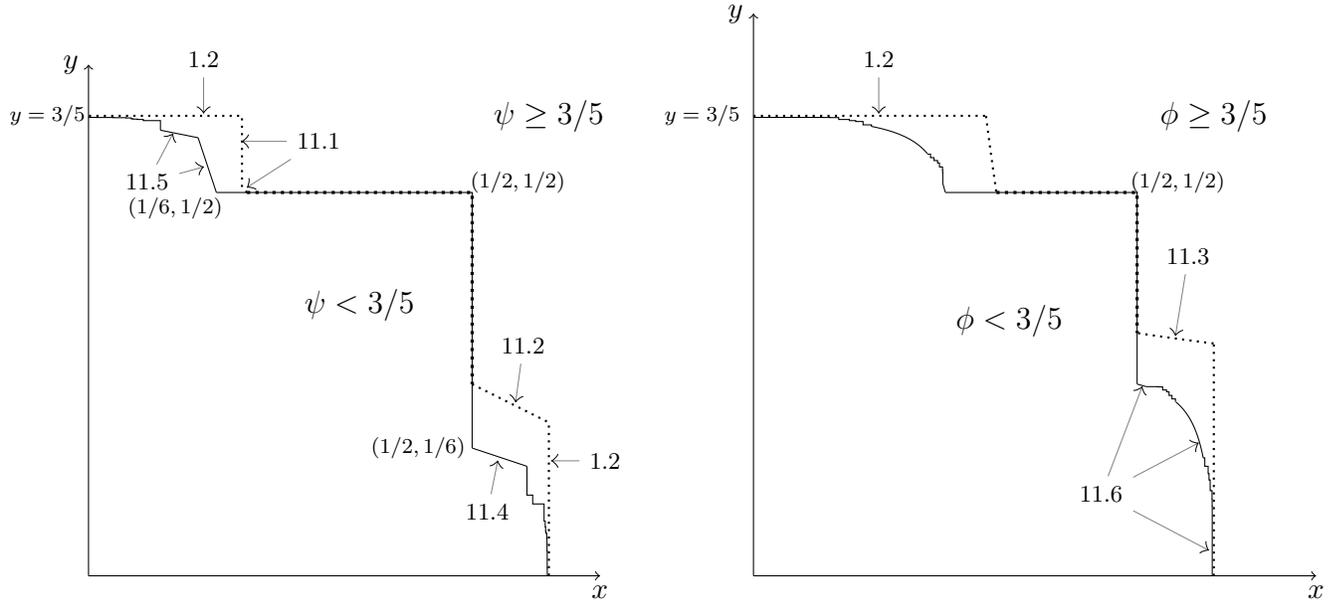

Next, we analyze when $\psi \ge 3/5$, and similarly for $\phi$. The results are shown in figure \ref{fig:phi35}.

\begin{thm}\label{3/5psiub1}
If $x,y\in (0,1]$ with $y > 1/2$, $x \ge 1/5$, and 
$$2y-\frac{3-5x}{3-3x}y^2>1-x,$$ 
then $\psi(x,y) \ge 3/5$.
\end{thm}
\Proof
Suppose not, and let $G$ be $(x,y)$-biconstrained via $(A,B,C)$, such that $|N_A^2(w)|<3|A|/5$ for all $w \in C$. We can 
assume that 
$x,y$ are rational, and by multiplying vertices if necessary, we can assume that both $x|B|$ and $|C|/5$ are integers. By averaging, 
there exists $u \in A$ such that $|N_C^2(u)|<3|C|/5$. Choose $B' \subseteq N(u)$ with $|B'| = x|B|$, and choose 
$C' \subseteq C$ with $N_C^2(u) \subseteq C'$ and $|C'|=3|C|/5$.
\\
\\
(1) {\em There exist $w_1\in C'$ and $w_2\in C\setminus C'$ such that $|N(w_1) \cup N(w_2)|>(1-x)|B|$.}
\\
\\
Choose $w_1 \in C'$ and $w_2 \in C \setminus C'$ uniformly and independently at random, and let 
$S = N(w_1) \cup N(w_2)$. We will compute the expectation of $|S|$.
Let $t=2/5$. For each $v \in B_1$, $v$ has at least $y|C|$ neighbours in $C'$, so the probability it is in $S$ is at least 
$y/(1-t)$. For each $v \in B \setminus B'$, define $d(v) = |N(v) \cap (C \setminus C_1)|/|C|$. Then the probability that $v$ 
is a neighbour of $w_2$ is $d(v)/t$, and the probability that $v$ is a neighbour of $w_1$ and not a neighbour of $w_2$
is at least 
$$\left(1-\frac{d(v)}{t}\right)\frac{y-d(v)}{1-t}.$$ Thus the 
expectation of $|S|$ is at least
$$\sum_{v \in B'} \frac{y}{1-t}+\sum_{v \in B \setminus B'} \left(\frac{d(v)}{t}+\left(1-\frac{d(v)}{t}\right)\frac{y-d(v)}{1-t}\right)
=\sum_{v \in B} \frac{y}{1-t}+ \sum_{v \in B \setminus B'} \frac{(1-y-2t)d(v)+d(v)^2}{t(1-t)}.$$
Choose $q$ with $$\sum_{v \in B \setminus B'} d(u) = qt|B|.$$ 
Each vertex $w \in C \setminus C'$ has at least $y|B|$ neighbours in $B \setminus B'$, so it follows that 
$q \ge y$. Thus, the expectation of $|S|$ is at least 
$$\frac{y}{1-t}|B|+\frac{1-y-2t}{1-t}q|B| +
\sum_{v \in B \setminus B'}\frac{d(v)^2}{t(1-t)}.$$
%\[qt|B|(1-2t-y)+\left(\sum_{u \in B \setminus B_1} d(u)^2 \right)+ yt|B| > t(1-t)(1-x)|B|\]
Since $\sum_{v \in B \setminus B'} d(v) = qt|B|$ and $|B \setminus B'| = (1-x)|B|$, it follows by Cauchy-Schwarz that 
$$\sum_{v \in B \setminus B'} d(v)^2 \ge \frac{q^2t^2}{1-x}|B|.$$ 
Thus the expectation of $|S|$ is at least
$$\frac{y}{1-t}|B|+\frac{1-y-2t}{1-t}q|B| +
 \frac{q^2t}{(1-x)(1-t)}|B|=\left(y+q(1-y-2t)+\frac{q^2t}{1-x}\right)\frac{|B|}{1-t}.$$

To prove (1), it suffices to show that the expectation of $|S|$ is more than $(1-x)|B|$, and so it suffices to show
that
$$y+q(1-y-2t)+\frac{q^2t}{1-x} > (1-t)(1-x).$$
Remembering that $t=2/5$, this is
$$2(1-5y)q+\frac{4q^2}{1-x}+10y > 6(1-x),$$
and the derivative of the left-hand side with respect to $q$ is
$$2(1-5y)+\frac{8q}{1-x} \ge 2(1-5y)+\frac{8y}{4/5} = 2 > 0$$
for $q \ge y$. It follows that the left-hand side is minimized when $q = y$, so it suffices to show
that
$$2(1-5y)y+\frac{4y^2}{1-x}+10y>6(1-x),$$
which is equivalent to the hypothesis. This proves (1).

\bigskip

Let $w_1,w_2$ be as in (1), and let $A_i = N_A^2(w_i)$ for $i = 1,2$. Then $A = A_1 \cup A_2$, and since $y > 1/2$ 
we have $N(v_1) \cap N(v_2) \ne \emptyset$, and consequently $|A_1 \cap A_2| \ge x|A| \ge |A|/5$. Then 
$$|A_1|+|A_2|=|A_1 \cup A_2| + |A_1 \cap A_2| \ge 6|A|/5$$ 
and so we have $|A_i| \ge 3|A|/5$ for some $i$, a contradiction. This proves \ref{3/5psiub1}.~\bbox

\begin{thm}\label{3/5psiub2}
If $x,y\in (0,1]$ with $x > 1/2$ and $x+2y>1$, then $\psi(x,y) \ge 3/5$.
\end{thm}
\Proof
Let $G$ be $(x,y)$-biconstrained via $(A,B,C)$, and suppose that $|N^2_A(w)|<(3/5)|A|$ for each $w\in C$.
Choose $\alpha$ maximum such that there exist $w_1\ll w_{\alpha}\in C$ where $N(w_i)\cap N(w_j)=\emptyset$ for $1\le i<j\le \alpha$.
\\
\\
(1) {\em $\alpha=2$.}
\\
\\
Since  $x+2y>1$ it follows that $N^2_A(w_i)\cup N^2_A(w_j)=A$ for all distinct $i,j\in \{1\ll \alpha$, and so if $\alpha\ge 3$
then every vertex in $A$ belongs to at least two of the sets $N^2_A(w_1),N^2_A(w_2),N^2_A(w_3)$, which is impossible since they each
have cardinality less than $(3/5)|A|$. So $\alpha\le 2$.

Suppose that $\alpha=1$; then 
every $w_2 \in C$ satisfies $N(w_1) \cap N(w_2) \ne \emptyset$. Let $B_1 = N(w_1)$ and $A_1 = N_A^2(w_1)$; thus $|A_1|<(3/5)|A|$ and $|B_1|\ge y|B|$. Choose $v \in B \setminus B_1$ with $v$ at least $x|A \setminus A_1|/(1-y) > 2x|A|/(5(1-y))$ neighbours in 
$A \setminus A_1$, and let $w_2 \in C$ be a neighbour of $v$. Then
$$(3/5)|A|>|N_A^2(w_2)| > x|A| + 2x|A|/(5(1-y)) > \left(x+\frac{4x}{5(x+1)}\right)|A| \ge 3|A|/5$$
since $x > 1/2$, a contradiction. This proves (1).

\bigskip

Since $\alpha=2$, every vertex $w \in C$ shares a neighbour with at least one of $w_1, w_2$. Let $A_i=N^2_A(w_i)$
for $i = 1,2$. Since $x+2y>1$, we have $A = A_1 \cup A_2$, and so $|A_1 \cap A_2| < |A|/5$ because $|A_1|, |A_2| < 3|A|/5$.
Then, if some $w \in C$ shares a neighbour with $w_1$ and shares a neighbour with $w_2$, it follows that 
$|N_A^2(w)| > 2x|A| - |A|/5 > 4|A|/5$, a contradiction.

Thus, every vertex in $C$ shares a neighbour with exactly one of $w_1$ and $w_2$. Let $H$ be the bipartite graph $G[B\cup C]$.
It follows that there are exactly two components of $H$, say $H_1, H_2$, where $w_i\in V(H_i)$ for $i = 1,2$.
Let $B_i = B \cap H_i$ and $C_i = C \cap H_i$ for $i = 1,2$. Without loss of generality we may assume that $|B_1| \ge |B|/2$. 
It follows that for each $u\in A$,  $u$ has a neighbour in $B_1$ 
and consequently 
$$|N^2_C(u)\cap C_1|\ge y|C|\le \frac{y}{1-y}|C_1|,$$
because $|C_2| \ge y|C|$, and thus $|C_1| \le (1-y)|C|$.
Since $|A\setminus A_1|>(2/5)|A|$, there exists $w \in C_1$ with more than $2|A|y/(5(1-y))$ neighbours in $A \setminus A_1$. 
But $w$ and $w_1$ share a common neighbour, so
$$|N_A^2(w)| > \frac{2y|A|}{5(1-y)}+x|A| > \left(\frac{2y}{5(1-y)}+(1-2y)\right)|A| \ge 3|A|/5$$
since the last inequality is equivalent to $5y^2-5y+1 \ge 0$, which is true because $y \le 1/4$ (since $x+2y>1$, and $x > 1/2$).
This proves \ref{3/5psiub2}.~\bbox

\begin{thm}\label{3/5phiub1}
If $x,y\in (0,1]$ with $y > 1/2$ and $40x^2y \ge (3-5x)^2$, then $\phi(x,y) \ge 3/5$. 
\end{thm}
\Proof
Apply \ref{hompe4} with $z = 3/5$. This proves \ref{3/5phiub1}.~\bbox

\begin{thm}\label{3/5psilb}
If $x,y\in (0,1]$ with $x\le 4/7$ and $y\le 1/2$ and $x+3y\le 1$, then $\psi(x,y)<3/5$.
\end{thm}
\Proof
We may assume that $x>1/2$ since $\psi(1/2,1/2)<3/5$; and so $y<1/6$ since $x+3y\le 1$.
The claim follows from applying \ref{psilb1} with $z$ slightly less than $3/5$ and $x'=y'=z'=1/4$.
This proves \ref{3/5psilb}.~\bbox

\begin{thm}\label{lostcurve35}
If $x,y\in (0,1]$, such that $3x+y\le 1$, and $x+5y\le 3$, with strict inequality in both if $x$ or $y$ is irrational, then
$\psi(x,y)<3/5$.
\end{thm}
\Proof 
By increasing $x, y$ if necessary, we may assume that $x,y$ are rational. Suppose that $\psi(x,y)\ge 3/5$. We claim first that:
\\
\\
(1) {\em $x<1/6$, and $y>1/2$, and $5xy+15y< 9$, and  $x< (1-y)/(5y)$, and $x\le 3(1-y)^2/(1+5y)$.}
\\
\\
Since $3x+y\le 1$, it follows that $x<1/3$, and $y>1/2$ since $\psi(1/2,1/2)=1/2<3/5$. Thus $x< 1/6$, since
$3x+y\le 1$. This proves the first two statements.
Since $x+5y\le 3$, it follows that $y<3/5$, and so $5xy+15y<3x+15y\le 9$. This proves the third statement.
For the fourth, $5x<3x/y$ (since $y<3/5$), and $3x\le 1-y$, and so $5x<(1-y)/y$. 
Finally, for the fifth statement, if $y\le 4/7$, then $1+5y\le 9-9y$, and so 
$$x\le 9x(1-y)/(1+5y)\le 3(1-y)^2/(1+5y);$$
and if $y\ge 4/7$, then 
$$(3-5y)(1+5y)=3+10y-25y^2\le 3-6y+3y^2= 3(1-y)^2$$ 
and so $$x\le 3-5y\le 3(1-y)^2/(1+5y).$$
This proves (1).

\bigskip
Since $x<1/6$, it follows that $x/(1-x)< (5x)/3$ and $(1-y)/(3y)<1/3$.
The hypotheses (via (1)) imply that
$$\frac{2x}{3-3y-x}\le \min\left(\frac3y-5, \frac{1-y}{3y}\right),$$
and 
$$\frac{5x}{3}<\min\left(\frac3y-5, \frac{1-y}{3y}\right).$$
Consequently there exists a rational $x'$ with $x/(1-x)<5x/3<x'$, and 
$$\frac{2x}{3-3y-x}\le x'\le \min\left(\frac3y-5, \frac{1-y}{3y},\frac13\right).$$
Thus 
$$\max\left(\frac{2y-1}{y}, 1-\frac{x'(1-y)}{x}\right)\le \min\left(1-3x', \frac{1-x'}{3}\right);$$
choose a  rational $y'$ between them. Then $x'+3y'\le 1$ and $3x'+y'\le 1$, and so $\psi(x',y')<1/3$, by (theorem 3.3 of the paper).
Let $\psi(x',y')=z'<1/3$, and choose $z<3/5$ with $(1-z)/z\le 1-z'$, and $(1-z)/(1-x)\le 1-z'$, and $z\ge x/x'$.
Then from \ref{psilb2}, $\psi(x,y)\le z<3/5$, a contradiction. This proves \ref{lostcurve35}.~\bbox

\begin{thm}\label{3/5philb}
If $x,y\in (0,1]$ with $y < 1/3$ and $\frac{2x-1}{2-3x}+\frac{y}{1-3y} \le 1$, and with strict inequality if $x$ or $y$ is irrational,
then $\phi(x,y) < 3/5$.
\end{thm}
\Proof
We may assume that $x,y$ are rational. Let $x'=2-1/x$ and $y'=y/(1-y)$; it follows that $x', y' < 1/2$, and
$\frac{x'}{1-2x'}+\frac{y'}{1-2y'} \le 1$. By \ref{phi12curve} for $k = 2$, it follows that
$\phi(x',y') < 1/3$. Then applying \ref{philb} with $z$ slightly less than $3/5$ gives the result.
This proves \ref{3/5philb}.~\bbox

\section{Peaceful coexistence}

We have not been able to evaluate $\phi(x,y)$ in general, but here is an easier question (that we also cannot do, but
it seems to be less far out of reach).
It is always true that $\phi(x,y)\ge x$, by \ref{maxbound}, but if $y$ is sufficiently small then equality may hold.
For fixed $x$, what is the largest $y$ such that $\phi(x,y)=x$?

Let $(G,w)$ be a weighted graph. We say it is {\em $x$-regular} via a bipartition $(A,B)$
if
\begin{itemize}
\item $|A|=|B|$, and $w(v)>0$ for each $v\in V(G)$;
\item the $0,1$-adjacent matrix between $A$ and $B$ is nonsingular;
\item $\sum_{u\in A} w(u) = \sum_{v\in B}w(v)=1$; and
\item for each $u\in V(G)$, $\sum_{v\in N(u)}w(v)=x$.
\end{itemize}
(Note that the fourth bullet is required to hold both for $u\in A$ and for $u\in B$.)
Its {\em order} is $|A|$, and its {\em min-weight} is $\min_{v\in B}w(v)$.
We will show:

\begin{thm}\label{peacephi}
For $x,y\in(0,1]$, $\phi(x,y) = x$ if and only if there is an $x$-regular bipartite weighted graph with order at
most $1/y$.
\end{thm}
\Proof
If there is a such a weighted graph $(G,w)$, via $(A,B)$, where $|A|=|B|=n$ say, let $C$ be a set of $n$ new vertices,
and add a perfect matching between $B$ and $C$. Extend $w$ to $C$ by defining $w(v)=1/n$ for each $v\in C$.
The weighted graph just made is $(x,1/n)$-constrained, and shows that $\phi(x,1/n)\le x$, and consequently $\phi(x,y)\le x$
(and so $\phi(x,y)=x$).

For the converse, suppose that $G$ is $(x,y)$-constrained via $(A,B,C)$, and $|N^2_A(v)|\le x|A|$ for each $v\in C$.
\\
\\
(1) {\em Each vertex in $A$ has exactly $x|B|$ neighbours in $B$, and each vertex in $B$ has exactly
$x|A|$ neighbours in $A$.}
\\
\\
Each vertex $u\in B$ has at most $x|A|$ neighbours in $A$, since $u$ has a neighbour $v\in C$ and
$|N^2_A(v)|\le x|A|$. Since each vertex in $A$ has at least $x|B|$ neighbours in $B$, averaging shows that equality holds
throughout. That proves (1).

\bigskip
Say two vertices in $A$ are {\em twins} if they have the same neighbour set in $B$, and two vertices in $B$ are
{\em twins} if they have the same neighbour set in $A$. This defines equivalence relations of $A$ and $B$, and we call
the equivalence classes {\em twin classes}.
\\
\\
(2) {\em For each vertex $v\in C$, all its neighbours in $B$ are twins, and so $N(v)$ is a subset of a twin class of $B$.}
\\
\\
By (1) each vertex in $N(v)$ has $x|A|$ neighbours in $A$, and all these vertices belong to $N^2_A(v)$;
and since $|N^2_A(v)|=x|A|$, equality holds, and in particular, all vertices in $N(v)$ are twins.
This proves (2).

\bigskip
Let $\mathcal{T}$ be the set of all twins classes of $B$.
For each $T\in \mathcal{T}$, let $C(T)$ be the set of all $v\in C$ with $N(v)\subseteq T$. Thus the sets $C(T)\;(T\in \mathcal{T})$
are nonempty, pairwise disjoint and have union $C$. There is one of cardinality at most $|C|/|\mathcal{T}|$, say
$C(T)$; and then each vertex in $T$ has only at most $|C|/|\mathcal{T}|$ neighbours in $C$, and so $y\le 1/|\mathcal{T}|$.

Choose one vertex from each twin class of $A$ and of $B$, and let $H$ be the subgraph induced on this set.
For each vertex $v$ of $H$, let $w(v)=|T|/|B|$ if $v\in T$ for some twin class $T$ of $B$, and $w(v)=|T|/|A|$
if $v\in T$ for some twin class $T$ of $A$. Then we have:

\begin{itemize}
\item $(H,w)$ is a bipartite graph, with bipartition $(A_0, B_0)$ say;
\item $\sum_{u\in A_0} w(u) = \sum_{v\in B_0}w(v)=1$;
\item for each $u\in V(H)$, $\sum_{v\in N(u)}w(v)=x$; and
\item $|B_0|\le 1/y$.
\end{itemize}

Let us choose a weighted graph $(H,w)$ and bipartition with these properties, with $|V(H)|$ minimum.
If there is a function $f:A\rightarrow \mathbb{R}$ such that $\sum_{u\in N(v)}f(u)=0$ for each $v\in B$, not identically zero,
then by adding a suitable multiple of $f$ to the restriction of $w$ to $A$, we can arrange that $w(u)=0$ for some $u\in A$,
and then $u$ can be deleted, contrary to the minimality of $|V(H)|$. Thus there is no such $f$, and similarly there
is no $f:B\rightarrow \mathbb{R}$ such that $\sum_{v\in N(u)}f(v)=0$ for each $u\in A$, not identically zero. Consequently
$|A_0|=|B_0|=n$ say, and the adjacency matrix between $A_0$ and $B_0$ is nonsingular. Moreover $w(v)>0$
for each $v\in V(H)$, from the minimality of $V(H)$. This proves \ref{peacephi}.~\bbox

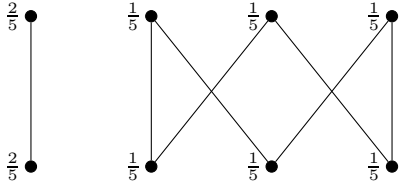
\begin{figure}[H]
\centering

\begin{tikzpicture}[yscale=1,xscale=0.8,auto=left]
\tikzstyle{every node}=[inner sep=1.5pt, fill=black,circle,draw]

\node (a1) at (10,0) {};
\node (a2) at (8,0) {};
\node (a3) at (6,0) {};
\node (a4) at (4,0) {};
\node (b1) at (10,2) {};
\node (b2) at (8,2) {};
\node (b3) at (6,2) {};
\node (b4) at (4,2) {};

\foreach \from/\to in {a1/b2,a1/b1,a2/b1,a2/b3,a3/b2,a3/b3,a4/b4}
\draw [-] (\from) -- (\to);

\tikzstyle{every node}=[left]
\draw (a1) node []           {\scriptsize$\frac15$};
\draw (a2) node []           {\scriptsize$\frac15$};
\draw (a3) node []           {\scriptsize$\frac15$};
\draw (a4) node []           {\scriptsize$\frac25$};
\draw (b1) node []           {\scriptsize$\frac15$};
\draw (b2) node []           {\scriptsize$\frac15$};
\draw (b3) node []           {\scriptsize$\frac15$};
\draw (b4) node []           {\scriptsize$\frac25$};
\end{tikzpicture}

\caption{A $2/5$-regular weighted bipartite graph of order four.} \label{fig:2-5regular}
\end{figure}

By \ref{permute}, $\phi(x,y) = x$ if and only $\phi(y,x)=x$, so this also answers the analogous question for $\phi(y,x)$.
If $x$ is irrational, there is no $x$-regular bipartite weighted graph, and so $\phi(x,y)>x$ for all $y>0$.
If $x\in (0,1]$ is rational,
let us define the {\em order} of $x\in (0,1]$ to be the minimum order of $x$-regular bipartite weighted graphs. If $x=p/q$
say where $p,q>0$ are integers, then the order of $x$ is at most $q$,
because one can construct an appropriate cyclic
shift graph.
But the order of $x$ can be strictly less than $q$. For instance, the top part of the graph of figure
\ref{fig:exactcount} is $13/27$-regular (take as vertex-weights the numbers given, divided by 27), and so the order
of $13/27$
is at most seven.
Figure \ref{fig:2-5regular} gives a smaller example, showing that the order of $2/5$ is at most four.

We can prove that the order is also bounded below by a function of $q$ that goes to infinity with $q$.
More exactly, if $G$ is $p/q$-regular (in lowest terms) and has order $n$, then $q$ is at most $(n+1)^{(n+1)/2}$.
This follows from a theorem of Hadamard~\cite{hadamard}, that
every $n\times n$ $0,1$-matrix has determinant at most $(n+1)^{(n+1)/2}2^{-n}$.
We do not know whether
there are weighted bipartite graphs with order $n$ that are $p/q$-regular (in lowest terms), where $q$ is exponentially large in $n$.
(Hadamard $n\times n$ $0,1$-matrices have determinant that achieve Hadamard's bound, and they exist when $n+1$ is a power of two,
but they give
weighted bipartite graphs that are vertex-transitive, and which therefore are $p/q$-regular with $q=n$.)

One could ask the same question for the biconstrained problem: given $x$, for which values of $y$ is it true that
$\psi(x,y)=x$? A similar analysis (we omit the details) shows:
\begin{thm}\label{peacepsi}
For $x,y\in(0,1]$, the following are equivalent:
\begin{itemize}
\item $\psi(x,y) = x$;
\item $\psi(y,x)=x$; and
\item there is an $x$-regular bipartite weighted graph with min-weight at
least $y$.
\end{itemize}
\end{thm}

\end{document}